\documentclass[a4paper,11pt,onecolumn,twoside]{article}
\usepackage{mathrsfs}
\usepackage{mathrsfs}
\usepackage{amsfonts}
\usepackage{fancyhdr}
\usepackage{amsmath,amsfonts,amssymb,graphicx,amsthm}
\allowdisplaybreaks

\usepackage{amsmath,amsfonts,amssymb,amscd,amsthm,amsbsy,bm,latexsym}
\usepackage{gensymb}
\usepackage[pdftex,bookmarksnumbered=true,bookmarksopen=true,          %bookmark
            colorlinks=true,pdfborder=001,citecolor=blue,
            linkcolor=red,anchorcolor=green,urlcolor=blue]{hyperref}
\usepackage{mathrsfs}
\allowdisplaybreaks

\usepackage{color}

\begin{document}

\title{\bf Global stability of combination of viscous contact wave with rarefaction waves  for the   compressible fluid models of Korteweg type}
\author{
{\bf Zhengzheng Chen}\thanks{Corresponding author. E-mail:
chenzzandu@163.com}\\
School of Mathematical Sciences,
Anhui University, Hefei 230601,  P. R. China\\[2mm]
{\bf Mengdi Sheng}\\
School of Mathematical Sciences,
Anhui University, Hefei 230601,  P. R. China\\[2mm]
}

\date{}

\vskip 0.2cm

\maketitle

\vskip 0.2cm \arraycolsep1.5pt
\newtheorem{Lemma}{Lemma}[section]
\newtheorem{Theorem}{Theorem}[section]
\newtheorem{Definition}{Definition}[section]
\newtheorem{Proposition}{Proposition}[section]
\newtheorem{Remark}{Remark}[section]
\newtheorem{Corollary}{Corollary}[section]

\begin{abstract}

This paper is concerned  with the large-time  behavior of solutions to the Cauchy problem of the one-dimensional compressible fluid models of Korteweg type with density- and  temperature-dependent viscosity, capillarity, and heat conductivity coefficients, which  models the motions of compressible  viscous fluids with internal capillarity. We  show that the combination of the viscous contact wave with two rarefaction waves  is asymptotically stable with a large  initial perturbation if the  strength of the composite wave  and the heat conductivity coefficient satisfy some smallness conditions. The proof   is based on  some   refined $L^2$-energy estimates  to control the  possible growth of the solutions caused by the highly nonlinearity of the system, the  interactions of waves from  different families and large data,  and the key ingredient  is to derive the uniform positive lower and upper bounds on the specific volume and the temperature.

\bigbreak
\noindent

\bigbreak

\end{abstract}

\section{Introduction }
\setcounter{equation}{0}
 This paper is concerned with the Cauchy problem of the  nonisothermal compressible fluid models of Korteweg type in Lagrangian coordinates:
\begin{equation}\label{1.1}
\left\{\begin{array}{l}
  v_t -u_x =0,\\[2mm]
  u_t +p(v,\theta)_x =\displaystyle \left(\frac{\mu(v,\theta)u_x}{v}\right)_x +K_x,\\[2mm]
  \displaystyle\left(C_v-\frac{\theta}{2}\kappa_{\theta\theta}\frac{v_x^2}{v^5}\right)\theta_t +p(v,\theta)u_x =\displaystyle\left(\frac{\tilde{\alpha}(v,\theta)\theta_x}{v}\right)_x +\frac{\mu(v,\theta)u_x^2}{v} +F
\end{array}\right.
\end{equation}
with the  initial  and  far field conditions:
\begin{equation}\label{1.2}
  \left\{\begin{array}{l}
  (v,u,\theta)(0,x)=(v_0,u_0,\theta_0)(x),\\[2mm]
  (v,u,\theta)(t,\pm\infty)=(v_{\pm},u_{\pm},\theta_{\pm}).\end{array}\right.
\end{equation}
Here the unknown functions are the specific volume  $v(x,t)>0$, the velocity $u(x,t)$, the temperature $\theta(x,t)>0$, and the pressure $p(v,\theta)$ of the fluid, respectively, while $\mu(v,\theta),\kappa(v,\theta),\tilde{\alpha}(v,\theta)$  denote the viscosity coefficient, the capillary coefficient and the  heat conductive coefficient respectively.  $C_v>0$, $v_{\pm}>0$, $u_{\pm}$  and $\theta_{\pm}>0$  are given constants, and we assume   $(v_0,u_0,\theta_0)(\pm\infty)=(v_{\pm},u_{\pm},\theta_{\pm})$ as compatibility conditions.  The Korteweg stress tensor  $K$  and the nonlinear terms $F$ are given by
\begin{equation}\label{1.3}
\left\{\begin{array}{l}
  K=\displaystyle\frac{-\kappa(v,\theta)v_{xx}}{v^5}+\frac{5\kappa(v,\theta)-v\kappa_v(v,\theta)}{2v^6}v_x^2-\frac{\kappa_\theta(v,\theta)v_x\theta_x}{v^5},\\[2mm]
  F=\displaystyle\frac{\theta\kappa_{\theta}v_xu_{xx}}{v^5}+\frac{v\kappa_{v\theta}(v,\theta)-\kappa_{\theta}(v,\theta)}{2v^6}\theta u_xv_x^2.
\end{array}\right.
\end{equation}
 Throughout  this paper,  we suppose that the pressure $p(v,\theta)$ and the constant $C_v$ are given by
\begin{equation}\label{1.4}
p(v,\theta)=\frac{R\theta}{v}=Av^{-\gamma}\exp\left(\frac{\gamma-1}{R}s\right),\quad  C_v=\frac{R}{\gamma-1},
\end{equation}
where $s$ is the entropy of the fluid and $\gamma>1, A$ and $R$ are positive constants.

System (\ref{1.1}) can be used to model the  motions of  compressible viscous fluids with internal capillarity.   The formulation of the theory of capillarity with diffuse interface was first studied by  Van der Waals \cite{Van der Waals} and Korteweg  \cite{Korteweg-1901},  and then derived rigorously by Dunn and Serrin \cite{Dunn-Serrin-1985}. Note that if  the capillary coefficient $\kappa=0$, the system (\ref{1.1}) is reduced to the compressible Navier-Stokes system.

There have been extensive studies on the mathematical theory of the compressible fluid models of Korteweg type.  For the case with small initial data,  the results available now are almost complete.  We briefly review some of them here. Hattori and  Li  \cite{H. Hattori-D. Li-1996,Hattori-Li-1996-2} proved the   global existence of   smooth solutions  around constant states in  Sobolev space. Wang and Tan \cite{Wang-Tan-2011} established the optimal time decay  rate of smooth solution obtained in \cite{H. Hattori-D. Li-1996}. Danchin and Desjardins \cite{R. Danchin-B. Desjardins-2001} and Haspot \cite{B. Haspot-2009,B. Haspot-2016} discussed  the global existence and uniqueness of strong  solutions in  Besov space. Kotschote \cite{Kotschote-2014} proved the exponential stability of a non-constant stationary solution  in the phase space.
The authors in \cite{Y.-P. Li-2012,Chen-2012-5,Wang-2015} obtained  the existence and nonlinear stability of non-constant   stationary  solutions  in  Sobolev space.  Chen et al. \cite{Chen-2012-1, Chen-2012-2, Chen-2012-3} studied  the nonlinear stability of some single basic   waves (such as rarefaction wave, viscous shock wave and viscous contact wave) in  Sobolev space.  And the global existence of weak  solutions  in the whole space $\mathbb{R}^2$ was obtained by Danchin and Desjardins \cite{R. Danchin-B. Desjardins-2001} and Haspot\cite{B. Haspot-2011}.

For the case with large initial data, Haspot \cite{B. Haspot-2017} proved the global  existence of strong solution for an isothermal fluid with  density-dependent  viscosity and capillary coefficients in the whole space $R^N$ with $N\geq2$.
Bresch, Desjardins, and Lin  \cite{D. Bresch-B. Desjardins-C. K. Lin-2003} studied  the global existence of weak solutions for an  isothermal Korteweg system  with  a linearly density-dependent viscosity    and a constant capillarity coefficient in  a periodic domain $\mathbb{T}^d$ with $d=2$ or $3$. Then  Haspot  \cite{B. Haspot-2011} and J\"{u}ngel  \cite{Jungel-2010} improved the results of \cite{D. Bresch-B. Desjardins-C. K. Lin-2003} to some other types of  density-dependent viscosity  and capillarity coefficients.  Tsyganov \cite{E-Tsyganov-2008} discussed  the global existence of weak solutions for an isothermal system  with the viscosity coefficient $\mu(\rho)\equiv1$  and the capillarity coefficient $\kappa(\rho)={\rho^{-5}}$ on the interval $[0,1]$.   Germain and LeFloch \cite{Germain-LeFloch-2012} investigated   the global existence of weak solutions  for the isothermal Korteweg system with general density-dependent viscosity  and  capillarity coefficients in $\mathbb{R}$.    The global existence of large strong solution to  an isothermal Korteweg system with the viscosity coefficient  $\mu(\rho)=\varepsilon\rho$ and the capillarity coefficient $\kappa(\rho)=\varepsilon^2{\rho^{-1}}$ in $\mathbb{R}$ ($\varepsilon$ is positive constant) was obtained by Charve and  Haspot \cite{F-Charve-B-Haspot-2011}. For  the global existence of  smooth large-amplitude solutions to the compressible  Korteweg system in the whole space $\mathbb{R}$, we refer to \cite{Chen-2012-5,Chen-2016} for  isothermal system with  general density-dependent viscosity and  capillarity coefficients, \cite{Chen-2017} for nonisothermal system with general density-dependent viscosity and  capillarity coefficients, and density- and temperature-dependent heat conductivity coefficient.
 The time-asymptotic nonlinear stability of strong rarefaction waves for the isothermal Korteweg system with large initial data  was also obtained in \cite{Chen-2012-5}.

From the above result,  it is seen that few results have been obtained on the global stability of basic  waves for the compressible fluid models of Korteweg system so far.  Here and hereafter,  {\it global stability} means the nonlinear stability result with large initial perturbation. And if the initial perturbation is small, the  nonlinear stability result is usually called {\it local stability}.  A  natural question is whether some global stability results  of the composite waves for the one-dimensional  compressible  Korteweg system hold or not? This paper is devoted to this problem and we are concerned with the global stability of the combination   of viscous contact waves with  two rarefaction waves for the Cauchy problem of the nonisothermal Korteweg system (\ref{1.1}) with general density- and  temperature-dependent viscosity, capillarity, and heat conductivity coefficients.

It is known \cite{Chen-2012-1,Chen-2012-2,Chen-2012-3} that the  large-time behavior of solutions of the Cauchy problem (\ref{1.1})-(\ref{1.2}) is closely related to the Riemann problem of the compressible  Euler system£º
\begin{equation}\label{1.5}
\left\{\begin{array}{l}
    v_t-u_x=0, \\[2mm]
    u_t+p(v,\theta)_x=0, \\[2mm]
    C_v\theta_t+p(v,\theta)u_x=0
  \end{array}\right.
\end{equation}
with the  Riemann  initial data
\begin{equation}\label{1.6}
(v,u,\theta)(0,x)=\left\{\begin{array}{c}
                             (v_-,u_-,\theta_-),\quad   x<0, \\[2mm]
                             (v_+,u_+,\theta_+),\quad   x>0.
                           \end{array}\right.
\end{equation}
The Euler system (\ref{1.5}) is a strict hyperbolic system of conservation laws with three distinct eigenvalues:
\[\lambda_1(v,\theta)=-\sqrt{\frac{\gamma p}{v}},\quad\lambda_2=0,\quad \lambda_3(v,\theta)=\sqrt{\frac{\gamma p}{v}},\]
which implies that the second characteristic field is linearly degenerate and the others are genuinely nonlinear.
Then it is well-known \cite{J. Smoller-1994} that the Riemann problem (\ref{1.5})-(\ref{1.6}) admits three basic wave patterns: the shock wave, rarefaction wave and contact discontinuity, and the Riemann solution to (\ref{1.5}) has a wave pattern consisting  of a linear combination of these three basic waves.
In particular,  the contact discontinuity solution of the Riemann problem (\ref{1.5})-(\ref{1.6}) takes the form \cite{J. Smoller-1994}
\begin{equation}\label{1.7}
(v^{cd},u^{cd},\theta^{cd})(t,x)=\left\{\begin{array}{c}
                                          (v_-,u_-,\theta_-),\quad   x<0,\ t>0, \\[2mm]
                                          (v_+,u_+,\theta_+),\quad   x>0,\ t<0,
                                        \end{array}\right.
\end{equation}
provided that
\begin{equation}\label{1.71}
u_-=u_+,\qquad p_-\triangleq\frac{R\theta_-}{v_-}=\frac{R\theta_+}{v_+}\triangleq p_+.
\end{equation}

The viscous contact wave $(V,U,\Theta)(t,x)$ corresponding to  the contact discontinuity  $(v^{cd},u^{cd},\theta^{cd})(t,x)$ for the compressible Navier-Stokes-Korteweg system (\ref{1.1}) becomes smooth and behaviors as a diffusion waves due to the  effect of heat conductivity. As \cite{Chen-2012-2}, we can define the  viscous contact wave $(V, U, \Theta)(t,x)$ as follows.

Since the properties of the contact discontinuity wave motivate us to expect that
\begin{equation}\label{1.8}
  P=\frac{R\Theta}{V}\approx p_+=p_-, \quad |U|\ll1,
\end{equation}
the leading part of the energy equation  $(\ref{1.1})_3$ is
\begin{equation}\label{1.9}
  C_v\Theta_t+p_+U_x=\left(\frac{\widetilde{\alpha}(V,\Theta)\Theta_x}{V}\right)_x.
\end{equation}
Using the equations (\ref{1.8}), $V_t=U_x$ and (\ref{1.9}), we get  a nonlinear diffusion equation
\begin{equation}\label{1.10}
\left\{\begin{array}{l}
         \Theta_t=a\displaystyle\left(\frac{\hat{\alpha}(\Theta)\Theta_x}{\Theta}\right)_x,\quad a=\displaystyle\frac{p_+(\gamma-1)}{\gamma R^2}, \\[3mm]
      \Theta(\pm\infty,t)=\theta_{\pm},
       \end{array}\right.
\end{equation}
where $\hat{\alpha}(\Theta)\triangleq\widetilde{\alpha}(\frac{R\Theta}{p_+},\Theta)$.
Due to \cite{C.J. van Duyn-L.A. Peletier-1996/1997}, (\ref{1.10})  has a unique self-similar  solution $\Theta(t,x)=\Theta(\xi),\xi=\frac{x}{\sqrt{1+t}}$, which is a monotone function, increasing if $\theta_+>\theta_-$ and decreasing
if $\theta_+<\theta_-$.

 Once $\Theta(t,x)$ is determined, the viscous contact wave $(V, U, \Theta)(t,x)$ is defined by
\begin{equation}\label{1.11}
  V=\frac{R{\Theta}}{p_+},\quad U=u_-+\frac{\gamma-1}{\gamma R}\frac{\hat{\alpha}(\Theta){\Theta_x}}{\Theta},\quad {\Theta}=\Theta(t,x),
\end{equation}
then it is easy to check  that  the viscous contact wave $(V, U, \Theta)(t, x)$ satisfies
\begin{eqnarray}\label{1.13}
  \left\{\aligned
          &V_t- U_x=0, \\[2mm]
          &U_t+\displaystyle\left(\frac{R\Theta}{V}\right)_x=\displaystyle\left(\frac{\mu(V,\Theta)U_x}{V}\right)_x+R_1, \\[2mm]
          &\displaystyle\frac{R}{\gamma-1}{\Theta}_t+p(V, \Theta)U_x=\displaystyle\left(\frac{\tilde{\alpha}(V, \Theta)}{V}{\Theta}_x\right)_x+\frac{\mu(V, \Theta)U_x^2}{V}+R_2
  \endaligned\right.
\end{eqnarray}
with
\begin{equation}\label{1.13-1}
\aligned
    R_1 &= U_t-\left(\frac{\mu(V,\Theta)U_x}{V}\right)_x=\displaystyle O(1)(\gamma-1)^{\frac{1}{2}}(1+t)^{-\frac{3}{2}}e^{-\frac{c_0x^2}{(\gamma-1)(1+t)}}, \\
    R_2 &= \displaystyle-\mu(V,\Theta)\frac{U_x^2}{V}=O(1)(\gamma-1)^2(1+t)^{-2}e^ {-\frac{c_0x^2}{(\gamma-1)(1+t)}}.
    \endaligned
\end{equation}

Our first theorem is concerned with  the global stability of the single viscous contact wave $(V,U,\Theta)(t,x)$,  which is stated   as follows.
\begin{Theorem}[Global stability of viscous contact wave]  Let $(V, U, \Theta)(t,x)$ be the viscous contact wave defined in (\ref{1.11}). Suppose that
\begin{itemize}
\item[(i)] The given constants $v_{\pm}$, $u_{\pm}$, $\theta_{\pm}$ do not depend on $\gamma-1$, and satisfy (\ref{1.71}). Moreover, $|\theta_+-\theta_-|\leq m_0(\gamma-1)$ for some positive constant $m_0$ independent of $\gamma-1$;
\item[(ii)]The initial data $\left(v_0(x)-V(0,x), u_0(x)-U(0,x), \theta_0(x)-\Theta(0,x)\right)\in H^2(\mathbb{R})\times H^1(\mathbb{R})\times H^1(\mathbb{R})$ and $N_0:=\|v_0(x)-V(0,x)\|_{H^2(\mathbb{R})}+\|u_0(x)-U(0,x)\|_{H^1(\mathbb{R})}+\left\| (\theta_0(x)-\Theta(0,x))/\sqrt{\gamma-1}\right\|_{H^1(\mathbb{R})}$ is bounded by some constant independent of $\gamma-1$;
\item[(iii)]
There exist positive constants $\underline{V},$ $\overline{V},$  $\underline{\Theta}$ and  $\overline{\Theta}$ independent of $\gamma-1$ such that
\begin{equation}\label{1.14}
\underline{V}\leq v_0(x), V(t,x)\leq\overline{V},\quad \underline{\Theta}\leq \theta_0(x), \Theta(t,x)\leq\overline{\Theta},\quad \forall\, (t,x)\in\mathbb{R}^+\times\mathbb{R};
\end{equation}
 \item[(iv)] The viscosity coefficient $\mu(v,\theta)$, the capillarity coefficient $\kappa(v,\theta)$  and the heat-conductivity coefficient  $\tilde{\alpha}(v,\theta)$  are  smooth positive functions of $v>0$ and $\theta>0$, and  the  following assumptions hold:\\
     (a) \begin{equation}\label{1.16-1}
\left\|\frac{\mu(v,\theta)}{\tilde{\alpha}(v,\theta)}\right\|_{L^\infty([0,\infty)\times\mathbb{R})}\leq M_0,\quad|\kappa_{\theta}(v,\theta)|<\varepsilon,\quad \kappa_{\theta\theta}(v,\theta)<0,
\end{equation}
where $M_0>0$ is a uniform constant independent of $\gamma-1$ and $\varepsilon>0$ is small positive constant whose precise range can be specified in the proof of Theorem 1.1.

   (b) One of the following  two conditions holds.

    \quad$(b)_1$ There exist constants $a\geq 0$ and  $b\le \frac{1}{2}$ such that
     \begin{equation}\label{1.15}
  \quad \mu_1(v)=\min_{\theta\in[\frac{\underline \Theta}{2},2\overline \Theta]}\{\mu(v,\theta)\}\sim\left\{
  \begin{array}{l}
   v^{-a},\quad v\longrightarrow 0^+, \\[2mm]
    v^{-b},\quad v\longrightarrow +\infty.
    \end{array}
    \right.
\end{equation}

\quad$(b)_2$ There exist constants $c\leq3$ and $d\geq2$ such that
\begin{equation}\label{1.16}
  \quad \kappa_1(v)=\min_{\theta\in[\frac{\underline \Theta}{2},2\overline \Theta]}\{\kappa(v,\theta)\}\sim\left\{
  \begin{array}{l}
     v^{-c},\quad v\longrightarrow 0^+, \\[2mm]
     v^{-d},\quad v\longrightarrow +\infty.
      \end{array}
      \right.
\end{equation}

(c) The viscosity and capillarity coefficients are coupled by
\begin{equation}\label{1.16-3}
\aligned
f(v,\theta):&=\displaystyle-\frac{2}{3}\left(\frac{\sqrt{\mu\kappa}}{v^3}
  \left(\frac{\sqrt{\mu\kappa}}{v^3}\right)_v\right)_v+\left[\left(\frac{\sqrt{\mu\kappa}}{v^3}\right)_v\right]^2+\displaystyle\frac{1}{3}\left(\frac{\kappa}{v^5}\left(\frac{\mu}{v}\right)_v\right)_v\\
  &\quad+\displaystyle\left(\frac{\mu}{v}\right)_v\frac{5\kappa-v\kappa_v}{2v^6}-\frac{1}{3}\left(\frac{\mu}{2v^7}(5\kappa-v\kappa_v)\right)_v\\
  &\leq0,
\endaligned
\end{equation}
\quad \quad or\begin{equation}\label{1.16-4}
\aligned
g(v,\theta):=3\kappa \mu+2v\kappa\mu_v-v\mu\kappa_v=0.
\endaligned
\end{equation}
\end{itemize}
Then there exist  positive constants $\varepsilon_0\ll1, \delta_0\ll1$ and  $C_0$ which depend only on $\underline{V},$ $\overline{V},$  $\underline{\Theta}$,  $\overline{\Theta}$ and the initial data $N_0$  such that if $0<\varepsilon<\varepsilon_0$ and $0<\delta:=\gamma-1\leq\delta_0$, the Cauchy problem (\ref{1.1})-(\ref{1.2}) admits a unique global-in-time solution $(v,u,\theta)(t,x)$ satisfying
\begin{equation}\label{1.17-1}
\left\{ \begin{array}{l} \aligned
&(v-V,u-U, \theta-\Theta)(t,x)\in C(0,T;H^2(\mathbb{R})\times H^1(\mathbb{R})\times H^1(\mathbb{R})),\\[1mm]
   &(v-V,u-U, \theta-\Theta)_x(t,x)\in L^2(0,T;H^2(\mathbb{R})\times H^1(\mathbb{R})\times H^1(\mathbb{R})),
   \endaligned\end{array}\right.
\end{equation}
\begin{equation}\label{1.17}
C_0^{-1}\leq v(t,x)\leq C_0,\quad \frac{\underline{\Theta}}{2}\leq \theta(t,x)\leq 2\overline{\Theta}, \quad \forall\,(t,x)\in[0,\infty)\times\mathbb{R},
\end{equation}
and  \begin{equation}\label{1.18}
\lim_{t\rightarrow+\infty}\sup_{x\in\mathbb{R}}\left\{\left|\left(v-V,u-U, \theta-\Theta\right)(t,x)\right|\right\}=0.
\end{equation}
\end{Theorem}

\begin{Remark}
Several remarks on  Theorem 1.1  are given  as follows.
\begin{itemize}

\item[$(1)$] Since $\theta_{\pm}=\frac{A}{R}v_{\pm}^{1-\gamma}\exp\left(\frac{\gamma-1}{R}s_{\pm}\right)$ and the constants  $\theta_{\pm},v_{\pm}, s_{\pm}$ are assumed independent of $\gamma-1$,  it is  easy to check  that $|\theta_+-\theta_-|\leq m_0(\gamma-1)$ for  some
constant $m_0$ independent of $\gamma-1$.

\item[$(2)$] The  assumption  (iv)(a) is used to control the possible growth of solutions to the Cauchy problem (\ref{1.1})- (\ref{1.2}) induced by the nonlinearity of the system, (iv)(b) is  used to derive the uniform-in-time lower and upper bounds for the specific  volume $v(t,x)$, and (iv)(c) is a technical condition in estimating  $\|\frac{\mu(v,\theta)\phi_x}{v}(\tau)\|$ (see the proof of Lemmas 3.3-3.5 for details).

\item[$(3)$] In Theorem 1.1,  although the initial perturbation $\|\theta_0(x)-\Theta(0,x)\|_{H^1(\mathbb{R})}$ is small when $\gamma>1$ is close to $1$,  the initial perturbations $\|v_0(x)-V(0,x)\|_{H^2(\mathbb{R})}$ and $\|u_0(x)-U(0,x)\|_{H^1(\mathbb{R})}$ can be arbitrarily  large.  This improves the main result of \cite{Chen-2012-2}, where the nonlinear stability of viscous contact wave for the one-dimensional compressible fluid models of Korteweg type was obtained with all the initial perturbations are sufficiently small.

\item[$(4)$]From the proof of Theorem 1.1, wee see that $\gamma-1$ needs to be sufficiently small such that $(\gamma-1)F(N_0)<1$ with $F(N_0)$ being a smooth increasing function on the initial data $N_0$ (see (\ref{3.1111})-(\ref{3.1112})). Thus in this sense,  Theorem 1.1 is a Nishida-Smoller type result \cite{Nishida-Smoller-1973} with large initial data.
\end{itemize}
\end{Remark}

When the relation (\ref{1.71}) fails,   the basic theory of hyperbolic systems of conservation laws \cite{J. Smoller-1994} tells us  that for any given constant state $(v_-, u_-,\theta_-)$ with $v_->0$, $u_-\in \mathbb{R}$ and $\theta_->0$,  there exists a neighborhood $\Omega(v_-, u_-,\theta_-)$ of $(v_-, u_-,\theta_-)$ such that for any $(v_+, u_+,\theta_+)\in \Omega(v_-, u_-,\theta_-)$, the Riemannn problem (\ref{1.5})-(\ref{1.6}) has a unique solution. In this paper, we only consider the stability of the combination of the viscous contact wave and rarefaction waves. Consequently, we assume that
\begin{equation}\label{1.19-1}
(v_+, u_+,\theta_+)\in R_-C_cR_+(v_-, u_-,\theta_-)\subset\Omega(v_-, u_-,\theta_-), \quad |\theta_+-\theta_-|\leq\delta_1,
\end{equation}
where $\delta_1$ is a positive constant, $R_-,R_+, C_c$ denote  the 1-rarefaction wave curve, 3-rarefaction wave curve, and  the contact wave curve respectively, and
\[
\begin{array}{l} R_-C_cR_+(v_-, u_-,\theta_-)\triangleq\Bigg\{(v, u,\theta)\in\Omega(v_-, u_-,\theta_-)\Bigg|
s\neq s_-,
\\[2mm]
\displaystyle u\geq u_--\int_{v_-}^{e^{\frac{\gamma-1}{R\gamma}(s_--s)}v}\lambda_-(\eta, s_-)\,d\eta,\quad
   u\geq u_--\int_{e^{\frac{\gamma-1}{R\gamma}(s-s_-)}v_-}^{v}\lambda_+(\eta, s)\,d\eta\Bigg\}\end{array}
\]

 Due to  \cite{J. Smoller-1994}, if $\delta_1$ in (\ref{1.19-1}) is suitably small,  then there exist a positive constant $C=C(\theta_-,\delta_1)$ and a pair of points $(v_-^m, u^m,\theta_-^m)$ and $(v_+^m, u^m,\theta_+^m)$ in $\Omega(v_-, u_-,\theta_-)$ such that
 \begin{equation}\label{1.19}
\frac{R\theta_-^m}{v_-^m}=\frac{R\theta_+^m}{v_+^m}\triangleq p^m, \quad |v_{\pm}^m-v_\pm|+|u^m-u_{\pm}|+|\theta_{\pm}^m-\theta_\pm|\leq C|\theta_--\theta_+|.
\end{equation}
Moreover, the states $(v_-^m, u^m,\theta_-^m)$ and $(v_+^m, u^m,\theta_+^m)$ belong to the  1-rarefaction wave curve $R_-(v_-, u_-,\theta_-)$ and  3-rarefaction wave curve $R_+(v_+, u_+,\theta_+)$  respectively, where
\begin{equation}\label{1.20}
R_\pm(v_\pm, u_\pm,\theta_\pm)=\left\{s=s_\pm, u=u_\pm-\int_{v_\pm}^v\lambda_\pm(\eta,s_\pm)\,d\eta,\,\, v>v_\pm\right\}
\end{equation}
with \[ \aligned
&s=\frac{R}{\gamma-1}\ln\frac{R\theta}{A}+R\ln v, \quad s_\pm=\frac{R}{\gamma-1}\ln\frac{R\theta_\pm}{A}+R\ln v_\pm,\\[2mm]
&\lambda_\pm(v,s)=\pm\sqrt{A\gamma v^{-\gamma-1}e^{(\gamma-1)s/R}}.
 \endaligned\]
The contact discontinuity wave curve $C_c$  is defined by
\begin{equation}\label{1.21}
{C_c}(v_-^m, u^m,\theta_-^m)=\left\{(v,u,\theta)(t,x)|u=u^m,\,
p=p^m,\,v\not\equiv v_-^m\right\}.
\end{equation}

The 1-rarefaction wave $(v_-^r, u^r_-,\theta_-^r)(\frac{x}{t})$ (respectively the  3-rarefaction wave $(v_+^r, u^r_+,\theta_+^r)(\frac{x}{t})$) connecting $(v_-, u_-,\theta_-)$ and $(v_-^m, u^m,\theta_-^m)$ (respectively $(v_+^m, u^m,\theta_+^m)$ and $(v_+, u_+,\theta_+)$) is the weak solution of the Riemann problem of the  Euler system (\ref{1.5}) with the Riemann initial data:

\begin{equation}\label{1.22}
(v_\pm^r, u^r_\pm,\theta_\pm^r)(0,x)=\left\{\begin{array}{c}
                             (v_\pm^m,u^m,\theta_\pm^m),\quad   \pm x<0, \\[2mm]
                             (v_\pm,u_\pm,\theta_\pm),\quad   \pm x>0.
                           \end{array}\right.
\end{equation}
The contact discontinuity  wave $(v^{cd},u^{cd},\theta^{cd})(t,x)$ connecting $(v_-^m, u^m,\theta_-^m)$ and $(v_+^m, u^m,\theta_+^m)$  takes  the form
\begin{equation}\label{1.23}
(v^{cd},u^{cd},\theta^{cd})(t,x)=\left\{\begin{array}{c}
                                          (v_-^m, u^m,\theta_-^m),\quad   x<0,\ t>0, \\[2mm]
                                          (v_+^m, u^m,\theta_+^m),\quad   x>0,\ t<0.
                                        \end{array}\right.
\end{equation}

Since the rarefaction waves $(v_\pm^r, u^r_\pm,\theta_\pm^r)(t,x)$ are not smooth enough, to study the stability problem, we need to construct their smooth approximations. As \cite{A. Matsumura-K. Nishihara-1986}, the smooth approximate rarefaction waves $(V_\pm^r, U^r_\pm,\Theta_\pm^r)(t,x)$ of $(v_\pm^r, u^r_\pm,\theta_\pm^r)(t,x)$ can be defined by
\begin{eqnarray}\label{1.25}
\left\{\begin{array}{ll}
        \lambda_{\pm}(V_\pm^r,s_\pm)=w_{\pm}(x,t),\\[2mm]
        U^{r}_\pm(t,x)=\displaystyle u_\pm-\int_{v_\pm}^{V^{r}(t,x)}\lambda_\pm(\eta, s_\pm)\,d\eta,\\[4mm]
        \Theta_\pm^r=\theta_\pm(v_\pm)^{\gamma-1}(V_\pm^r)^{1-\gamma},
 \end{array}\right.
\end{eqnarray}
where $w_-$ (respectively $w_+$) is the solution of the Cauchy  problem of the Burger equation:
\begin{eqnarray}\label{1.26}
\left\{\begin{array}{ll}
         w_{t}+ ww_{x}=0,\,\, x\in\mathbb{R}, t>0,\\[2mm]
 w(0,x)=\displaystyle\frac{w_++w_-}{2}+\displaystyle\frac{w_+-w_-}{2} \tanh x
 \end{array}\right.
\end{eqnarray}
with $w_-=\lambda_-(v_-,s_-)$ and $w_+=\lambda_-(v_-^m,s_-)$ (respectively $w_-=\lambda_+(v_+^m,s_+)$ and $w_r=\lambda_+(v_+,s_+)$).

Let $(V^{c},  U^{c},  \Theta^{c})(t,x)$ be  the viscous contact wave defined in (\ref{1.11}) with $(v_\pm, u_\pm,\theta_\pm)$ replaced by $(v^m_\pm, u^m_\pm,\theta_\pm^m)(t,x)$ respectively. Set
\begin{eqnarray}\label{1.27}
\left(\begin{array}{ll}
        V\\[1mm]
         U\\[1mm]
     \Theta\\[1mm]
 \end{array}\right)(t,x)=\left(\begin{array}{c}
        V^{r}_-+V^c+V^{r}_+\\[1mm]
        U^{r}_-+U^c+U^{r}_+\\[1mm]
        \Theta^{r}_-+\Theta^c+\Theta^{r}_+\\[1mm]
 \end{array}\right)(t,x)-\left(\begin{array}{c}
        v_{-}^m+v_+^m\\[1mm]
       2u^m\\[1mm]
        \theta_{-}^m+\theta_{+}^m\\[1mm]
 \end{array}\right),
 \end{eqnarray}
then our second main  result  is as follows.

\begin{Theorem}[Global stability of composite waves]
Suppose that the  constant states $(v_{\pm}, u_{\pm}, \theta_{\pm})$ satisfy (\ref{1.19-1}) for some  small constant $\delta_1>0$, and  $|\theta_+-\theta_-|\leq m_0(\gamma-1)$ for some $(\gamma-1)-$independent positive constant $m_0$.  Let $(V, U, \Theta)(t,x)$ be the  combination of the viscous contact wave and approximate rarefaction waves defined in (\ref{1.27}), and the conditions (ii)-(iv) of Theorem 1.1 hold. Then there exist  positive constants $\varepsilon_1\ll1, \delta_2\ll1$ and  $C_1$ which  depend only on $\underline{V},$ $\overline{V},$  $\underline{\Theta}$,  $\overline{\Theta}$ and the initial data $N_0$  such that if $0<\varepsilon\leq\varepsilon_1$ and  $0<\delta:=\gamma-1\leq\delta_2$, the Cauchy problem (\ref{1.1})-(\ref{1.2}) admits a unique global-in-time  solution $(v,u,\theta)(t,x)$ satisfying
\begin{equation}\label{1.28-1}
\left\{ \begin{array}{l} \aligned
&(v-V,u-U, \theta-\Theta)(t,x)\in C(0,T;H^2(\mathbb{R})\times H^1(\mathbb{R})\times H^1(\mathbb{R})),\\
   &(v-V,u-U, \theta-\Theta)_x(t,x)\in L^2(0,T;H^2(\mathbb{R})\times H^1(\mathbb{R})\times H^1(\mathbb{R})),
   \endaligned\end{array}\right.
\end{equation}
\begin{equation}\label{1.28}
C_1^{-1}\leq v(t,x)\leq C_1,\quad \frac{\underline{\Theta}}{2}\leq \theta(t,x)\leq 2\overline{\Theta}, \quad \forall\,(t,x)\in[0,\infty)\times\mathbb{R},
\end{equation}
and   \begin{equation}\label{1.29}
\lim_{t\rightarrow+\infty}\sup_{x\in\mathbb{R}}\left\{\left|\left(v-V,u-U, \theta-\Theta\right)(t,x)\right|\right\}=0.
\end{equation}
\end{Theorem}

\begin{Remark} Two remarks on  Theorems 1.1-1.2  are listed below.
\begin{itemize}
\item[$(1)$] From Lemma 2.3 and  (\ref{1.29}), we  have also the following asymptotic behavior of solutions:
\[\label{1.30}
\lim_{t\rightarrow+\infty}\sup_{x\in\mathbb{R}}\left(\begin{array}{c}
       \left|\left(v-v_-^r-V^{c}-v_+^r+v_-^m+v_+^m\right)(t,x)\right|\\[2mm]
       \left|\left(u-u_-^r-u_+^r+u^m\right)(t,x)\right|\\[2mm]
        \left|\left(\theta-\theta_-^r-\Theta^{c}-\theta_+^r+\theta_-^m+\theta_+^m\right)(t,x)\right|\end{array}\right)=0,
\]
where $(v_-^r, u_-^r, \theta_-^r)(t,x)$ and $(v_-^r, u_-^r, \theta_-^r)(t,x)$ are the 1-rarefaction wave and 3-rarefaction wave uniquely determined by (\ref{1.5}), (\ref{1.22}), respectively.

\item[$(2)$]In Theorems 1.1-1.2, the smallness of  $\gamma-1$ plays an important role in our analysis.  Recently,  Huang and Wang \cite{Huang6} studied the  global stability of   the combination  of  viscous contact wave with  rarefaction waves for the Cauchy problem of the 1-D compressible Navier-Stokes system without the conditions that $\gamma$ is  close to $1$. However, it seems that the method of \cite{Huang6}  can not be applied to the nonisothermal compressible Navier-Stokes-Korteweg system (\ref{1.1})  because of some difficult nonlinear terms caused by the Korteweg tensor. The problem on how to get the  global stability of basic waves for the nonisothermal compressible fluid models of Korteweg type  with general constant $\gamma>1$  is  under our current research.
\end{itemize}
\end{Remark}

Now we outline the main ideas used in proving Theorems 1.1-1.2. The key ingredient in the proof of Theorem 1.1 is to deducing the uniform-in-time  positive lower and upper bounds on the specific volume $v(t,x)$ and the  temperature $\theta(t,x)$. To achieve this, we make the a priori assumption $\|(\theta-\Theta)(t)/\sqrt{\gamma-1}\|_1\leq N_1,\,\,\forall\, t\in[0,T]$ for some positive constants $N_1$ and $T$. Then the lower and upper bounds for the  temperature $\theta(t,x)$ follow easily by
 the smallness assumption of  $\gamma-1$ and the  Sobolev inequality (see (\ref{3.11})-(\ref{3.12}) for details). The bounds for   the specific volume $v(t,x)$  from below and above were established  by using Kanel's technique \cite{Y. Kanel}, which is based on the basic energy estimates of solutions $(\phi,\psi,\zeta)$  to the reformulated system (\ref{3.1})  and the estimate of   $\|\frac{\mu(v,\theta)\phi_x}{v}(\tau)\|$ (see Lemma 3.4 and Corollary 3.1). Here we remark that even in the case of constant viscosity coefficient, the classical method of Kazhikhov and Shelukhin \cite{Kazhikhov-Shelukhin-1977} can't yields the desired lower and upper bounds on $v(t,x)$ and $\theta(t,x)$ for the compressible Navier-Stokes-Korteweg system (\ref{1.1})  due to the appearance of the Korteweg tensor. Once the uniform lower and upper bounds on $v(t,x)$ and $\theta(t,x)$ are obtained, the higher order energy estimates of  solutions  can be deduced by using the lower order energy estimates and Gronwall's inequality.  Compared with former results \cite{E-Tsyganov-2008,Germain-LeFloch-2012,F-Charve-B-Haspot-2011,Chen-2012-5,Chen-2017} on the construction of global large solutions to the one-dimensional compressible fluid models of Korteweg type,  there are two additional  difficulties in  our analysis.
\begin{itemize}
\item[$1.$] The first  one lies in dealing with  the  highly  nonlinear  terms in  system (\ref{1.1}) under large initial perturbation, such as the terms $K$ and $F$ in (\ref{1.3}).
The system (\ref{1.1}) under consideration here is more complex than those in
\cite{E-Tsyganov-2008,Germain-LeFloch-2012,F-Charve-B-Haspot-2011,Chen-2012-5,Chen-2017}. In particular,  the viscosity coefficient  $\mu(v,\theta)$ and  the capillarity  coefficient $\kappa(v,\theta)$ here   can  depend on both the specific volume  $v$ and the temperature $\theta$. As far as we know,  temperature-dependent viscosity and capillarity  coefficients are not considered for the global large solutions to the compressible fluid models of Korteweg type in the literatures available now.
These density- and  temperature-dependent physical coefficients will enhance  the nonlinearity of the system and thus leads to difficulty in analysis for global solvability with large data. To control the possible  growth of solutions caused by the highly nonlinearity of the system,  we mainly use the smallness assumptions on $\gamma-1$ and  $\kappa_\theta(v,\theta)$ (see (\ref{1.16-1})),  and some  elaborate  analysis.
\item[$2.$] The second difficulty is to control the growth of solutions caused  by the viscous contact wave, which is quite different from the case of the compressible Navier-Stokes system \cite{Hakho-Hong-2010,Huang-Liao-2017}. Due to the effect of the Korteweg tensor, a third order spatial derivative of the viscous contact wave $\Theta_{xxx}(t,x)$ appears in the estimate of $\int_0^t\|\phi_{xxx}(\tau)\|^2d\tau$. Consequently, $\int_0^t\|\phi_{xxx}(\tau)\|^2d\tau$ may be  bounded by $C(\gamma-1)^{-\frac{1}{2}}$ for  some positive constant $C$ independent of $t$ and $\gamma-1$ (see (\ref{3.108})). Moreover,  $\int_0^t\|\phi_{xxx}(\tau)\|^2d\tau$ also presents as a remainder term in the estimate of $\|\psi_x(t)\|$ (see (\ref{3.1041})). Since $\gamma-1$ should be small in our setting, the terms $\int_0^t\|\phi_{xxx}(\tau)\|^2d\tau$ and $\|\psi_x(t)\|$ will  grow as $\gamma-1$ decrease. Thus a difficulty problem is how  to control the growth of solutions induced by $\Theta_{xxx}(t,x)$.  To over come such a difficulty,  we make the a priori  assumption (\ref{3.10})-(\ref{3.10-1}), which together with  a careful continuation argument can  yield  the desired energy-type  estimates of solutions to the Cauchy problem (\ref{3.1})-(\ref{3.2}).
\end{itemize}

The proof of Theorem 1.2 is similar to that of Theorem 1.1, but with an  additional difficulty to control the  interactions of wave from different families. With the aid of a domain  decomposition technique developed in \cite{F.M. Huang-J. Li-A. Matsumura-2010} and the properties of the approximate rarefaction waves and  viscous contact wave, we can successfully  overcome this difficulty and finally get the desired a priori estimates for  solutions of the Cauchy problem (\ref{4.5})-(\ref{4.6}).

Before concluding this section, we should  mention that the nonlinear stability of basic waves  for the compressible Navier-Stokes equations has been studied by many authors. We  refer to \cite{A. Matsumura-K. Nishihara-1985,Matsumura-Mei-1999, Liu-1986,F.-M. Huang-Matsumura-2009} for the nonlinear  stability of viscous shock waves, \cite{A. Matsumura-K. Nishihara-1986,Liu-Xin-1988,K. Nishihara-T.Yong-H. J. Zhao-2004}  for the nonlinear stability of rarefaction waves, \cite{F.M. Huang-J. Li-A. Matsumura-2010,Huang2,Huang3,Huang4,Huang5,Huang6,Hakho-Hong-2010}  for the nonlinear stability of contact discontinuity, and  \cite{F.M. Huang-J. Li-A. Matsumura-2010,Huang6,Hakho-Hong-2017,Huang-Liao-2017,Fan-2016,Wan-2016} for the nonlinear stability of composite waves.

The rest of this paper is organized as follows.  In Section 2, we list some basic properties of the viscous contact wave and rarefaction waves  for later use. An important lemma concerning  the heat kernel and a domain decomposition technique were also presented in this section.  The main theorems 1.1 and 1.2 will be proved in Sections 3 and 4, respectively.

{\bf Notations:} Throughout this paper, we denote $\delta:=\gamma-1$ for notational simplicity.  $c, C$ and $O(1)$ stand for some generic positive
constants which may depend on $\underline V,\overline V,\underline \Theta,\overline \Theta, m_0$ and $M_0$ ($m_0,M_0$ are positive constants given in Proposition 3.2 ),  but are independent of $t$ and $\delta$. If the dependence need to be explicitly pointed out,
the natation  $C(\cdot,\cdots,\cdot), c_i(\cdot,\cdots,\cdot)$ or $C_i(\cdot,\cdots,\cdot)(i\in
{\mathbb{N}})$ is used. For function spaces, $L^p(\mathbb{R})$($1\leq p\leq+\infty$) denotes   the standard  Lebesgue space with
the norm
\[\|f\|_{L^p(\mathbb{R})}=\left(\int_\mathbb{R}\left|f(x)\right|^pdx\right)^{\frac{1}{p}},\] and
 $H^l(\mathbb{R})$ is the usual $l$-th order Sobolev space with its norm
\[\|f\|_{l}=\left(\sum_{i=0}^{l}\|\partial_x^if\|^2\right)^{\frac{1}{2}}\quad \mathrm{with}\quad \|\cdot\|:=\|\cdot\|_{L^2(\mathbb{R})}.\]
Finally,  $\|\cdot\|_{L_{T,x}^\infty}$ stands for the norm $\|\cdot\|_{L^\infty([0,T]\times\mathbb{R})}$.

\section{Preliminaries}
\setcounter{equation}{0}
The viscous contact wave $(V,U,\Theta)(t,x)$ defined in (\ref{1.11})  has the following properties.
\begin{Lemma}[\cite{Hakho-Hong-2010}]
Let $|\theta_+-\theta_-|\le m_0\delta$, where $\delta=\gamma-1$ and $m_0>0$ is a positive constant independent of $\delta$, then it holds that
\begin{itemize}
\item[(i)]\quad $|V-v_{\pm}|+|\Theta-\theta_{\pm}|\le c_1\delta e^{-\frac{c_0x^2}{\delta(1+t)}}$,
\item[(ii)]\quad $|\partial_x^kV|+|\partial_x^k\Theta|\le c_2\delta^{\frac{2-k}{2}}(1+t)^{-\frac{k}{2}}e^{-\frac{c_0x^2}{\delta(1+t)}},\quad \forall \,k\in\mathbb{Z}^+,$
\item[(iii)]\quad $|\partial_x^{k-1}U|\le c_3\delta|\partial_x^k\Theta|\le c_4\delta^\frac{4-k}{2}(1+t)^{-\frac{k}{2}}e^{-\frac{c_0x^2}{\delta(1+t)}},\quad \forall\, k\in\mathbb{Z}^+$.
\end{itemize}
where $c_i, i=0,1,2,3,4$ are positive constants depending only on $\theta_{\pm}$.
\end{Lemma}

The following lemma on the heat kernel will play an important role in the analysis of this paper, whose proof can be found in \cite{F.M. Huang-J. Li-A. Matsumura-2010}.

For $\alpha\in(0,\frac{c_0}{4}]$, we define
\begin{equation}\label{2.1}
w(t,x)=(1+t)^{-\frac{1}{2}}\exp\left\{-\frac{\alpha x^2}{\delta(1+t)}\right\},\quad g(t,x)=\int_{-\infty}^x w(t,y)\,dy.
\end{equation}
Then it is easy to check that
\begin{equation}\label{2.2}
4\alpha g_t=\delta w_x,\qquad \|g(t,\cdot)\|_{L^\infty}=\sqrt{\pi}\alpha^{-1/2}\delta^{\frac{1}{2}},
\end{equation}
and  we have
\begin{Lemma}[\cite{F.M. Huang-J. Li-A. Matsumura-2010}]
For any $0<T\leq\infty$, suppose that $h(t,x)$ satisfies
\[h\in L^\infty(0,T; L^2(\mathbb{R})),\quad h_x\in L^2(0,T; L^2(\mathbb{R})),\quad h_t\in L^2(0,T; H^{-1}(\mathbb{R})).\]
Then the following estimate holds:
\begin{equation}\label{2.3}
\int_0^T\int_{\mathbb{R}}h^2w^2dxdt\leq4\pi\|h(x,0)\|^2+4\pi\delta\alpha^{-1}\int_0^T\|h_x(\tau)\|^2d\tau+\frac{8\alpha}{\delta}\int_0^T\langle h_t, hg^2\rangle\,d\tau,
\end{equation}
where $\langle\cdot,\cdot\rangle$ denotes the inner product on $H^{-1}\times H^1$.
\end{Lemma}

The solution $w(t,x)$ has the following properties.
\begin{Lemma}
For given $w_-\in\mathbb{R}$ and $\bar{w}>0$, let $w_+\in\{w|0<\tilde{w}\triangleq w-w_-<\bar{w}\}$. Then the problem $(\ref{1.26})$ has a unique global  smooth solution satisfying the following
\begin{itemize}
\item[(i)]$w_-<w(t,x)<w_+$, $w_x>0$, $x\in\mathbb{R}, t>0.$
\item[(ii)] For any $p\in[1,+\infty]$, there  exists some positive constant $C=C(p,w_-,\bar{w})$ such that for $\tilde{w}\geq0$ and $t\geq0$,
\[\|w_x(t)\|_{L^p}\leq C\min\{\tilde{w}, \tilde{w}^{\frac{1}{p}}t^{-1+\frac{1}{p}}\},\quad
\|w_{xx}(t)\|_{L^p}\leq C\min\{\tilde{w}, t^{-1}\}.\]
\item[(iii)] If $w_->0$, for any $(t,x)\in[0,+\infty)\times(-\infty,0]$,
\[\displaystyle|w(t,x)-w_-|\leq \tilde{w}e^{-2(|x|+w_-t)},\quad\displaystyle|w_x(t,x)|\leq 2\tilde{w}e^{-2(|x|+w_-t)}.
\]
\item[(iv)]If $w_+<0$, for any $(t,x)\in[0,+\infty)\times[0,+\infty)$,
\[\displaystyle|w(t,x)-w_+|\leq \tilde{w}e^{-2(|x|+|w_+|t)},\quad\displaystyle|w_x(t,x)|\leq 2\tilde{w}e^{-2(|x|+|w_+|t)}.
\]
\item[(v)] Let $w^{r}(\frac{x}{t})$ be the Riemann solution of the scalar equation $(\ref{1.26})_1$ with the Riemann initial data
\[w(0,x)=\left\{\begin{array}{c}
       w_-, \quad x<0,\\[2mm]
        w_+, \quad x>0,\end{array}
  \right.\]
  then we have \[\lim_{t\rightarrow+\infty}\sup_{x\in\mathbb{R}}\left|w(t,x)-w^r\left(\frac{x}{t}\right)\right|=0.\]
\end{itemize}
\end{Lemma}

In order to use Lemma 2.3 to study the properties of the smooth rarefaction waves $(V_{\pm}^r, U_{\pm}^r, \Theta_{\pm}^r)$ constructed in (\ref{1.25}) and the viscous contact wave  $(V^{c}, U^{c}, \Theta^{c})(t,x)$, we divided the the domain $\mathbb{R}\times(0,t)$ into three parts, that is $\mathbb{R}\times(0,t)=\Omega_{-}\cup\Omega_{c}\cup\Omega_{+}$ with
\[
\Omega_{\pm}=\{(x,t)|\pm2x>\pm\lambda_\pm(v_\pm^m,s_\pm)t\},\quad
\Omega_{c}=\{(x,t)|\lambda_-(v_-^m,s_-)t\leq2x\leq\lambda_+(v_+^m,s_+)t\}.
\]

\begin{Lemma}
Assume that (\ref{1.19}) holds. Then the smooth rarefaction waves $(V_{\pm}^r, U_{\pm}^r, \Theta_{\pm}^r)$ constructed in (\ref{1.25}) and the viscous contact wave  $(V^{c}, U^{c}, \Theta^{c})(t,x)$  satisfy the following
\begin{itemize}
\item[(i)]$ \left(U_{\pm}^r\right)_x\geq0$, $x\in\mathbb{R},t>0.$
\item[(ii)] For any $p\in[1,+\infty]$, there  exists a positive constant $C=C(p,v_-,u_-, \theta_-, \delta_1, m_0)$ such that for $\delta=\gamma-1$ and $t\geq0$,
\[ \aligned
&\left\|\left((V_{\pm}^r)_x, (U_{\pm}^r)_x, (\Theta_{\pm}^r)_x\right)(t)\right\|_{L^p}\leq C\min\{\delta, \delta^{\frac{1}{p}}t^{-1+\frac{1}{p}}\},\\
&\left\|\left(\partial_x^kV_{\pm}^r, \partial_x^kU_{\pm}^r, \partial_x^k\Theta_{\pm}^r\right)(t)\right\|_{L^p}\leq C\min\{\delta, t^{-1}\}, \quad k=2,3.
 \endaligned\]
\item[(iii)] There exists some positive constant $C=C(p,v_-,u_-, \theta_-, \delta_1, m_0)$ such that for $\delta=\gamma-1$ and
\[c_1=\frac{1}{20}\min\left\{\left|\lambda_-(v_-^m,s_-)\right|, \lambda_+(v_+^m,s_+), c_0\lambda_-^2(v_-^m,s_-),c_0\lambda_+^2(v_+^m,s_+), 1\right\},\]
we have in $\Omega_{c}$ that
\[
\left|\left((V_{\pm}^r)_x, (U_{\pm}^r)_x, (\Theta_{\pm}^r)_x\right)\right|+\left|V_{\pm}^r-v_{\pm}^m\right|+\left|\Theta_{\pm}^r-\theta_{\pm}^m\right|\leq C\delta e^{-c_1(|x|+t)},
\]
and in  $\Omega_{\mp}$,
\[
\displaystyle\left|V^{c}_x\right|+\left|\Theta^{c}_x\right|\leq C \displaystyle\delta^{\frac{1}{2}} e^{-c_1(|x|+t)},
\]
\[
\displaystyle\left|V^{c}-v_{\mp}^m\right|+\left|\Theta^{c}-\theta_{\mp}^m\right|+\left|U^{c}_x\right|\leq\displaystyle C\delta e^{-c_1(|x|+t)},
\]
\[
\left|\left((V_{\pm}^r)_x, (U_{\pm}^r)_x, (\Theta_{\pm}^r)_x\right)\right|+\left|V_{\pm}^r-v_{\pm}^m\right|+\left|\Theta_{\pm}^r-\theta_{\pm}^m\right|\leq C\delta e^{-c_1(|x|+t)}.
\]
\item[(iv)]  It holds that \[\lim_{t\rightarrow+\infty}\sup_{x\in\mathbb{R}}\left|(V_{\pm}^r, U_{\pm}^r, \Theta_{\pm}^r)(t,x)-(v_{\pm}^r, u_{\pm}^r, \theta_{\pm}^r)\left(\frac{x}{t}\right)\right|=0.\]
\end{itemize}
\end{Lemma}

\section{Proof of Theorem 1.1}
\setcounter{equation}{0}
This section is devoted to proving Theorem 1.1. To do so, we first  reformulate  the original problem,  and then perform energy  estimates on  solutions to the reformulated system.

\subsection{Reformulation of the  problem}
 First, we define the perturbation $(\phi,\psi,\zeta)(t,x)$ by
\[
  \phi(t,x)=v(t,x)-V(t,x),\, \psi(t,x)=u(t,x)-U(t,x),\, \zeta(t,x)=\theta(t,x)-\Theta(t,x),
\]
then it follows from (\ref{1.1}) and (\ref{1.13}) that
\begin{equation}\label{3.1}
  \left\{\begin{array}{ll}
           \phi_t-\psi_x=0, \\[2mm]
           \psi_t+ \displaystyle\left(\frac{R(\Theta+\zeta)}{v}-\frac{R\Theta}{V}\right)_x=\displaystyle\left(\frac{\mu(v,\theta)u_x}{v}-\frac{\mu(V,\Theta)U_x}{V}\right)_x+K_x-R_1, \\[2mm]
           \displaystyle\left(\frac{R}{\gamma-1}-\frac{\theta}{2}k_{\theta\theta}(v,\theta)\frac{v_x^2}{v^5}\right)\zeta_t+p(v,\theta)u_x-P(V,\Theta)U_x\\[3mm]
           =\displaystyle\left(\frac{\tilde{\alpha}(v,\theta)\theta_x}{v}-\frac{\tilde{\alpha}(V,\Theta)}{V}\Theta_x\right)_x+\frac{\mu(v,\theta)u_x^2}{v}
           -\frac{\mu(V,\Theta)U_x^2}{V}-R_2+F+\frac{\theta}{2}k_{\theta\theta}(v,\theta)\frac{v_x^2}{v^5}\Theta_t,
         \end{array}\right.
\end{equation}
where $K$ and $F$ are defined in (\ref{1.3}).
System (\ref{3.1}) is  supplemented with the following initial data and far-field end state:
\begin{equation}\label{3.2}
  \left\{\begin{array}{l}
    (\phi,\psi,\zeta)(0,x)=(\phi_0,\psi_0,\zeta_0)(x)=(v-V, u-U, \theta-\Theta)(0,x), \\[2mm]
    (\phi,\psi,\zeta)(t,\pm\infty)=0.
  \end{array}\right.
\end{equation}

 We seek  the solutions of  the Cauchy problem (\ref{3.1})-(\ref{3.2}) in the following set of functions:
\[\aligned
  &X(0,T;m_1,M_1,m_2,M_2)\\
  &=\left\{(\phi,\psi,\zeta)(t,x)\left|\begin{array}{c}
   (\phi,\psi,\zeta)(t,x)\in C(0,T;H^2(\mathbb{R})\times H^1(\mathbb{R})\times H^1(\mathbb{R})),\\[1mm]
   (\phi_x,\psi_x,\zeta_x)(t,x)\in L^2(0,T;H^2(\mathbb{R})\times H^1(\mathbb{R})\times H^1(\mathbb{R})), \\[1mm]
   m_1\le \phi(t,x)+V((t,x)\le M_1,\quad m_2\le \zeta(t,x)+\Theta(t,x)\le M_2,\\[1mm]
   \end{array}\right.\right\}
\endaligned\]
where $m_1,m_2, M_1, M_2$ and $0\le T\le+\infty$ are some positive constants. Then to prove Theorem 1.1, it suffices to show  the following theorem.
\begin{Theorem}
Under the assumptions of Theorem 1.1,  there exist two small positive constant $\varepsilon_0$ and $\delta_0$ depending  only on $\underline{V},$ $\overline{V},$  $\underline{\Theta}$,  $\overline{\Theta}$, $\|\phi_0\|_2$ and $\left\|\left(\psi_0, \frac{\zeta_0}{\sqrt{\delta}}\right)\right\|_{1}$  such that if $0<\varepsilon\leq\varepsilon_0$ and $0<\delta:=\gamma-1\leq\delta_0$, the Cauchy problem (\ref{3.1})-(\ref{3.2}) admits a unique global-in-time solution $(\phi, \psi, \zeta)(t,x)$ satisfying
\begin{equation}\label{3.3}
C_0^{-1}\leq v(t,x)\leq C_0,\quad \frac{\underline{\Theta}}{2}\leq \theta(t,x)\leq 2\overline{\Theta}, \quad \forall\,(t,x)\in[0,\infty)\times\mathbb{R},
\end{equation}

\begin{equation}\label{3.4}\aligned
&\|\phi(t)\|^2_{2}+\left\|\left(\psi, \frac{\zeta}{\sqrt{\delta}}\right)(t)\right\|^2_{1}+\int_0^t\left(\left\|\left(\phi_x, \psi_x,\zeta_x\right)(\tau)\right\|^2_{1}\right)d\tau\\
&\leq C_3\left(\|\phi_0\|^2_{2}+\left\|\left(\psi_0, \frac{\zeta_0}{\sqrt{\delta}}\right)\right\|^2_{1}\right), \quad \forall\,t>0,\endaligned
\end{equation}
and
\begin{equation}\label{3.4-1}\aligned
&\int_0^t\left\|\phi_{xxx}(\tau)\right\|^2d\tau
\leq C_4\left(1+\delta^{-\frac{1}{2}}\right), \quad \forall\,t>0.\endaligned
\end{equation}
Moreover, the following large-time behavior of solutions hold:
 \begin{equation}\label{3.5}
\lim_{t\rightarrow+\infty}\sup_{x\in\mathbb{R}}\left\{\left|\left(\phi,\psi, \zeta\right)(t,x)\right|\right\}=0.
\end{equation}
Here $C_0$ is a positive constants  depending  only on $\underline{V},$ $\overline{V},$  $\underline{\Theta}$,  $\overline{\Theta}$, $\|\phi_0\|_1$ and $\left\|\left(\psi_0, \frac{\zeta_0}{\sqrt{\delta}}\right)\right\|$,  and $C_3, C_4$ are  positive constants  depending  only on $\underline{V},$ $\overline{V},$  $\underline{\Theta}$,  $\overline{\Theta}$, $\|\phi_0\|_2$ and $\left\|\left(\psi_0, \frac{\zeta_0}{\sqrt{\delta}}\right)\right\|_{1}$.
\end{Theorem}

In order to prove Theorem 3.1,  we first give  the following local existence result.
\begin{Proposition} [Local existence]
 Under the assumptions of  Theorem 1.1, there exists a sufficiently small positive constant $t_1$ depending only on $\underline V,$ $\overline V,$ $\underline\Theta,$ $\overline\Theta$ and $\|\phi_0\|_2,$ $\|(\psi_0, \frac{\zeta_0}{\sqrt{\delta}})\|_1$ such that the Cauchy  problem (\ref{3.1})-(\ref{3.2}) admits a unique smooth solution $(\phi, \psi, \zeta)(t,x)\in X \left(0,t_1; \frac{\underline{V}}{2}, 2\overline{V}, \frac{\underline{\Theta}}{2}, 2\overline{\Theta}\right)$, and
 \[\aligned
   &\displaystyle \sup_{t\in[0,t_1]}\left\{\left\|\phi(t)\right\|^2_2+\left\|\left(\psi,\frac{\zeta}{\sqrt{\delta}}\right)(t)\right\|^2_1\right\}
  +\int_0^{t_1}\left(\|\phi_x(\tau)\|_2^2+\|(\psi_x,\zeta_x)(\tau)\|_1^2\right)d\tau\\
  &\le b\left(\left\|\phi_0\right\|^2_2+\left\|\left(\psi_0,\frac{\zeta_0}{\sqrt{\delta}}\right)\right\|^2_1\right),
 \endaligned\]
 where $b>1$ is a positive constant depending only on  $\underline{V}$ and $\overline{V}$.
 \end{Proposition}

Proposition 3.1 can be obtained  by using the dual argument and iteration technique,  the  proof of which is similar to that of Theorem 2.1 in
\cite{H. Hattori-D. Li-1996}  and  thus omitted here for brevity.

Suppose that the local solution $(\phi, \psi, \zeta)(t,x)$ obtained in Proposition 3.1 has been extended to the time step $t=T\geq t_1$ for some positive constant $T>0$. To prove  the global existence of solutions to the Cauchy problem  (\ref{3.1})-(\ref{3.2}),  by the standard continuation argument, we need to establish the following a priori estimates.
\begin{Proposition} [A priori estimates]
 Under the assumptions of Theorem 3.1, suppose that $(\phi, \psi, \zeta)\\(t,x)\in X(0, T; m_0, M_0, \Theta_0, \Theta_1)$ is a solution of the Cauchy problem (\ref{3.1})-(\ref{3.2}) for some  positive constants $T$, $m_0$, $M_0$, $\Theta_0,$ and $\Theta_1$, and satisfies the following a priori assumptions:
 \begin{equation}\label{3.10}
   \displaystyle \sup_{t\in[0,T]}\left\{\left\|\phi(t)\right\|^2_2+\left\|\left(\psi,\frac{\zeta}{\sqrt{\delta}}\right)(t)\right\|^2_1\right\}
  +\int_0^T\left\|(\phi_x,\psi_x,\zeta_x)(\tau)\right\|_1^2d\tau\le N_1^2,
\end{equation}
\begin{equation}\label{3.10-1}
   \displaystyle \int_0^T\|\phi_{xxx}(\tau)\|^2\,d\tau\le N_2^2.
\end{equation}
Then there exist a smooth positive  function $\Xi_1(m_0,M_0; \underline{V}, \overline{V}, \underline{\Theta}, \overline{\Theta},N_{01})$ which is increasing on both $(m_0)^{-1}$ and $M_0$, and a positive constant $\Xi_2(\underline{V}, \overline{V}, \underline{\Theta}, \overline{\Theta},N_{01})$  with $N_{01}:=\|\phi_0\|_1+\|(\psi_0,\frac{\zeta_0}{\sqrt{\delta}})\|$ such that if
\begin{equation}\label{3.10-2}
\left\{\begin{array}{ll}
    \displaystyle \Xi_1(m_0,M_0; \underline{V}, \overline{V}, \underline{\Theta}, \overline{\Theta},N_{01})N_1^{11}\delta^{\frac{1}{4}}<\frac{1}{3},\\[2mm]
 \displaystyle \Xi_2(\underline{V}, \overline{V}, \underline{\Theta}, \overline{\Theta},N_{01})N_1^{4}\varepsilon<\frac{1}{3},\\[2mm]
 N_2^2\delta^{\frac{3}{4}}<1,
  \end{array}\right.
\end{equation}
then the inequalities in (\ref{3.3})-(\ref{3.4-1}) hold for all $(t,x)\in [0,T]\times\mathbb{R}$.
\end{Proposition}
\noindent{\bf Proof of Theorem 3.1.} Based on Propositions 3.1-3.2,  we now use the continuation argument to extend the unique  local solution $(\phi, \psi, \zeta)(t,x)$ to be a global one, i.e., $T=+\infty$. First,  we have from Proposition 3.1 that  $(\phi, \psi, \zeta)(t,x)\in X(0, t_1; m_0, M_0, \Theta_0, \Theta_1)$ with $m_0=\frac{\underline{V}}{2}, M_0=2\overline{V}, \Theta_0=\frac{\underline{\Theta}}{2}, \Theta_1=2\overline{\Theta}
$ and the a priori assumption (\ref{3.10})-(\ref{3.10-1}) hold with
\[N_1=N_2=\sqrt{b}\left(\left\|\phi_0\right\|_2+\left\|\left(\psi_0,\frac{\zeta_0}{\sqrt{\delta}}\right)\right\|_1\right):=\sqrt{b}N_0
\]
for all $t\in[0,t_1]$, where $t_1>0$ is a small positive constant given in Proposition 3.1.
 Then it is easy to find two small positive constants  $\delta_1>0$ and $\varepsilon_1>0$ depending only on $\underline{V}, \overline{V}, \underline{\Theta}, \overline{\Theta}$ and $N_0$   such that
\begin{equation}\label{3.1111}
\left\{\begin{array}{ll}
    \displaystyle \Xi_1\left(\frac{\underline{V}}{2}, 2\overline{V}; \underline{V}, \overline{V}, \underline{\Theta}, \overline{\Theta},N_{01}\right)(\sqrt{b}N_0)^{11}\delta_1^{\frac{1}{4}}<\frac{1}{3},\\[2mm]
 \displaystyle \Xi_2(\underline{V}, \overline{V}, \underline{\Theta}, \overline{\Theta},N_{01})(\sqrt{b}N_0)^{4}\varepsilon_0<\frac{1}{3},\\[2mm]
(\sqrt{b}N_0)^2\delta_1^{\frac{3}{4}}<1.
  \end{array}\right.
\end{equation}
Thus if $0<\delta\leq\delta_1$ and $0<\varepsilon\leq\varepsilon_1$, then the inequalities in (\ref{3.3})-(\ref{3.4-1}) hold for all $(t,x)\in [0,t_1]\times\mathbb{R}$.

Now we take $(\phi, \psi, \zeta)(t_1,x)$ as initial data, then by Proposition 3.1, we can extend the local solution $(\phi, \psi, \zeta)(t,x)$ to the time step $t=t_1+t_2$ for some suitably small constant $t_2>0$ depending only on $\underline{V}, \overline{V}, \underline{\Theta}, \overline{\Theta}$ and $N_{0}$.  Moreover, $(\phi, \psi, \zeta)(t,x)\in X(t_1,t_1+t_2; m_0, M_0, \Theta_0, \Theta_1)$ with $m_0=\frac{C_0^{-1}}{2}, M_0=2C_0, \Theta_0=\frac{\underline{\Theta}}{4}, \Theta_1=4\overline{\Theta}
$ and the a priori assumption (\ref{3.10})-(\ref{3.10-1}) hold with
\[N_1=\sqrt{C_3}N_0,\quad N_2=\sqrt{C_4\left(1+\delta^{-\frac{1}{2}}\right)}
\]
for all $t\in[t_1, t_1+t_2]$.
 Then  there exist two small positive constants  $\delta_2>0$ and $\varepsilon_2>0$ depending only on $\underline{V}, \overline{V}, \underline{\Theta}, \overline{\Theta}$ and $N_0$  such that
\begin{equation}\label{3.1112}
\left\{\begin{array}{ll}
    \displaystyle \Xi_1\left(\frac{{C_0^{-1}}}{2}, 2C_0; \underline{V}, \overline{V}, \underline{\Theta}, \overline{\Theta},N_{01}\right)(\sqrt{C_3}N_0)^{11}\delta_2^{\frac{1}{4}}<\frac{1}{3},\\[2mm]
 \displaystyle \Xi_2(\underline{V}, \overline{V}, \underline{\Theta}, \overline{\Theta},N_{01})(\sqrt{C_3}N_0)^{4}\varepsilon_2<\frac{1}{3},\\[2mm]
C_4\left(1+\delta_2^{-\frac{1}{2}}\right)\delta_2^{\frac{3}{4}}<1.
  \end{array}\right.
\end{equation}
Consequently, if $0<\delta\leq\delta_2$ and $0<\varepsilon\leq\varepsilon_2$,  the inequalities in (\ref{3.3})-(\ref{3.4-1}) hold for all $(t,x)\in [t_1,t_1+t_2]\times\mathbb{R}$. Letting $\delta_0=\min\{\delta_1, \delta_2\}$ and $\varepsilon_0=\min\{\varepsilon_1, \varepsilon_2\}$, we see that if $0<\delta\leq\delta_0$ and $0<\varepsilon\leq\varepsilon_0$, then the local solution $(\phi, \psi, \zeta)(t,x)\in  X(0,t_1+t_2; C_0^{-1}, C_0, \frac{\underline{\Theta}}{2}, 2\overline{\Theta})$.

Next, taking $(\phi, \psi, \zeta)(t_1+t_2,x)$ as initial data and  exploiting Proposition 3.1 again, we can extend the local solution $(\phi, \psi, \zeta)(t,x)$ to the time step $t=t_1+2t_2$. By repeating the above procedure, we
can thus extend the local solution $(\phi, \psi, \zeta)(t,x)$ step by step to a global
one provided that $0<\delta<\delta_0$ and $0<\varepsilon<\varepsilon_0$. And as a by-product, the inequalities in (\ref{3.3})-(\ref{3.4-1}) hold for all $(t,x)\in [0,+\infty)\times\mathbb{R}$.

Finally, the estimates (\ref{3.4}) and the system (\ref{3.1}) imply that
\begin{equation}\label{3.8}
\int_0^{+\infty}\left(\|\left(\phi_x, \psi_x,\zeta_x\right)(t)\|^2+\left|\frac{d}{dt}\left(\|\left(\phi_x, \psi_x,\zeta_x\right)(t)\|^2\right)\right|\right)dt<\infty,
\end{equation}
which, together with (\ref{3.4}) and the Sobolev inequality:
\begin{equation}\label{3.9}
\|f(t)\|_{L^\infty}\leq\|f(t)\|^{\frac{1}{2}}\|f_x(t)\|^{\frac{1}{2}},\quad\forall\,f(t,\cdot)\in H^1(\mathbb{R})
 \end{equation}
leads to the  asymptotic behaviors (\ref{3.5}). This completes the proof of Theorem 3.1.

\subsection{Energy  estimates}
 In this subsection, we shall prove  Proposition 3.2.  First of all,  notice that (\ref{3.10}) implies  $\|\zeta(t)\|_1\leq N_1\delta^{\frac{1}{2}}$ for all $t\in[0,T]$, thus if $\delta>0$ is sufficiently small such that $ N_1\delta^{\frac{1}{2}}<\frac{\underline\Theta}{2}$, then we have
\begin{equation}\label{3.11}
  |(\theta-\Theta)(t,x)|\le\|\zeta(t)\|_{L^{\infty}(\mathbb{R})}\le\|\zeta(t)\|^\frac{1}{2}\|\zeta_x(t)\|^\frac{1}{2}\le\|\zeta(t)\|_{1}\le\frac{\underline\Theta}{2},
  \quad\forall \,(t,x)\in [0,T]\times\mathbb{R}.
\end{equation}
Consequently
\begin{equation}\label{3.12}
  \frac{\underline\Theta}{2}\le\Theta(t,x)-\frac{\underline\Theta}{2}\le\theta(t,x)\le\Theta(t,x)+\frac{\underline\Theta}{2}\le2\overline\Theta,\quad \forall(t,x)\in[0,T]\times \mathbb{R}.
\end{equation}

Throughout of this subsection, we always assume  $N_1\delta^{\frac{1}{2}}<\frac{\underline\Theta}{2}$ so that (\ref{3.11})-(\ref{3.12})  hold. Moreover, we  denote
\[N_0:=\|\phi_0\|_2+\left\|\left(\psi_0,\frac{\zeta_0}{\sqrt{\delta}}\right)\right\|_1, \quad N_{01}:=\|\phi_0\|_1+\left\|\left(\psi_0,\frac{\zeta_0}{\sqrt{\delta}}\right)\right\|,\]
and assume that $N_1\geq N_{0}\gg1,  N_{01}\gg1$ and $\delta<1$ without loss of generality.

 Proposition 3.2 will obtained by a  series of lemmas below. The following basic energy estimates is key for the proof of Proposition 3.2.
\begin{Lemma}[Basic energy estimates]
 Under the assumptions of Proposition 3.2, there exist a positive constant $C(\underline V,\overline V,\underline\Theta,\overline \Theta)$ and a positive constant $C_5$ depending only on $\underline V,\overline V,\underline \Theta,\overline \Theta, m_0, M_0$ such that
\begin{eqnarray} \label{3.13}
     && \quad\int_{\mathbb{R}}\left[R\Theta\Phi\left(\frac{v}{V}\right)+\frac{\psi^2}{2}+\frac{R}{\delta}\Theta\Phi\left(\frac{\theta}{\Theta}\right)+\frac{\kappa(v,\theta)v_x^2}{v^5}\right] dx +\int_{0}^{t}\int_{\mathbb{R}}\left(\frac{\mu(v,\theta)\Theta\psi_x^2}{\theta v}+\frac{\tilde{\alpha}(v,\theta)\Theta\zeta_x^2}{v\theta^2}\right)dxd\tau \nonumber\\
      && \le C(\underline V,\overline V,\underline\Theta,\overline \Theta)\left\|\left(\phi_0,\phi_{0x}, \psi_0, \frac{\zeta_0}{\sqrt{\delta}}\right)\right\|^2 +C_5\delta\int_0^t (1+\tau)^{-1}\int_{\mathbb{R}}\left(\phi^2+\frac{\zeta^2}{\delta}\right)e^{-\frac{c_0x^2}{\delta(1+\tau)}}\ dxd\tau \nonumber\\
      &&\quad +C_5\left(N_1^6\delta^{\frac{1}{2}}\int_0^t\|(\phi_x,\phi_{xx},\psi_x,\zeta_x)(\tau)\|^2\ d\tau+N_1^5\delta\right).
\end{eqnarray}

 \end{Lemma}
\noindent{\bf Proof.}~~Multiplying $(\ref{3.1})_1$ by $-R\Theta(v^{-1}-V^{-1})$, $(\ref{3.1})_2$ by $\psi$, and  $(\ref{3.1})_3$ by $\displaystyle \zeta\theta^{-1}$,  then adding the resulting equations together, we have
\begin{equation}\label{3.14}
  \begin{split}
     & \left\{R\Theta\Phi\left(\frac{v}{V}\right)+\frac{\psi^2}{2}+\frac{R}{\delta}\Theta\Phi\left(\frac{\theta}{\Theta}\right)\right\}_t+E_x+\frac{\mu(v,\theta)\Theta\psi_x^2}{v\theta}+\frac{\tilde{\alpha}(v,\theta)\Theta\zeta_x^2}{v\theta^2} \\
     = & Q_0+Q_1+Q_2+Q_3,
  \end{split}
\end{equation}
where
\begin{equation}\label{3.15}
  \begin{split}
     E& =(p-p_+)\psi-\left(\frac{\mu(v,\theta)}{v}u_x-\frac{\mu(V,\Theta)}{V}U_x\right)\psi-\left(\frac{\tilde{\alpha}(v,\theta)\theta_x}{v}-\frac{\tilde{\alpha}(V,\Theta)}{V}\Theta_x\right)\frac{\zeta}{\theta}, \\
     Q_0 & =-p_+\Phi\left(\frac{V}{v}\right)U_x-\frac{p_+}{\delta}\Phi\left(\frac{\Theta}{\theta}\right)U_x-\left(\frac{\mu(v,\theta)}{v}-\frac{\mu(V,\Theta)}{V}\right)U_x\psi_x, \\
     & \quad+(p_+-p)\frac{\zeta}{\theta}U_x+\frac{\tilde{\alpha}(v,\theta)\zeta_x\zeta\Theta_x}{\theta^2}-\left(\frac{\tilde{\alpha}(v,\theta)}{v}-\frac{\tilde{\alpha}(V,\Theta)}{V}\right)\left(\frac{\zeta_x}{\theta}-\frac{\zeta\theta_x}{\theta^2}\right)\Theta_x \\
     &\quad +\left(\frac{\mu(v,\theta)(U_x^2+2\psi_xU_x)}{v}-\frac{\mu(V,\Theta)}{V}U_x^2\right)\frac{\zeta}{\theta}, \\
    Q_1&=K_x\psi, \quad  Q_2=-R_1\psi+\frac{F}{\theta}\zeta-R_2\frac{\zeta}{\theta}, \quad Q_3=\frac{\kappa_{\theta\theta}(v,\theta)}{2v^5}v_x^2\theta_t\zeta,\quad\\
     \Phi(s)&=s-1-\ln s.
  \end{split}
\end{equation}
Integrating (\ref{3.14})  over $[0,t]\times\mathbb{R}$ yields
\begin{equation}\label{3.16}
  \begin{split}
      & \int_{\mathbb{R}}\left\{R\Theta\Phi\left(\frac{v}{V}\right)+\frac{\psi^2}{2}+\frac{R}{\delta}\Theta\Phi\left(\frac{\theta}{\Theta}\right)\right\}dx+\int_0^t\int_{\mathbb{R}}\left(\frac{\mu(v,\theta)\Theta}{\theta v}\psi_x^2+\frac{\tilde{\alpha}(v,\theta)\Theta}{v\theta^2}\zeta_x^2\right)dxd\tau \\
       =& \int_{\mathbb{R}}\left[R\Theta_0\Psi\left(\frac{v_0}{V_0}\right)
       +\frac{\psi_0^2}{2}+\frac{R}{\delta}\Theta_0\Phi\left(\frac{\theta_0}{\Theta_0}\right)\right]\,dx+\int_0^t\int_{\mathbb{R}}(Q_0+Q_1+Q_2+Q_3)\,dxd\tau,
   \end{split}
\end{equation}
where $V_0=V(x,t)|_{t=0},\ \Theta_0=\Theta(t,x)|_{t=0}$.

By the convexity  of $\Phi(s)$ and the Cauchy inequality,  we have
\begin{equation}\label{3.17}
  \begin{split}
     |Q_0| & \le C(\underline V,\overline V,\underline \Theta,\overline \Theta,m_0,M_0)\left(|\phi^2U_x|+\left|\frac{\zeta^2}{\delta}U_x\right|+|(\phi,\zeta)U_x\psi_x|+|(\phi,\zeta)\zeta U_x|\right. \\
       & \left.\quad+|\zeta_x\zeta\Theta_x|+|(\phi,\zeta)\Theta_x||(\Theta\zeta_x,\zeta\Theta_x)|+|(\phi,\zeta)\zeta||U_x^2|+|\psi_xU_x\zeta|\right)\\
       & \le \frac{\mu(v,\theta)\Theta}{4v}\psi_x^2+\frac{\tilde{\alpha}(v,\theta)\Theta}{4v\theta^2}\zeta_x^2+C(\underline V,\overline V,\underline \Theta,\overline \Theta,m_0,M_0)\left(\phi^2+\frac{\zeta^2}{\delta}\right)(|U_x|+\Theta_x^2|),
   \end{split}
\end{equation}
thus
\begin{equation}\label{3.18}
  \begin{split}
    \left|\int_0^t\int_{\mathbb{R}}Q_0\ dxd\tau\right| & \le \frac{1}{4}\int_0^t\int_{\mathbb{R}}\frac{\mu(v,\theta)\Theta}{v\theta}\psi_x^2\   dxd\tau+\frac{1}{4}\int_0^t\int_{\mathbb{R}}\frac{\tilde{\alpha}(v,\theta)\Theta}{v\theta^2}\zeta_x^2\,dxd\tau \\
      &\quad +C(\underline V,\overline V,\underline \Theta,\overline \Theta,m_0,M_0)\int_0^t\int_{\mathbb{R}}\left(\phi^2+\frac{\zeta^2}{\delta}\right)(|U_x|+\Theta_x^2)\,dxd\tau.
  \end{split}
\end{equation}
Using $(\ref{3.1})_1$, we have by a direct computation that
\begin{equation}\label{3.19}
 \begin{split}
     Q_1 & =K_x\psi=\{K\psi\}_x-K\psi_x \\
       & =\{\cdots\}_x-\left\{-\frac{\kappa(v,\theta)v_{xx}}{v^5}+\frac{5\kappa(v,\theta)-v\kappa_{v}(v,\theta)}{2v^6}v_x^2-\frac{\kappa_{\theta}(v,\theta)v_x\theta_x}{v^5}\right\}\psi_x\\
       & =\{\cdots\}_x-\left\{-\frac{\kappa(v,\theta)v_{xx}}{v^5}+\frac{5\kappa(v,\theta)-v\kappa_{v}(v,\theta)}{2v^6}v_x^2-\frac{\kappa_{\theta}(v,\theta)v_x\theta_x}{v^5}\right\}\phi_t\\
       & =\{\cdots\}_x-\left\{-\frac{\kappa(v,\theta)v_{xx}}{v^5}+\frac{5\kappa(v,\theta)-v\kappa_{v}(v,\theta)}{2v^6}v_x^2-\frac{\kappa_{\theta}(v,\theta)v_x\theta_x}{v^5}\right\}v_t\\
       & \quad+\left\{-\frac{\kappa(v,\theta)v_{xx}}{v^5}+\frac{5\kappa(v,\theta)-v\kappa_{v}(v,\theta)}{2v^6}v_x^2-\frac{\kappa_{\theta}(v,\theta)v_x\theta_x}{v^5}\right\}V_t\\
       &=\{\cdots\}_x-\left(\frac{\kappa(v,\theta)v_x^2}{2v^5}\right)_t+\frac{\kappa_{\theta}\theta_tv_x^2}{2v^5}\\
       &\quad +\left\{-\frac{\kappa(v,\theta)v_{xx}}{v^5}+\frac{5\kappa(v,\theta)-v\kappa_{v}(v,\theta)}{2v^6}v_x^2-\frac{\kappa_{\theta}(v,\theta)v_x\theta_x}{v^5}\right\}V_t.
   \end{split}
\end{equation}
Here and hereafter, $\{\cdots\}_x$ denotes the terms which will disappear after integrating with respect to $x$.
Thus we have
\begin{equation}\label{3.20}
  \begin{split}
     \int_0^t\int_{\mathbb{R}}Q_1\,dxd\tau =&-\int_{R}\frac{\kappa(v,\theta)v_x^2}{2v^5}\,dx
     +\int_{\mathbb{R}}\frac{\kappa(v_0,\theta_0)v_{0x}^2}{2v_0^5}\,dx+\int_0^t\int_{\mathbb{R}}\frac{\kappa_{\theta}\theta_tv_x^2}{2v^5}\,dxd\tau  \\
       & + \int_0^t\int_{\mathbb{R}}\left\{-\frac{\kappa(v,\theta)v_{xx}}{v^5}+\frac{5\kappa(v,\theta)-v\kappa_{v}(v,\theta)}{2v^6}v_x^2-\frac{\kappa_{\theta}(v,\theta)v_x\theta_x}{v^5}\right\}{ V}_t\,dxd\tau\\
       :=&-\int_{\mathbb{R}}\frac{\kappa(v,\theta)v_x^2}{2v^5}\,dx+\int_{\mathbb{R}}\frac{\kappa(v_0,\theta_0)v_{0x}^2}{2v_0^5}\,dx+I_1+I_2.
   \end{split}
\end{equation}
For $I_1$,  notice that  $(\ref{3.3})_3$ implies
\begin{eqnarray}\label{3.21}
     \frac{\kappa_{\theta}\theta_tv_x^2}{2v^5} && = \frac{\kappa_{\theta}v_x^2}{2v^5}\frac{1}{C_v-\frac{\theta}{2}\kappa_{\theta\theta}(v,\theta)\frac{v_x^2}{v^5}}\left[-pu_x+\left(\frac{\tilde{\alpha}(v,\theta)\theta_x}{\theta}\right)_x+\frac{\mu(v,\theta)u_x^2}{v}+F\right] \nonumber\\
       && =\frac{\kappa_{\theta}v_x^2}{2v^5}\frac{1}{C_v-\frac{\theta}{2}\kappa_{\theta\theta}\frac{v_x^2}{v^5}}\left(-pu_x+\frac{\mu(v,\theta)u_x^2}{v}+\frac{v\kappa_{\theta  v}-\kappa_{\theta}}{2v^6}\theta u_xv_x^2\right)\nonumber\\
       && \quad-\frac{1}{2}\left(\frac{\theta\kappa_{\theta}^2v_x^3}{(C_v-\frac{\theta}{2}\kappa_{\theta\theta}\frac{v_x^2}{v^5})v^{10}}\right)_x u_x-\frac{\tilde{\alpha}\theta_x\left(2\kappa_{\theta}v_xv_{xx}+\kappa_{\theta v}v_x+\kappa_{\theta\theta}\theta_xv_x^2\right)}{2v^6(C_v-\frac{\theta}{2}\kappa_{\theta\theta}\frac{v_x^2}{v^5})}
       \nonumber\\
       && \quad+\frac{5\kappa_{\theta}\tilde{\alpha}\theta_xv_x^3}{2v^7(C_v-\frac{\theta}{2}\kappa_{\theta\theta}\frac{v_x^2}{v^5})}-\frac{\kappa_{\theta}v_x^2}{2v^5}\frac{1}{(C_v-\frac{\theta}{2}\kappa_{\theta\theta}\frac{v_x^2}{v^5})^2}
       \left[\frac{\kappa_{\theta\theta}\theta_xv_x^2}{2v^5}+\frac{\theta\kappa_{\theta\theta v}v_x^3}{2v^5}\right.\nonumber\\
       && \quad\left.+\frac{\theta\kappa_{\theta\theta\theta}\theta_xv_x^2}{2v^5}+\frac{\theta\kappa_{\theta\theta}v_xv_{xx}}{v^5}
       -\frac{5\theta\kappa_{\theta\theta}v_x^3}{2v^6}\right]\frac{\tilde{\alpha}\theta_x}{v}+\{\cdots\}_x\\
       && \le C\delta\left(|v_x^2u_x|+|v_x^2u_x^2|+|u_xv_x^4|+|v_xv_{xx}\theta_x|+|v_x^3\theta_x|+|\theta_x^2v_x^2|+|\theta_x^2v_x^4|+|\theta_xv_x^5|\right.\nonumber\\
       && \quad\left.+|\theta_xv_x^3v_{xx}|+\left(|\theta_xv_x^3|+|v_x^4|+|v_x^2v_{xx}|+|\theta_xv_x^5|+|v_x^6|+|v_x^4v_{xx}|\right)|u_x|\right)+\{\cdots\}_x,\nonumber
\end{eqnarray}
where we have used $\frac{1}{C_v-\frac{\theta}{2}\kappa_{\theta\theta}\frac{v_x^2}{v^5}}\le\frac{1}{C_v}\le C\delta$.
Therefore,
\begin{equation}\label{3.22}
  \begin{split}
    I_1 & \le C\delta\int_0^t\int_{\mathbb{R}}|v_x^2u_x|\left(1+v_x^2+v_x^4\right)dxd\tau \\
      &\quad+C\delta\int_0^t\int_{\mathbb{R}}(v_x^2+1)\left(|v_xv_{xx}\theta_x|+|v_x^3\theta_x|+|\theta_x^2v_x^2|+|\theta_xu_xv_x^3|\right)dxd\tau \\
      &\quad+C\delta\int_0^t\int_{\mathbb{R}}\left(v_x^2u_x^2+(v_x^2+v_x^4)|u_xv_{xx}|\right)dxd\tau \\
      &:=I^1_{1}+I^2_{1}+I^3_{1}.
        \end{split}
    \end{equation}
It follows from the Cauchy inequality, Young  inequality, Sobolev  inequality  and Lemma 2.1 that
  \begin{eqnarray}\label{3.23}
     I_{1}^1 && \le C\delta\int_0^t\left(1+\|v_x\|_{L^{\infty}}^2+\|v_x(\tau)\|_{L^{\infty}}^4\right)\int_{\mathbb{R}}\left(\|\phi_x(\tau)\|_{L^{\infty}}\phi_x\psi_x+\phi_x^2U_x^2+V_x^4+\psi_x^2+V_x^2U_x\right) dxd\tau\nonumber\\
     &&\le C N_1^5\delta\int_0^t\|(\phi_x,\psi_x)(\tau)\|^2\ d\tau+O(1)\delta^{\frac{7}{2}}N_1^4,
\end{eqnarray}
  \begin{eqnarray}\label{3.24}
     I_{1}^2 && \le C\delta\int_0^t\left(\|v_x(\tau)\|_{L^{\infty}}^2+1\right)\int_{\mathbb{R}}\left\{\|\phi_x(\tau)\|_{L^{\infty}}\phi_{xx}\zeta_x+|\phi_x\zeta_xV_{xx}|+|\phi_{xx}\zeta_xV_x|+|\zeta_xV_xV_{xx}|\right.\nonumber\\
     && \quad+\|\phi_x(\tau)\|_{L^{\infty}}^2|\phi_x\zeta_x|+|\zeta_xV_x^3|+|\phi_x\phi_{xx}\Theta_x|+|\phi_xV_{xx}\Theta_x|+|\phi_{xx}V_x\Theta_x|+|V_xV_{xx}\Theta_x|\nonumber\\
     && \quad+\|\phi_x(\tau)\|_{L^{\infty}}^2\phi_x^2|\Theta_x|+|V_x^3\Theta_x|+\|\phi_x(\tau)\|_{L^{\infty}}^2\zeta_x^2+\zeta_x^2V_x^2+\Theta_x^2\phi_x^2+\Theta_x^2V_x^2+\|\phi_x(\tau)\|_{L^{\infty}}\zeta_x\psi_x\nonumber\\
    &&\quad\left. +|\zeta_x\phi_xU_x|+|\phi_x\psi_x\Theta_x|+|\phi_x\Theta_xU_x|+|\zeta_x\psi_xV_x|+|\zeta_xU_xV_x|+|\psi_x\Theta_xV_x|+|\Theta_xU_xV_x|\right\}\,dxd\tau\nonumber\\
    && \le C N_1^4\delta \int_0^t\|(\phi_x,\phi_{xx},\psi_x,\zeta_x)(\tau)\|^2\,d\tau+CN_1^2\delta^{\frac{5}{2}},
\end{eqnarray}
  \begin{eqnarray}\label{3.25}
    I_{1}^3 && \le C\delta\int_0^t\int_{\mathbb{R}}\{\|\phi_x(\tau)\|_{L^{\infty}}^2\psi_x^2+\phi_x^2U_x^2+V_x^2\psi_x^2+V_x^2U_x^2\nonumber\\
    && +(\|v_x(\tau)\|_{L^{\infty}}^2+\|v_x(\tau)\|_{L^{\infty}}^4)(|\psi_x\phi_{xx}|+|\psi_xV_{xx}|+|U_x\phi_{xx}|+|U_xV_{xx}|)\}\,dxd\tau \nonumber\\
    &&\le C N_1^4\delta \int_0^t\|(\phi_{xx},\psi_x,\phi_x)(\tau)\|^2\,d\tau+CN_1^4\delta^{\frac{5}{2}}.
\end{eqnarray}
Combining (\ref{3.22})-(\ref{3.25}) yields
\begin{equation}\label{3.26}
  I_1\le C N_1^5\delta\int_0^t\|(\phi_{xx},\psi_x,\zeta_x,\phi_x)(\tau)\|^2\,d\tau+CN_1^4\delta^{\frac{5}{2}}.
\end{equation}
Similarly, it holds
\begin{equation}\label{3.27}
  \begin{split}
    I_2 & \le C\delta\int_0^t\int_{\mathbb{R}}\left(|\phi_{xx}|+|V_{xx}|+|\phi_x^2|+V_x^2+\zeta_x^2+\Theta_x^2\right)|\Theta_t|\,dxd\tau \\
      & \le C\delta\int_0^t\|(\phi_{xx},\phi_x,\zeta_x)(\tau)\|^2\,d\tau+C\delta^{\frac{3}{2}}.
  \end{split}
\end{equation}
Substituting (\ref{3.26}) and (\ref{3.27}) into (\ref{3.20}), we have
\begin{equation}\label{3.28}
  \begin{split}
    \int_0^t\int_{\mathbb{R}}Q_1\,dxd\tau \le& -\int_{R}\frac{\kappa(v,\theta)v_x^2}{2v^5}\,dx+\int_{\mathbb{R}}\frac{\kappa(v_0,\theta_0)v_{0x}^2}{2v_0^5}\,dx \\
      & +C N_1^5\delta\int_0^t\|(\phi_{xx},\phi_x,\psi_x,\zeta_x)(\tau)\|^2\,d\tau+CN_1^4\delta^{\frac{3}{2}}.
  \end{split}
\end{equation}
By using the estimate of $I_1$ and the fact that $\|\zeta(t)\|_{L^\infty_{T,x}}\leq\displaystyle\sup_{t\in[0,T]}\|\zeta(t)\|_1\leq N_1\sqrt{\delta}$,
\begin{equation}\label{3.29}
  \begin{split}
    \left|\int_0^t\int_{\mathbb{R}}Q_3\,dxd\tau\right| & \le C\left|\int_0^t\int_{\mathbb{R}}\frac{\kappa_{\theta}\theta_tv_x^2}{2v^5}\,dxd\tau\right|\cdot\|\zeta\|_{L_{T,x}^{\infty}}\\
      & \le C N_1^6\delta^{\frac{3}{2}}\int_0^t\|(\phi_{xx},\phi_x,\psi_x)(\tau)\|^2\,d\tau+CN_1^5\delta^{3}.
  \end{split}
\end{equation}
Finally,  for the estimate of $\displaystyle\int_0^t\int_{\mathbb{R}}Q_2\,dxd\tau$, notice that
\begin{equation*}
  \begin{split}
     \frac{F}{\theta}\zeta =&\left(\frac{\kappa_{\theta}v_xu_{xx}}{v^5}+\frac{v\kappa_{\theta v}(v,\theta)-\kappa_{\theta}(v,\theta)}{2v^6}u_xv_x^2\right)\zeta \\
       =&\{\cdots\}_x-\left(\frac{\kappa_{\theta}v_x}{v^5}\zeta\right)_xu_x+\frac{v\kappa_{\theta v}(v,\theta)-\kappa_{\theta}(v,\theta)}{2v^6}u_xv_x^2\zeta\\
       =&\{\cdots\}_x-v^{-5}(\kappa_{\theta v}v_x^2\zeta+\kappa_{\theta\theta}\theta_xv_x\zeta+\kappa_{\theta}v_{xx}\zeta+\kappa_{\theta}v_x\zeta_x)u_x\\
       &+5v^{-6}\kappa_{\theta}v_x^2\zeta u_x+\frac{1}{2}v^{-6}[v\kappa_{\theta v}(v,\theta)-\kappa_{\theta}(v,\theta)]v_x^2\zeta u_x\\
=&\{\cdots\}_x+O(1)\left(|v_x^2u_x\zeta|+|\theta_xv_xu_x\zeta|+|v_{xx}u_x\zeta|+|v_x\zeta_xu_x|\right),
   \end{split}
\end{equation*}
then similar to the estimates of (\ref{3.23})-(\ref{3.25}), we obtain
\begin{eqnarray}\label{3.30}
     \left|\int_0^t\int_{\mathbb{R}}\frac{F}{\theta}\zeta\,dxd\tau\right| &\le& C\|\zeta\|_{L_{T,x}^{\infty}}\int_0^t\int_{\mathbb{R}}(|v_x^2u_x|+|\theta_xv_xu_x|+|v_{xx}u_x|)\,dxd\tau \nonumber\\
       &&+C\int_0^t\int_{\mathbb{R}}|v_x\zeta_xu_x|\,dxd\tau\nonumber\\
       &\le& C N_1\delta^{\frac{1}{2}}\int_0^t\int_{\mathbb{R}}(u_x^2+v_x^4+{\Theta}_x^4+v_{xx}^2)\,dxd\tau+C\int_0^t\int_{\mathbb{R}}|\zeta_x||v_xu_x|\,dx\nonumber\\
       &\le& C N_1\delta^{\frac{1}{2}}\int_0^t\int_{\mathbb{R}}(\psi_x^2+{U}_x^2+\phi_x^4+{V}_x^4+{ \Theta}_x^4+\phi_{xx}^2+{ V}_{xx}^2)\,dxd\tau\nonumber\\
       &&+C\int_0^t\int_{\mathbb{R}}(|\zeta_x\phi_x\psi_x|+|\zeta_x\phi_x{ U}_x|+|\zeta_x{ V}_x\psi_x|+|\zeta_x{ V}_x{ U}_x|)\,dxd\tau\nonumber\\
&\le& C N_1\delta^{\frac{1}{2}}\int_0^t\int_{\mathbb{R}}\left(\psi_x^2+{ U}_x^2+\|\phi_x\|_{L^{\infty}}^2\phi_x^2+{ V}_x^4+\phi_{xx}^2+{\Theta}_x^4+{ V}_{xx}^2\right)\,dxd\tau\nonumber\\
       &&+C\int_0^t\left\{\|\phi_x(\tau)\|_{L^{\infty}}\|\zeta_x(\tau)\|\|\psi_x(\tau)\|+\|{ U}_x\|_{L^{\infty}}\|(\zeta_x,\phi_x)(\tau)\|^2\right.\nonumber\\
       &&\left. \quad+\|{V}_x\|_{L^{\infty}}\|(\zeta_x,\psi_x)(\tau)\|
       +\|{V}_x\|_{L^{\infty}}\|\zeta_x(\tau)\|\|{ U}_x(\tau)\|\right\}\,d\tau\nonumber\\
&\le&C N_1^3\delta^{\frac{1}{2}} \int_0^t\|(\psi_x,\phi_x,\zeta_x,\phi_{xx})(\tau)\|^2\,d\tau+C N_1\delta.
\end{eqnarray}
On the other hand, we derive from (\ref{1.13-1}) that
\begin{equation} \label{3.31}
  \begin{split}
    & \left|\int_0^t\int_{\mathbb{R}}\left(-R_1\psi-\frac{R_2}{\theta}\zeta\right)dxd\tau\right| \\
      \le & C\left|\int_0^t\int_{\mathbb{R}}(|R_1\psi|+|R_2\zeta|)\,dxd\tau\right|\\
      \le & C\int_0^t\left(\|\psi(\tau)\|^{\frac{1}{2}}\|\psi_x(\tau)\|^{\frac{1}{2}}\|R_1\|_{L^1}+\|\zeta(\tau)\|^{\frac{1}{2}}\|\zeta_x(\tau)\|^{\frac{1}{2}}\|R_2\|_{L^1}\right)\,d\tau\\
      \le & \frac{1}{4}\int_0^t\int_{\mathbb{R}}\frac{\mu(v,\theta){\Theta}}{v\theta}\psi_x^2\,dxd\tau+\frac{1}{4}\int_0^t\int_{\mathbb{R}}\frac{\tilde{\alpha}(v,\theta){ \Theta}}{v\theta^2}\zeta_x^2\,dxd\tau \\
      & +C(m_0,M_0,\underline {\Theta},\overline { \Theta})\int_0^t\left(\|\psi(\tau)\|^{\frac{2}{3}}\|R_1\|_{L^1}^{\frac{4}{3}}+\|\zeta(\tau)\|^{\frac{2}{3}}\|R_2\|_{L^1}^{\frac{4}{3}}\right)\,d\tau\\
  \le& \frac{1}{4}\int_0^t\int_{\mathbb{R}}\frac{\mu(v,\theta){\Theta}}{v\theta}\psi_x^2\,dxd\tau
  +\frac{1}{4}\int_0^t\int_{\mathbb{R}}\frac{\tilde{\alpha}(v,\theta){\Theta}}{v\theta^2}\zeta_x^2\,dxd\tau+CN_1^{\frac{2}{3}}\delta^{\frac{4}{3}},
  \end{split}
\end{equation}
Consequently,
\begin{equation}\label{3.32}
  \begin{split}
    \left|\int_0^t\int_{\mathbb{R}}Q_2\,dxd\tau\right| \le& \left|\int_0^t\int_{\mathbb{R}}\left(-R_1\psi-\frac{R_2}{\theta}\zeta\right)dxd\tau\right|
    +\left|\int_0^t\int_{\mathbb{R}}\frac{F}{\theta}\,dxd\tau\right| \\
      \le &\frac{1}{4}\int_0^t\int_{\mathbb{R}}\frac{\mu(v,\theta){\Theta}}{\theta v}\psi_x^2\,dxd\tau+\frac{1}{4}\int_0^t\int_{\mathbb{R}}\frac{\tilde{\alpha}(v,\theta){ \Theta}}{v\theta^2}\zeta_x^2\,dxd\tau\\
      & +C N_1^3\delta^{\frac{1}{2}}\int_0^t\|(\phi_x,\phi_{xx},\psi_x,\zeta_x)(\tau)\|^2\,d\tau+CN_1\delta.
  \end{split}
\end{equation}
Putting (\ref{3.18}), (\ref{3.28}), (\ref{3.29}), (\ref{3.32}) into (\ref{3.16}), we get (\ref{3.13}) immediately. This finishes  the proof of Lemma 3.1.

For the remainder term  $\int_0^t (1+\tau)^{-1}\int_{\mathbb{R}}(\phi^2+\frac{\zeta^2}{\gamma-1})e^{-\frac{c_0x^2}{\delta(1+\tau)}}\ dxd\tau$ in (\ref{3.13}), we establish the following:

\begin{Lemma}
 Under the assumptions of Proposition 3.2, there exist two positive constants $C_6, C_7>1$ depending only on $\underline V,\overline V,\underline \Theta,\overline \Theta, m_0, M_0$ such that if
\begin{equation}\label{3.33-1}
C_6N_1^2\delta^{\frac{1}{2}}<\frac{1}{2}\min\left\{2p_+^2, R^2, \underline{\Theta}\right\},
\end{equation}
 then it holds for all $t\in[0,T]$ that
\begin{equation}\label{3.33}
  \int_0^t\int_{\mathbb{R}}\left(\phi^2+\psi^2+\frac{\zeta^2}{\delta}\right)w^2\,dxd\tau\le C_7 N_1^6\delta^{-\frac{3}{4}}+C_7N_1^7\delta^{-\frac{3}{4}}\int_0^t\|(\phi_x,\phi_{xx},\psi_x,\zeta_x)(\tau)\|^2\,d\tau.
\end{equation}
 \end{Lemma}

The proof of Lemma 3.2 is given in the Appendix, which is technique but similar to that of Lemma 5 in \cite{F.M. Huang-J. Li-A. Matsumura-2010}.

Next, we estimate $\left\|\frac{\mu(v,\theta)\phi_x}{v}(t)\right\|$.
\begin{Lemma}
 Under the assumptions of Proposition 3.2, there exist a positive constant $C(\underline V,\overline V,\underline\Theta,\overline \Theta)$ and a positive constant $C_9$ depending only on $\underline V,\overline V,\underline \Theta,\overline \Theta, m_0, M_0$ such that
 \begin{equation}\label{3.67}
  \begin{split}
     & \int_{\mathbb{R}}\frac{\mu^2(v,\theta)}{v^2}\phi_x^2\,dx
     +\int_0^t\int_{\mathbb{R}}\frac{\mu(v,\theta)\theta}{v^3}\phi_x^2\,dxd\tau\\
     &
     +\int_0^t\int_{\mathbb{R}}\left(\left[\left(\frac{\sqrt{\mu(v,\theta)\kappa(v,\theta)}}{v^3}\phi_{x}\right)_x\right]^2+
     \frac{\kappa(v,\theta)}{\mu(v,\theta) v^4}\left[\left(\frac{\mu(v,\theta)}{v}\phi_x\right)_x\right]^2\right)\,dxd\tau \\
      \le &C(\underline V,\overline V,\underline\Theta,\overline \Theta)\left\|\left(\phi_0,\psi_0,\frac{\zeta_0}{\sqrt{\delta}},\phi_{0x}\right)\right\|^2+C_9\left(N_1^6\delta^{\frac{1}{4}}
      +N_1^7\delta^{\frac{1}{4}}\int_0^t\|(\phi_x,\phi_{xx},\psi_x,\zeta_x)(\tau)\|^2\,d\tau\right).
  \end{split}
\end{equation}
 \end{Lemma}
\noindent{\bf Proof.}~~Rewriting $(\ref{3.1})_2$ as
\begin{equation}\label{3.68}
  \begin{split}
    \left(\frac{\mu(v,\theta)\phi_x}{v}\right)_t-\psi_t+\frac{R\theta}{v^2}\phi_x=& -\left(\frac{\mu(v,\theta){ V}_x}{v}\right)_t+\frac{R\zeta_x}{v}-\frac{R\zeta-p_+\phi}{v^2}{V}_x \\
      & +\frac{\mu_{\theta}(v,\theta)(v_x\theta_t-\theta_xu_x)}{v}-K_x+{U}_t.
  \end{split}
\end{equation}
Multiplying  (\ref{3.68}) by $\displaystyle\frac{\mu(v,\theta)\phi_x}{v}$ and integrating the resulting equation over $[0,t]\times \mathbb{R}$, we have
\begin{equation}\label{3.69}
  \begin{split}
     & \int_{\mathbb{R}}\frac{\mu^2(v,\theta)}{4v^2}\phi_x^2\,dx+\int_0^t\int_{\mathbb{R}}\frac{R\mu(v,\theta)\theta}{v^3}\phi_x^2\,dxd\tau \\
     \le &\int_{\mathbb{R}}\frac{\mu^2(v_0,\theta_0)}{2v_0^2}\phi_{0x}^2\,dx+\|\psi(t)\|^2+\|\psi_0\|^2+\int_0^t\int_{\mathbb{R}}\frac{\mu(v,\theta)\psi_x^2}{v}\,dxd\tau +\sum_{i=0}^{3}\int_0^t\int_{\mathbb{R}}J_i\,dxd\tau,
  \end{split}
\end{equation}
where
\begin{equation*}
  \begin{split}
    J_0&=-\psi\phi_x\theta_t\frac{\mu_{\theta}(v,\theta)}{v}+\frac{\mu\mu_{\theta}}{v^2}v_x\phi_x\theta_t-\frac{\mu\mu_{\theta}}{v^2}\theta_t\phi_x{V}_x , \\
      J_1&=-\frac{\mu\mu_{\theta}\theta_xu_x\phi_x}{v^2}-\frac{\mu_{\theta}\theta_x\psi\psi_x}{v},\\
      J_2&=\left(\frac{\mu_v}{v}-\frac{\mu}{v^2}\right)(-\psi\phi_x{ U}_x+\psi\psi_x{ V}_x)-\frac{\mu}{v^3}(R\zeta-p_+\phi){V}_x\phi_x\\
      &\quad-\left(\frac{\mu(v,\theta)}{v}\right)_v\frac{\mu}{v}\phi_x{ V}_x{ V}_t-\frac{\mu^2}{v^2}\phi_x{ V}_{xt}+\frac{\mu R}{v^2}\zeta_x\phi_x+\frac{\mu}{v}\phi_x{ U}_t,\\
      J_3&=-\left(\frac{\mu(v,\theta)}{v}\right)_v\frac{\mu}{v}\psi_x\phi_x{ V}_x,\quad J_4=-\frac{\mu}{v}\phi_xK_x.
  \end{split}
\end{equation*}
and we have used the the following Cauchy inequality:
\begin{equation*}
  \left|\int_{\mathbb{R}}\psi\frac{\mu\phi_x}{v}\,dx\right|\le \|\psi(t)\|^2+\frac{1}{4}\int_{\mathbb{R}}\frac{\mu^2(v,\theta)\phi_x^2}{v^2}\,dx.
\end{equation*}

Now we estimate the terms $\displaystyle\int_0^t\int_{\mathbb{R}}|J_i|\ dxd\tau, \,i=0,1,2,3$ one by one.
First, using $(\ref{1.1})_3$, we have
\begin{eqnarray}
    J_0 &=&\{\cdots\}_x+\left(-\psi\phi_x\frac{\mu_{\theta}}{v}+\frac{\mu\mu_{\theta}}{v^2}v_x\phi_x-\frac{\mu\mu_{\theta}}{v^2}\phi_x{ V}_x \right ) \frac{1}{(C_v-\frac{\theta}{2}\kappa_{\theta\theta}\frac{v_x^2}{v^5})}\left(-pu_x+\frac{\mu(v,\theta)u_x^2}{v}\right)\nonumber\\
      &&-\left(-\psi\phi_x\frac{\mu_{\theta}}{v}+\frac{\mu\mu_{\theta}}{v^2}v_x\phi_x-\frac{\mu\mu_{\theta}}{v^2}\phi_x{ V}_x \right)_x \frac{1}{(C_v-\frac{\theta}{2}\kappa_{\theta\theta}\frac{v_x^2}{v^5})}\left(\frac{\tilde{\alpha}(v,\theta)\theta_x}{\theta}+\frac{\theta\kappa_{\theta}v_x}{v^5}u_x\right)\nonumber\\
      &&+\left(\psi\phi_x\frac{\mu_{\theta}}{v}-\frac{\mu\mu_{\theta}}{v^2}v_x\phi_x+\frac{\mu\mu_{\theta}}{v^2}\phi_x{V}_x\right) \left[\left(\frac{1}{C_v-\frac{\theta}{2}\kappa_{\theta\theta}\frac{v_x^2}{v^5}}\right)_x\frac{\tilde{\alpha}\theta_x}{\theta}+\left(\frac{\theta\kappa_{\theta}v_x}{(C_v-\frac{\theta}{2}\kappa_{\theta\theta}\frac{v_x^2}{v^5})v^5}\right)_xu_x\right]
      \nonumber\\
&=&\{\cdots\}_x+O(1)\delta\left(|\psi\phi_x|+|{ V}_x\phi_x|+\phi_x^2)(|\psi_x|+|{ U}_x|+|\psi_x|^2\right)\nonumber\\
      &&+O(1)\delta\left(|\psi_x\phi_x|+|\psi\phi_{xx}|+|\psi\phi_x(\phi_x+{V}_x)|+|\psi\phi_x(\zeta_x+{\Theta}_x)|+|\phi_x^3|+|\phi_x^2{ V}_x|+|\zeta_x\phi_x^2|\right.\nonumber\\
      &&\left.+|\zeta_x\phi_x{V}_x|+|{\Theta}_x\phi_x^2|+|{ \Theta}_x{ V}_x\phi_x|+|\phi_x\phi_{xx}|+|{ V}_{xx}\phi_x|+|{V}_x\phi_{xx}|+|\phi_x{ V}_x^2|\right)\left(|\theta_x|+|v_xu_x|\right)\nonumber\\
    &&+O(1)\delta\left(|\psi\phi_x|+|\phi_x^2|+|\phi_x{ V}_x|\right)\left(|(v_x\theta_x,v_x^2,v_{xx})|+|(v_x^2,v_x^3,v_x^4,v_x^2v_{xx})|\right)|u_x|\nonumber\\
      &&+O(1)\delta\left(|\psi\phi_x|+|\phi_x^2|+|\phi_x{ V}_x|\right)\left(|\theta_x^2v_x^2|+|\theta_xv_x^3|+|\theta_xv_xv_{xx}|\right).\nonumber
\end{eqnarray}
Then similar to the estimate of $I_1$, we obtain

\begin{equation}\label{3.70}
  \int_0^t\int_{\mathbb{R}}|J_0|\,dxd\tau\leq CN_1^5\delta\int_0^t\|(\phi_{xx},\phi_x,\psi_x,\zeta_x)(\tau)\|^2\,d\tau
      +C\delta^2\int_0^t\int_{\mathbb{R}}\psi^2w^2\,dx
      +CN_1^5\delta^{\frac{3}{2}}.
\end{equation}
The terms $\displaystyle\int_0^t\int_{\mathbb{R}}|J_i|\,dxd\tau, i=1,2,3$ can be controlled by
\begin{equation}\label{3.71}
  \begin{split}
     \int_0^t\int_{\mathbb{R}}|J_1|\ dxd\tau\le & C\int_0^t\int_{\mathbb{R}}(|\theta_xu_x\phi_x|+|\theta_x\psi\psi_x|)\,dxd\tau \\
\le & \int_0^t\int_{\mathbb{R}}(|\zeta_x\psi_x\phi_x|+|\zeta_x\phi_x{ U}_x|+|\psi_x\phi_x{\Theta}_x|+|\phi_x{ \Theta}_x{ U}_x|+|\zeta_x\psi\psi_x|+|\psi\psi_x{\Theta}_x|)\,dxd\tau\\
\le& C\int_0^t\left\{\|\phi_x(\tau)\|^{\frac{1}{2}}\|\phi_{xx}(\tau)\|^{\frac{1}{2}}\sup_{0\le\tau\le t}\|\zeta_x(\tau)\|\|\psi_x(\tau)\|+\delta^{\frac{1}{2}}\|(\zeta_x,\psi_x,\phi_x)(\tau)\|^2\right.\\
      &\left.+\|\psi_{x}(\tau)\|^{\frac{1}{2}}\|\psi_{x}(\tau)\|^{\frac{3}{2}}\|\zeta_{x}(\tau)\|\right\}\,d\tau+\int_0^t\int_{\mathbb{R}}\frac{\mu(v,\theta){ \Theta}\psi_x^2}{v\theta}\,dxd\tau\\
      &+C\delta\int_0^t\int_{\mathbb{R}}\psi^2w^2\,dxd\tau+C\delta^2\int_0^t(1+\tau)^{-\frac{3}{2}}\,d\tau\\
\le & CN_1^{\frac{3}{2}}\delta^{\frac{1}{2}}\int_0^t\|(\phi_x,\phi_{xx},\psi_x,\zeta_x)(\tau)\|^2\,d\tau+\int_0^t\int_{\mathbb{R}}\frac{\mu(v,\theta){ \Theta}\psi_x^2}{v\theta}\,dxd\tau\\
      &+C\delta\int_0^t\int_{\mathbb{R}}\psi^2w^2\,dxd\tau+C\delta^2,
  \end{split}
\end{equation}

\begin{equation}\label{3.72}
  \begin{split}
    \int_0^t\int_{\mathbb{R}}|J_2|\ dxd\tau \le& C\int_0^t\int_{\mathbb{R}}(|\psi\phi_x{U}_x|+|\psi\psi_x{ V}_x|+|(\zeta,\phi){ V}_x\phi_x|+|\phi_x{ V}_x{ V}_t|+|\phi_x{ V}_{xt}| \\
      & +|\phi_x{ U}_t|)\,dxd\tau+\int_0^t\int_{\mathbb{R}}\left|\frac{\mu R}{v^2}\zeta_x\phi_x\right|\,dxd\tau\\
\le&\frac{1}{2}\int_0^t\int_{\mathbb{R}}\frac{R\mu\theta}{v^3}\phi_x^2\ dxd\tau+C(\underline{\Theta},\overline{ \Theta})\left\|\frac{\mu(v,\theta)}{\tilde{\alpha}(v,\theta)}\right\|_{L_{T,x}^{\infty}}\int_0^t\int_{\mathbb{R}}\frac{\tilde{\alpha}(v,\theta){ \Theta}\zeta_x^2}{v\theta^2}\,dxd\tau  \\
      &+\int_0^t\int_{\mathbb{R}}\frac{\mu(v,\theta){ \Theta}\psi_x^2}{v\theta}\,dxd\tau+C\delta\int_0^t\int_{\mathbb{R}}\psi^2w^2\,dxd\tau+C\delta,
  \end{split}
\end{equation}
and
\begin{equation}\label{3.73}
  \int_0^t\int_{\mathbb{R}}|J_3|\,dxd\tau\le \int_0^t\int_{\mathbb{R}}\frac{\mu(v,\theta){ \Theta}\psi_x^2}{v\theta}\,dxd\tau+C\delta^{\frac{1}{2}}\int_0^t\|(\phi_x,\psi_x)(\tau)\|^2\,d\tau.
\end{equation}
where in (\ref{3.72}),  we have used the assumption (\ref{1.16-1}).

For $\displaystyle\int_0^t\int_{\mathbb{R}}J_4\,dxd\tau$,  we have by a direct computation that
\begin{eqnarray}\label{3.73-1}
   -\frac{\mu}{v}\phi_xK_x
      &=&\{\cdots\}_x+\left(\frac{\mu}{v}\phi_x\right)_x\left(\frac{-\kappa v_{xx}}{v^5}+\frac{5\kappa-v\kappa_v}{2v^6}v_x^2-\frac{\kappa_{\theta}v_x\theta_x}{v^5}\right)\nonumber\\
      &=&\{\cdots\}_x-\frac{\mu\kappa}{v^6}\phi_{xx}^2+f(v,\theta)\phi_x^4+\frac{1}{3}\left(\frac{\kappa}{v^5}\left(\frac{\mu}{v}\right)_v\right)_v{ V}_x\phi_x^3+\frac{1}{3}\left(\frac{\kappa}{v^5}\left(\frac{\mu}{v}\right)_v\right)_{\theta}\theta_x\phi_x^3\nonumber\\
      &&-\frac{1}{6}\left(\frac{\mu}{v^7}(5\kappa-v\kappa_v)\right)_v{ V}_x\phi_x^3-\frac{1}{6}\left(\frac{\mu}{v^7}(5\kappa-v\kappa_v)\right)_{\theta}(\zeta_x+{ \Theta}_x)\phi_x^3\nonumber\\
      &&+\frac{\mu_{\theta}(5\kappa-v\kappa_v)}{2v^7}(\zeta_x+{ \Theta}_x)\phi_x^3-\left(\frac{\mu}{v}\right)_v(\phi_x+{ V}_x)\phi_x\frac{\kappa_{\theta}(\phi_x+{ V}_x)(\zeta_x+{ \Theta}_x)}{v^5}\nonumber\\
      &&+\left[\left(\frac{\mu}{v}\right)_v\left(\frac{3}{2}\phi_x+{ V}_x\right)+\frac{\mu_{\theta}}{v}(\zeta_x+{\Theta}_x)\right]\phi_x\frac{5\kappa-v\kappa_v}{2v^6}(2\phi_x{ V}_x+{V}_x^2)\nonumber\\
      && +\frac{\mu(5\kappa-v\kappa_v)}{2v^7}(2\phi_x{ V}_x+{ V}_x^2)\phi_{xx}-\frac{\kappa_{\theta}\mu\phi_{xx}(\phi_x+{ V}_x)(\zeta_x+{ \Theta}_x)}{v^6}\nonumber\\
      &&-\left(\frac{\mu}{v}\right)_v\frac{\kappa {V}_x\phi_{xx}\phi_x}{v^5}-\frac{\kappa\mu_{\theta}}{v^6}(\zeta_x+{ \Theta}_x)\phi_x\phi_{xx}-\left(\frac{\mu}{v}\right)_v(\phi_x+{V}_x)\phi_x\frac{\kappa { V}_{xx}}{v^5}\nonumber\\
      &&-\frac{\kappa {V}_{xx}\mu_{\theta}}{v^6}(\zeta_x+{ \Theta}_x)\phi_x-\frac{\mu_{\theta}\kappa_{\theta}}{v^6}(\zeta_x+{\Theta}_x)^2\phi_x(\phi_x+{ V}_x)-\frac{\mu\kappa}{v^6}\phi_{xx}{V}_{xx}\nonumber\\
      &=&\{\cdots\}_x-\frac{\mu\kappa}{v^6}\phi_{xx}^2+\left\{\displaystyle\frac{1}{3}\left(\frac{\kappa}{v^5}\left(\frac{\mu}{v}\right)_v\right)_v+\displaystyle
      \left(\frac{\mu}{v}\right)_v\frac{5\kappa-v\kappa_v}{2v^6}-\frac{1}{3}\left(\frac{\mu}{2v^7}(5\kappa-v\kappa_v)\right)_v\right\}\phi_x^4\nonumber\\
      &&+O(1)\left\{|{ \Theta}_x\phi_x^3|+|\zeta_x\phi_x^3|+(|{\Theta}_x|^2+|{ V}_{xx}|)\phi_x^2+(|{ \Theta}_x^3|+|{\Theta}_x{V}_{xx}|)|\phi_x|\right.\nonumber\\
      &&+({\Theta}_x^2+|{ V}_{xx}|)|\phi_x\zeta_x|+|{\Theta}_x\zeta_x\phi_x^2|+|{ V}_x\zeta_x^2\phi_x|+(|{ V}_x^2|+|{V}_{xx}|)|\phi_{xx}|\nonumber\\
      &&\left.+|\phi_{xx}\phi_x\zeta_x|+|{ \Theta}_x\phi_x\phi_{xx}|+|\phi_{xx}\zeta_x{V}_x|+\zeta_x^2\phi_x^2\right\},
\end{eqnarray}
and
\begin{eqnarray}\label{3.73-2}
  -\frac{\mu\kappa}{v^6}\phi_{xx}^2&=&-\left[\left(\frac{\sqrt{\mu\kappa}}{v^3}\phi_x\right)_x\right]^2+\left\{-\frac{2}{3}\left(\frac{\sqrt{\mu\kappa}}{v^3}
  \left(\frac{\sqrt{\mu\kappa}}{v^3}\right)_v\right)_v+\left[\left(\frac{\sqrt{\mu\kappa}}{v^3}\right)_v\right]^2\right\}\phi_x^4\nonumber\\
  &&+2\left[\left(\frac{\sqrt{\mu\kappa}}{v^3}\right)_\theta\right]^2\theta_x^2\phi_x^2
  +2\left(\frac{\sqrt{\mu\kappa}}{v^3}\right)_\theta\left(\frac{\sqrt{\mu\kappa}}{v^3}\right)_v v_x\theta_x\phi_x^2+2\frac{\sqrt{\mu\kappa}}{v^3}\left(\frac{\sqrt{\mu\kappa}}{v^3}\right)_\theta\theta_x\phi_x\phi_{xx}\nonumber\\
  &&+2\frac{\sqrt{\mu\kappa}}{v^3}\left(\frac{\sqrt{\mu\kappa}}{v^3}\right)_v V_x\phi_x\phi_{xx}+\left[\left(\frac{\sqrt{\mu\kappa}}{v^3}\right)_v\right]^2(2\phi_xV_x+V_x^2)\phi_x^2\nonumber\\
  &&-\frac{2}{3}\left(\frac{\sqrt{\mu\kappa}}{v^3}
  \left(\frac{\sqrt{\mu\kappa}}{v^3}\right)_v\right)_v V_x\phi_x^3-\frac{2}{3}\left(\frac{\sqrt{\mu\kappa}}{v^3}
  \left(\frac{\sqrt{\mu\kappa}}{v^3}\right)_v\right)_\theta \theta_x\phi_x^3+\{\cdots\}_x.
\end{eqnarray}
Substituting  (\ref{3.73-2}) into (\ref{3.73-1}) leads to
 \begin{eqnarray}\label{3.73-3}
   -\frac{\mu}{v}\phi_xK_x
      &=&\{\cdots\}_x-\left[\left(\frac{\sqrt{\mu\kappa}}{v^3}\phi_x\right)_x\right]^2+f(v,\theta)\phi_x^4+O(1)\left\{|{ (|\Theta}_x|+|\zeta_x|)|\phi_x^3|+|{\Theta}_x|^2\phi_x^2\right.\nonumber\\
      &&+|{ V}_{xx}|\phi_x^2+(|{ \Theta}_x^3|+|{\Theta}_x{V}_{xx}|)|\phi_x|+({\Theta}_x^2+|{ V}_{xx}|)|\phi_x\zeta_x|+|{\Theta}_x\zeta_x\phi_x^2|+|{ V}_x\zeta_x^2\phi_x|\nonumber\\
      &&\left.+(|{ V}_x^2|+|{V}_{xx}|)|\phi_{xx}|+|\phi_{xx}\phi_x\zeta_x|+|{\Theta}_x\phi_x\phi_{xx}|+|\phi_{xx}\zeta_x{V}_x|+\zeta_x^2\phi_x^2\right\},
\end{eqnarray}
where the function $f(v,\theta)$ is defined in (\ref{1.16-3}).

Thus by the assumption (\ref{1.16-3}) and some similar estimates as (\ref{3.71})-(\ref{3.73}), we obtain
\begin{equation}\label{3.74}
  \begin{split}
    \int_0^t\int_{\mathbb{R}}J_4\,dxd\tau\le-\int_0^t\int_{\mathbb{R}}\left[\left(\frac{\sqrt{\mu\kappa}}{v^3}\phi_{x}\right)_x\right]^2\,dxd\tau +CN_1^2\delta^{\frac{1}{4}}\int_0^t\|(\phi_x,\phi_{xx},\zeta_x)(\tau)\|^2\,d\tau+C\delta^{\frac{1}{4}}.
  \end{split}
\end{equation}

On the other hand,  we can also deal with $\displaystyle-\frac{\mu}{v}\phi_x K_x$ by
\begin{eqnarray}\label{3.74-1}
   \displaystyle-\frac{\mu}{v}\phi_x K_x
      &=&\{\cdots\}_x+\left(\frac{\mu}{v}\phi_x\right)_x\left(\frac{-\kappa v_{xx}}{v^5}+\frac{5\kappa-v\kappa_v}{2v^6}v_x^2-\frac{\kappa_{\theta}v_x\theta_x}{v^5}\right)\nonumber\\
      &=&\{\cdots\}_x+\left(\frac{\mu}{v}\phi_x\right)_x\left(\frac{-\kappa}{v^2}\left(\frac{\mu v_x}{v}\cdot\frac{1}{\mu v^2}\right)_x+\frac{5\kappa-v\kappa_v}{2v^6}v_x^2-\frac{\kappa_{\theta}v_x\theta_x}{v^5}\right)\nonumber\\
      &=&\{\cdots\}_x-\frac{\kappa}{\mu v^4}\left[\left(\frac{\mu}{v}\phi_x\right)_x\right]^2+\frac{g(v,\theta)}{2\mu v^6}v_x^2\left(\frac{\mu}{v}\phi_x\right)_x\nonumber\\
      &&-\frac{\kappa}{\mu v^4}\left(\frac{\mu}{v}V_x\right)_x\left(\frac{\mu}{v}\phi_x\right)_x+\left\{-\frac{\kappa\mu}{v^3}\left(\frac{1}{\mu v^2}\right)_\theta-\frac{\kappa_\theta}{v^5}\right\}v_x\theta_x\left(\frac{\mu}{v}\phi_x\right)_x
      \nonumber\\
       &=&\{\cdots\}_x-\frac{\kappa}{\mu v^4}\left[\left(\frac{\mu}{v}\phi_x\right)_x\right]^2+\frac{g(v,\theta)}{2\mu v^6}v_x^2\left(\frac{\mu}{v}\phi_x\right)_x+O(1)\left\{|(v_x,\theta_x)|^2|V_x\phi_x|\right. \nonumber\\
       &&\left.+|(v_x,\theta_x)V_x\phi_{xx}|+|(v_x,\theta_x)V_{xx}\phi_{x}|+|V_{xx}\phi_{xx}|+|v_x\theta_x(v_x,\theta_x)\phi_{x}|+|v_x\theta_x\phi_{xx}|\right\},\nonumber
\end{eqnarray}
where the function $g(v,\theta)$ is defined in  (\ref{1.16-4}). Thus if the condition (\ref{1.16-4}) hold, then
\begin{equation}\label{3.74-2}
  \begin{split}
    \int_0^t\int_{\mathbb{R}}J_4\,dxd\tau\le-\int_0^t\int_{\mathbb{R}}\frac{\kappa}{\mu v^4}\left[\left(\frac{\mu}{v}\phi_x\right)_x\right]^2\,dxd\tau +CN_1^2\delta^{\frac{1}{4}}\int_0^t\|(\phi_x,\phi_{xx},\zeta_x)(\tau)\|^2\,d\tau+C\delta^{\frac{1}{4}}.
  \end{split}
\end{equation}

(\ref{3.67}) thus follows from (\ref{3.69})-(\ref{3.73}), (\ref{3.74}) and Lemmas 3.1-3.2. This completes the proof of Lemma 3.3.

As a direct consequence of Lemmas 3.1-3.3,  we have
\begin{Corollary}
There exist a  constant $C_{10}>0$ depending only on $\underline {V},\overline { V},\underline {\Theta},\overline{ \Theta}, m_0, M_0$ and a constant $C_{11}>0$ depending only on $\underline {V},\overline { V},\underline {\Theta},\overline{ \Theta}$ such that if
\begin{equation}\label{3.75}
C_{10}N_{1}^{11}\delta^{\frac{1}{4}}<1,
\end{equation}
then it holds for all $t\in[0,T]$ that
\begin{eqnarray}\label{3.75-1}
    &&\int_{\mathbb{R}}\left[R{\Theta}\Phi\left(\frac{v}{{V}}\right)+\frac{\psi^2}{2}+\frac{R}{\delta}{ \Theta}\Phi\left(\frac{\theta}{{ \Theta}}\right)\right]dx+\int_{\mathbb{R}}\left(\frac{\kappa(v,\theta)v_x^2}{v^5}+\frac{\mu^2(v,\theta)}{v^2}\phi_x^2\right)\,dx  \nonumber\\
      && +\int_0^t\int_{\mathbb{R}}\left(\frac{\mu{\Theta}\psi_x^2}{v\theta}+\frac{\tilde{\alpha}{\Theta}\zeta_x^2}{v\theta^2}+\frac{\mu\theta}{v^3}\phi_x^2+\left[\left(\frac{\sqrt{\mu\kappa}}{v^3}\phi_{x}\right)_x\right]^2
      +\frac{\kappa}{\mu v^4}\left[\left(\frac{\mu}{v}\phi_x\right)_x\right]^2\right)\,dxd\tau\nonumber\\
&\le& C_{11}\left\|\left(\phi_0,\psi_0,\frac{\zeta_0}{\sqrt{\delta}},\phi_{0x}\right)\right\|^2.
\end{eqnarray}
\end{Corollary}
\noindent{\bf Proof.}~~Notice that the a priori assumption $m_0\leq v(t,x)\leq M_0$ and (\ref{3.12}) imply that
\begin{equation}\label{3.75-2}
 \int_0^t\int_{\mathbb{R}}\left(\frac{\mu(v,\theta){\Theta}\psi_x^2}{v\theta}+\frac{\tilde{\alpha}(v,\theta){\Theta}\zeta_x^2}{v\theta^2}+\frac{\mu(v,\theta)\theta}{v^3}\phi_x^2\right)\,d\tau
 \geq C_{12}\int_0^t\|(\phi_x,\psi_x,\zeta_x)(\tau)\|^2\,d\tau,
\end{equation}
where $C_{12}$ is a positive constant depending only $\underline {V},\overline { V},\underline {\Theta},\overline{ \Theta}, m_0, M_0$. Without loss of generality, we can assume $C_{12}<1$.

Thus by adding (\ref{3.13}), (\ref{3.33}) and (\ref{3.67}) together,  and choosing $\delta>0$ sufficiently small such that
 \begin{equation}\label{3.75-3}
C_{13}N_1^7\delta^{\frac{1}{4}}\leq\frac{C_{12}}{2},\quad C_{13}:=\max\{2C_5C_7+C_9\},
\end{equation}
 we have
\begin{equation}\label{3.75-4}
A(t)\leq C(\underline {V},\overline {V},\underline {\Theta},\overline{\Theta} )\left\|\left(\phi_0,\psi_0,\frac{\zeta_0}{\sqrt{\delta}},\phi_{0x}\right)\right\|^2+2C_{13}N_1^7\delta^{\frac{1}{4}}\int_0^t\|\phi_{xx}(\tau)\|^2\,d\tau,
\end{equation}
where $A(t)$ denote the formula on the left hand side of (\ref{3.75-1}).

To estimate the reminder term $\int_0^t\|\phi_{xx}(\tau)\|^2d\tau$ in (\ref{3.75-4}), we rewrite
\begin{eqnarray}\label{3.75-5}
\phi_{xx}&=&\frac{v^6}{\mu\kappa}\left\{\left[\left(\frac{\sqrt{\mu\kappa}}{v^3}\phi_x\right)_x\right]^2
-\left[\left(\frac{\sqrt{\mu\kappa}}{v^3}\right)_v\right]^2v_x^2\phi_x^2-\left[\left(\frac{\sqrt{\mu\kappa}}{v^3}\right)_\theta\right]^2\theta_x^2\phi_x^2
\right.\nonumber\\
&&\left.-2\left(\frac{\sqrt{\mu\kappa}}{v^3}\right)_v\left(\frac{\sqrt{\mu\kappa}}{v^3}\right)_\theta v_x\theta_x\phi_x^2-2\frac{\sqrt{\mu\kappa}}{v^3}\left(\frac{\sqrt{\mu\kappa}}{v^3}\right)_v v_x\phi_x\phi_{xx}-2\frac{\sqrt{\mu\kappa}}{v^3}\left(\frac{\sqrt{\mu\kappa}}{v^3}\right)_\theta \phi_x\phi_{xx}\theta_x\right\}\nonumber\\
&=&O(1)\left\{\left[\left(\frac{\sqrt{\mu\kappa}}{v^3}\phi_x\right)_x\right]^2+\phi_x^4+\Theta_x^2\phi_x^2+\zeta_x^2\phi_x^2+|\phi_x^2\phi_{xx}|+|\Theta_x\phi_x\phi_{xx}|+|\zeta_x\phi_x\phi_{xx}|\right\}.
\end{eqnarray}
Consequently, it follows from (\ref{3.75-5}), the Cauchy inequality and the Sobolev inequality that
\begin{equation}\label{3.75-6}
 \aligned
\int_0^t\|\phi_{xx}(\tau)\|^2\,d\tau\leq& C\int_0^t\int_{\mathbb{R}}\left[\left(\frac{\sqrt{\mu\kappa}}{v^3}\phi_x\right)_x\right]^2dxd\tau+\frac{1}{4}\int_0^t\|\phi_{xx}(\tau)\|^2\,d\tau\\
&+C\int_0^t\left(\left\|\Theta_x^2(\tau)\right\|_{L^\infty}\|\phi_x(\tau)\|^2+\left(\|\zeta_x(\tau)\|^2\|\phi_x(\tau)\|+\|\phi_x(\tau)\|^3\right)\|\phi_{xx}(\tau)\|\right)d\tau\\
\leq& C\int_0^t\int_{\mathbb{R}}\left[\left(\frac{\sqrt{\mu\kappa}}{v^3}\phi_x\right)_x\right]^2\,dxd\tau+\frac{1}{2}\int_0^t\|\phi_{xx}(\tau)\|^2\,d\tau+
CN_1^4\int_0^t\|\phi_x(\tau)\|^2\,d\tau,
\endaligned
\end{equation}
which implies that
\begin{equation}\label{3.75-7}
\int_0^t\|\phi_{xx}(\tau)\|^2\,d\tau
\leq C_{14}\int_0^t\int_{\mathbb{R}}\left[\left(\frac{\sqrt{\mu\kappa}}{v^3}\phi_x\right)_x\right]^2\,dxd\tau+
C_{14}N_1^4\int_0^t\|\phi_x(\tau)\|^2\,d\tau,
\end{equation}
where $C_{14}$ is a positive constant depending only on $\underline {V},\overline { V},\underline {\Theta},\overline{ \Theta}, m_0, M_0$.

Inserting (\ref{3.75-7}) into (\ref{3.75-4}), then  (\ref{3.75-1}) holds    provided that $\delta$ is sufficiently small such that
\begin{equation}\label{3.75-8}
2C_{13}C_{14}N_1^{11}\delta^{\frac{1}{4}}\leq\frac{C_{12}}{4}.
\end{equation}
Letting $C_{10}=\max\{8C_{12}^{-1}C_{13}C_{14},2C_6\left(\min\{2p_+^2,R^2,\underline{\Theta}\}\right)^{-1}\}$, then  we  finish  the proof of Corollary 3.1.

Based on Corollary 3.1, we now show the uniform lower and upper bounds on $v(x,t)$ by using Y. Kanel's method \cite{Y. Kanel}.
\begin{Lemma}
 Under the assumptions of Corollary 3.1, there exist a positive constant $C_{15}$ depending only on $\underline V,\overline V,\underline \Theta,\overline \Theta, m_0, M_0$ and a constant $C_{0}>0$ depending only on $\underline V,\overline V,\underline \Theta,\overline \Theta$ and $\|(\phi_0,\psi_0,\frac{\zeta_0}{\sqrt{\delta}},\phi_{0x})\|$ such if
  \begin{equation}\label{3.75-8}
  C_{15}N_{1}^{11}\delta^{\frac{1}{4}}<1,
\end{equation}
 then it holds
 \begin{equation}\label{3.76}
  C_0^{-1}\le v(t,x)\le C_{0},\quad \forall \,(t,x)\in[0,T]\times{\mathbb{R}}.
\end{equation}
 \end{Lemma}
\noindent{\bf Proof.}~~First,  (\ref{3.75-1}) imply that
\begin{equation}\label{3.77}
  \int_{\mathbb{R}}\frac{\mu_1^2(v)}{v^2}\phi_x^2\,dx+\int_{\mathbb{R}}\frac{\kappa_1^2(v)}{v^5}v_x^2\,dx\le C_{11}\left\|\left(\phi_0,\psi_0,\frac{\zeta_0}{\sqrt{\delta}},\phi_{0x}\right)\right\|^2,
\end{equation}
 where $\mu_1(v):=\min_{\theta\in[\frac{\underline {\Theta}}{2},2\overline { \Theta}]}\{\mu(v,\theta)\}$ and $\kappa_1(v)=\min_{\theta\in[\frac{\underline { \Theta}}{2},2\overline { \Theta}]}\{\kappa(v,\theta)\}$.

Let $\displaystyle\tilde{v}=\frac{v}{{ V}}$,  then under the assumption (\ref{1.15}) and (\ref{1.16}), there exists a positive constant $C(\underline { V},\overline {V},\underline {\Theta},\overline {\Theta})$ such that
\begin{equation}\label{3.79}
    \mu_1(v)\ge C(\underline { V},\overline { V},\underline {\Theta},\overline {\Theta})\mu_1(\tilde{v}) ,\quad
    \kappa_1(v) \ge C(\underline {V},\overline {V},\underline {\Theta},\overline {\Theta})\kappa_1(\tilde{v}).
\end{equation}
Consequently, it follows from (\ref{3.77}) that
\begin{equation}\label{3.80}
  \int_{\mathbb{R}}\frac{\mu_1^2(\tilde{v})}{v^2}\phi_x^2\,dx+\int_{\mathbb{R}}\frac{\kappa_1^2(\tilde{v})}{v^5}v_x^2\,dx\le C(\underline { V},\overline { V},\underline { \Theta},\overline {\Theta})\left\|\left(\phi_0,\psi_0,\frac{\zeta_0}{\sqrt{\delta}},\phi_{0x}\right)\right\|^2.
\end{equation}

Set
\begin{equation}\label{3.82}
  \overline \Phi(\tilde{v}(t,x))=\int_{1}^{\tilde{v}}\frac{\sqrt{\Phi(\eta)}}{\eta}\mu_1(\eta)\,d\eta, \quad \Phi(\eta)=\eta-1-\ln{\eta},
\end{equation}
then by the assumption (\ref{1.15}),
\begin{equation}\label{3.83}
  \left|\overline \Phi(\tilde{v}(t,x))\right|\ge \left\{
  \begin{array}{l}
    \displaystyle A_1|\ln{\tilde{v}}|-A_2,\quad \tilde{v}\longrightarrow 0^+ ,\\[2mm]
   \displaystyle A_1 \tilde{v}^{\frac{1}{2}-b}-A_2,\quad \tilde{v}\longrightarrow +\infty,
  \end{array}\right.
\end{equation}
where $A_1>0,A_2>0$ are positive constants.
On the other hand, it holds
\begin{equation}\label{3.84}
  \begin{split}
    \left|\overline \Phi(\tilde{v}(t,x))\right| & =\left|\int_{-\infty}^{x}\overline \Phi(\tilde{v}(t,y))_y\,dy\right| \\
      & \le \int_{\mathbb{R}}\left|\frac{\sqrt{\Phi(\tilde{v})}}{\tilde{v}}\mu_1(\tilde{v})\tilde{v}_x\right|\,dx\\
      & \le\|\sqrt{\Phi(\tilde{v})}(t)\|\left\|\frac{\mu_1(\tilde{v})\tilde{v}_x}{\tilde{v}}(t)\right\|\\
      & \le C(\underline {V},\overline {V},\underline {\Theta},\overline {\Theta})\left\|\left(\phi_0,\psi_0,\frac{\zeta_0}{\sqrt{\delta}},\phi_{0x}\right)\right\|^2,
  \end{split}
\end{equation}
where we have used (\ref{3.75-1}) and  the fact that
\begin{equation}\label{3.85}
\aligned
  \left\|\frac{\mu_1(\tilde{v})\tilde{v}_x}{\tilde{v}}(t)\right\|&=\left\|\left(\frac{\mu_1(\tilde{v})\phi_x}{v}+\mu_1(\tilde{v})\left(\frac{{ V}_x}{v}-\frac{{ V}_x}{{ V}}\right)\right)(t)\right\|\\
  &\le \left\|\frac{\mu_1(\tilde{v})\phi_x}{v}(t)\right\|+C_{16}(\underline {V},\overline {V}, m_0,M_0)\delta^{\frac{3}{4}}\\
  &\le  C(\underline {V},\overline { V},\underline { \Theta},\overline {\Theta})\left\|\left(\phi_0,\psi_0,\frac{\zeta_0}{\sqrt{\delta}},\phi_{0x}\right)\right\|
  \endaligned
\end{equation}
due to (\ref{3.80}) and the smallness of $\delta$ such that $C_{16}(\underline {V},\overline {V}, m_0,M_0)\delta^{\frac{3}{4}}<1$.

Then (\ref{3.83}) and (\ref{3.84}) lead to
\begin{equation}\label{3.86}
  C_{17}\exp\left\{-C_{18}\left\|\left(\phi_0,\psi_0,\frac{\zeta_0}{\sqrt{\delta}},\phi_{0x}\right)\right\|^2\right\}\le v(t,x)\le C_{19}\left\|\left(\phi_0,\psi_0,\frac{\zeta_0}{\sqrt{\delta}},\phi_{0x}\right)\right\|^{\frac{4}{1-2b}}
\end{equation}
for all $(t,x)\in[0,T]\times{\mathbb{R}}$, where $C_{17}=C_{17}(\underline {V})>0$ and $C_{18}, C_{19}$ are positive constants depending only on $\underline { V},\overline {V},\underline {\Theta},\overline {\Theta}$.

Now we suppose the condition (\ref{1.16}) holds. Since
\begin{equation*}
  \frac{\sqrt{\kappa_1(\tilde{v})}\tilde{v}_x}{\tilde{v}^{\frac{5}{2}}}=\frac{\sqrt{\kappa_1(\tilde{v})}v_x{ V}^{\frac{3}{2}}}{v^{\frac{5}{2}}}-\frac{\sqrt{\kappa_1(\tilde{v})}{ V}^{\frac{1}{2}}{V}_x}{v^{\frac{3}{2}}},
\end{equation*}
we have from (\ref{3.80}) that
\begin{equation}\label{3.87}
  \begin{split}
    \left\|\frac{\sqrt{\kappa_1(\tilde{v})}\tilde{v}_x}{\tilde{v}^{\frac{5}{2}}}(t)\right\| & \le C(\underline { V},\overline { V})\left\|\frac{\sqrt{\kappa_1(\tilde{v})}\tilde{v}_x}{\tilde{v}^{\frac{5}{2}}}(t)\right\|+C(\underline {V},\overline {V}, m_0,M_0)\|{ V}_x(t)\| \\
      &\le C(\underline { V},\overline {V})\left\|\frac{\sqrt{\kappa_1(\tilde{v})}\tilde{v}_x}{\tilde{v}^{\frac{5}{2}}}(t)\right\|+C_{20}(\underline {V},\overline {V}, m_0,M_0)\delta^{\frac{3}{4}}\\
      &\le C(\underline { V},\overline { V},\underline { \Theta},\overline {\Theta})\left\|\left(\phi_0,\psi_0,\frac{\zeta_0}{\sqrt{\delta}},\phi_{0x}\right)\right\|,
  \end{split}
\end{equation}
provided that $\delta$ is sufficiently small such that $C_{20}(\underline {V},\overline {V}, m_0,M_0)\delta^{\frac{3}{4}}<1$.

Define
\begin{equation*}
   \Psi(\tilde{v}(t,x))=\int_{1}^{\tilde{v}}\frac{\sqrt{\Phi(\eta)}}{\eta^{\frac{5}{2}}}\sqrt{\kappa_1(\eta)}\ d\eta,
\end{equation*}
then  the assumption  (\ref{1.16}) implies
\begin{equation}\label{3.88}
  | \Psi(\tilde{v}(t,x))|\ge \left\{
  \begin{array}{l}
    A_3|\ln{\tilde{v}}|-A_4,\quad \tilde{v}\longrightarrow 0^+ ,\\[2mm]
    A_3 \tilde{v}^{-1-\frac{d}{2}}-A_4,\quad \tilde{v}\longrightarrow +\infty,
  \end{array}\right.
\end{equation}
where $A_1>0,A_2>0$ are positive constants.
On the other hand, it follows from (\ref{3.75-1}), (\ref{3.80}) and  (\ref{3.87}) that
\begin{equation}\label{3.89}
  \begin{split}
    \left|\Psi(\tilde{v}(t,x))\right| & =\left|\int_{-\infty}^{x} \Psi(\tilde{v}(t,x))_y\,dy\right| \\
      & \le \int_{\mathbb{R}}\left|\frac{\sqrt{\Phi(\tilde{v})}}{\tilde{v}^{\frac{5}{2}}}\sqrt{\kappa_1(\tilde{v})}\tilde{v}_x\right|\,dx\\
      & \le\|\sqrt{\Phi(\tilde{v})}(t)\|\left\|\frac{\sqrt{\kappa_1(\tilde{v})}\tilde{v}_x}{\tilde{v}^{\frac{5}{2}}}(t)\right\|\\
      & \le C(\underline V,\overline V,\underline \Theta,\overline \Theta)\left\|\left(\phi_0,\psi_0,\frac{\zeta_0}{\sqrt{\delta}},\phi_{0x}\right)\right\|^2.
  \end{split}
\end{equation}
(\ref{3.88}) together with (\ref{3.89}) implies
\begin{equation}\label{3.90}
  C_{21}\exp\left\{-C_{22}\left\|\left(\phi_0,\psi_0,\frac{\zeta_0}{\sqrt{\delta}},\phi_{0x}\right)\right\|^2\right\}\le v(t,x)\le C_{23}\left\|\left(\phi_0,\psi_0,\frac{\zeta_0}{\sqrt{\delta}},\phi_{0x}\right)\right\|^{\frac{-1}{2+d}},
\end{equation}
where $C_{21}=C_{21}(\underline {V})>0$ and $C_{22}, C_{23}$  are positive constants depending only on $\underline {V},\overline {V},\underline {\Theta},\overline{ \Theta}$.  Letting
 $$C_1=\max\left\{C_{17}\exp\left\{C_{18}N_{01}^2\right\},\, C_{19}N_{01}^{\frac{4}{1-2b}},\, C_{21}\exp\left\{C_{22}N_{01}^2\right\},\, C_{23}N_{01}^{\frac{-1}{2+d}}\right\},$$
where $N_{01}:=\|(\phi_0,\psi_0,\frac{\zeta_0}{\sqrt{\delta}},\phi_{0x})(t)\|$, and $C_{15}=\max\left\{C_{10},C_{16},C_{20}\right\}$,  then we can get (\ref{3.76}).  This completes the proof of Lemma 3.4.

 Lemmas 3.1-3.4 imply the following corollary.
\begin{Corollary}
Under the assumptions of Lemma 3.4, there exists a positive constant $C_{24}$ depending only on $\underline {V}, \overline {V}, \underline {\Theta}, \overline{ \Theta}$ and $\|(\phi_0,\psi_0,\frac{\zeta_0}{\sqrt{\delta}},\phi_{0x})\|$ such that
\begin{equation}\label{3.91}
  \begin{split}
     & \left\|\left(\phi,\psi,\frac{\zeta}{\sqrt{\delta}},\phi_{x}\right)(t)\right\|^2+\int_0^t\|(\phi_x,\phi_{xx},\psi_x,\zeta_{x})(\tau)\| ^2\,d\tau \\
      \le & C_{24}\left\|\left(\phi_0,\psi_0,\frac{\zeta_0}{\sqrt{\delta}},\phi_{0x}\right)\right\|^2.
  \end{split}
\end{equation}
\end{Corollary}
\noindent{\bf Proof.}~~First, it is easy to see from Corollary 3.1 and Lemma 3.4 that
\begin{equation}\label{3.91-1}
  \begin{split}
     & \left\|\left(\phi,\psi,\frac{\zeta}{\sqrt{\delta}},\phi_x\right)(t)\right\|^2+\int_0^t\|(\psi_x,\zeta_x)(\tau)\| ^2\,d\tau \\
      \le & C_{25}\left\|\left(\phi_0,\psi_0,\frac{\zeta_0}{\sqrt{\delta}},\phi_{0x}\right)\right\|^2,
  \end{split}
\end{equation}
where $C_{25}$ is a positive constant depending only on  $\underline {V}, \overline {V}, \underline {\Theta}, \overline{ \Theta}$ and $\|(\phi_0,\psi_0,\frac{\zeta_0}{\sqrt{\delta}},\phi_{0x})\|$.

 On the other hand,
it follows from   (\ref{3.75-1}), (\ref{3.91-1}),  the Cauchy inequality, the Young inequality and the Sobolev inequality that
\begin{eqnarray}\label{3.91-2}
     &&\left\|\phi_x(t)\right\|^2+\int_0^t\|\phi_{x}(\tau)\|_1^2\,d\tau \nonumber\\
      &\le& C_{26}\left\|(\phi_0,\psi_0,\frac{\zeta_0}{\sqrt{\delta}},\phi_{0x})\right\|^2+ C_{26}\int_0^t\int_{\mathbb{R}}\left(\phi_x^4+|\phi_x^3\Theta_x|+|\phi_x^3\zeta_x|
      +|\phi_x^2\Theta_x^2|\right.\nonumber\\
      &&\left.+|\phi_x^2\zeta_x^2|+|\phi_x^2\zeta_x\Theta_x|+|\phi_x\zeta_x\phi_{xx}|+|\Theta_x\phi_x\phi_{xx}|\right)\,dxd\tau\nonumber\\
      &\le& C_{26}\left\|(\phi_0,\psi_0,\frac{\zeta_0}{\sqrt{\delta}},\phi_{0x})\right\|^2+\frac{1}{4}\int_0^t\|\phi_{xx}(\tau)\| ^2\,d\tau+C_\eta\int_0^t\left(\|\phi_x(\tau)\|^3\|\phi_{xx}(\tau)\|\right.\nonumber\\
      &&\quad\left.+\|\phi_x(\tau)\|^2\|\phi_{xx}(\tau)\|\|\zeta_{x}(\tau)\|+\|\phi_x(\tau)\|\|\phi_{xx}(\tau)\|\|\zeta_{x}(\tau)\|^2+\|(\phi_x,\zeta_x)(\tau)\|^2\right)d\tau.\nonumber\\
       &\le& C_{26}\left\|(\phi_0,\psi_0,\frac{\zeta_0}{\sqrt{\delta}},\phi_{0x})\right\|^2+\frac{1}{2}\int_0^t\|\phi_{xx}(\tau)\| ^2\,d\tau+C\int_0^t\left(\sup_{\tau\in[0,T]}\|\phi_x(\tau)\|^4+1\right)\|(\phi_x,\zeta_x)(\tau)\|^2\,d\tau\nonumber\\
       &&\quad+C\int_0^t\sup_{\tau\in[0,T]}\{\|\phi_x(\tau)\|^2\|\zeta_x(\tau)\|^2\}
       \|\zeta_x(\tau)\|^2\,d\tau  \nonumber\\
       &\le& C_{26}\left\|(\phi_0,\psi_0,\frac{\zeta_0}{\sqrt{\delta}},\phi_{0x})\right\|^2+\frac{1}{2}\int_0^t\|\phi_{xx}(\tau)\| ^2 d\tau+C(N_{01}^4+N_{01}^2N_1^2\delta)\int_0^t\|(\phi_x,\zeta_x)(\tau)\|^2\,d\tau\nonumber\\
       &\le& C_{26}\left\|(\phi_0,\psi_0,\frac{\zeta_0}{\sqrt{\delta}},\phi_{0x})\right\|^2+\frac{1}{2}\int_0^t\|\phi_{xx}(\tau)\| ^2\,d\tau+C N_{01}^6,
\end{eqnarray}
which implies that
\begin{equation}\label{3.91-3}
 \left\|\phi_x(t)\right\|^2+\int_0^t\|\phi_x(\tau)\|_1^2\,d\tau
      \le  C_{27}N_{01}^6.
\end{equation}
Here $C_{26}, C_{27}$ are positive constants depending only on  $\underline {V}, \overline {V}, \underline {\Theta}, \overline{ \Theta}$ and $N_{01}$,  and in the last step of (\ref{3.91-2}), we have used the smallness of $\delta$ such that $N_1^2\delta<1$.

Letting $C_{24}=\max\{C_{25}, C_{27}N_{01}^4\}$, then  (\ref{3.91}) follows  from (\ref{3.91-1}) and (\ref{3.91-3}) immediately. This completes the proof of Corollary 3.2.

The following lemma give the estimate on $\|(\psi_x,\zeta_x/\sqrt{\delta})(t)\|$.
\begin{Lemma}
 Under the assumptions of Lemma 3.4, there exists a positive constant $C_{28}$ depending only on $\underline V,\overline V,\underline \Theta,\overline \Theta$ and $N_{01}$ and a positive constant $C_{29}$ depending only on $\underline V,\overline V,\underline \Theta,\overline \Theta$ and $N_{0}$ such that if
 \begin{equation}\label{3.92}
 C_{28}N_1^4\delta^{\frac{1}{4}}<\frac{1}{3},\quad C_{28}N_1^4\varepsilon<\frac{1}{3},\quad N_2^2\delta^{\frac{3}{4}}<1,
\end{equation}
then it holds  for all $t\in[0,T]$ that
 \begin{equation}\label{3.922}
  \left\|\left(\phi_{xx},\psi_x, \frac{\zeta_x}{\sqrt{\delta}}\right)(t)\right\|^2+\int_0^t\|(\psi_{xx},\zeta_{xx})(\tau)\|^2\,d\tau\le C_{29}N_{0}^2.
\end{equation}
  \end{Lemma}

\noindent{\bf Proof.}~~Multiplying $(\ref{3.1})_2$ by $-\psi_{xx}$ and using $(\ref{3.1})_1$, we have
\begin{eqnarray}\label{3.93}
     && \left(\frac{\psi_x^2}{2}\right)_t+\left(\frac{\kappa(v,\theta)\phi_{xx}^2}{2v^5}\right)_t+\frac{\mu(v,\theta)}{v}\psi_{xx}^2\nonumber \\
      &=&\left\{-\left(\frac{\mu(v,\theta)}{v}\right)_vv_x\psi_x\psi_{xx}-\left(\frac{\mu(v,\theta)}{v}\right)_{\theta}\theta_x\psi_x\psi_{xx}-\left[\left(\frac{\mu(v,\theta)}{v}-\frac{\mu({ V},{ \Theta})}{{ V}}\right){ U}_x\right]_x\psi_{xx}\right. \\
       &&\left. +\frac{R\zeta_x-p_+\phi_x}{v}\psi_{xx}+\frac{R\zeta-p_+\phi}{v^2}v_x\psi_{xx}+R_1\psi_{xx}\right\}+\left\{\frac{\kappa_v-5\kappa}{2v^6}u_x\phi_{xx}^2-\left(\frac{5\kappa-v\kappa_v}{2v^6}\right)_vv_x^3\psi_{xx}\right. \nonumber \\
       &&\left. -\left(\frac{5\kappa-v\kappa_v}{2v^6}\right)_{\theta}\theta_xv_x^2\psi_{xx} -\frac{5\kappa-v\kappa_v}{v^6}v_xv_{xx}\psi_{xx}+\left(\frac{\kappa(v,\theta)}{v^5}\right)_vv_x^2\theta_x\psi_{xx}+\left(\frac{\kappa_{\theta}(v,\theta)}{v^5}\right)_{\theta}v_x\theta_x^2\psi_{xx}\right\}\nonumber \\
       &&+\frac{\kappa_{\theta}\theta_t}{2v^5}\phi_{xx}^2+\left \{\frac{\kappa_{\theta}(v,\theta)}{v^5}v_{xx}\theta_x\psi_{xx}+\frac{\kappa_{\theta}(v,\theta)}{v^5}v_x\theta_{xx}\psi_{xx}\right\}\nonumber \\
       &=& J_4+J_5+J_6+J_7.\nonumber
\end{eqnarray}
Integrating (\ref{3.93}) over  $[0,t]\times \mathbb{R}$ gives
\begin{equation}\label{3.99}
  \begin{split}
     & \frac{1}{2}\int_{\mathbb{R}}\left(\psi_x^2+\frac{\kappa(v,\theta)\phi_{xx}^2}{v^5}\right)\ dx+\int_0^t\int_{\mathbb{R}}\frac{\mu(v,\theta)}{v}\psi_{xx}^2\,dxd\tau \\
      =& \frac{1}{2}\int_{\mathbb{R}}\left(\psi_{0x}^2+\frac{\kappa(v_0,\theta_0)\phi_{0xx}^2}{v_0^5}\right)\ dx+\int_0^t\int_{\mathbb{R}}\sum_{j=4}^{7}J_i\,dxd\tau.
  \end{split}
\end{equation}

We derive from the Cauchy inequality, the Sobolev inequality and Lemmas 2.1 and 3.4 that
\begin{eqnarray}\label{3.100}
    \left|\int_0^t\int_{\mathbb{R}}J_4\ dxd\tau\right| &\le& \eta \int_0^t\|\psi_{xx}(\tau)\|^2\,d\tau+C_{\eta}\int_0^t\int_{\mathbb{R}}\left(v_x^2\psi_x^2+\theta_x^2\psi_x^2+|(\phi_x,\zeta_x)|^2{ U}_x^2+|(\phi,\zeta){ \Theta}_x|^2{U}_x^2\right.
     \nonumber\\
      && \left.+|(\phi,\zeta){ U}_{xx}|^2+\zeta_x^2+\phi_x^2+|(\phi,\zeta)|^2v_x^2+R_1^2\right)dxd\tau \nonumber\\
&\le& \eta \int_0^t\|\psi_{xx}(\tau)\|^2\,d\tau+C_{\eta}\int_0^t\left\{\|\psi_x(\tau)\|\|\psi_{xx}(\tau)\|\sup_{0\le\tau\le t}\|\phi_x(\tau)\|^2+\|(\phi_x,\psi_x,\zeta_x)(\tau)\|^2\right.\nonumber\\
      &&\left.+\|\psi_x(\tau)\|\|\psi_{xx}(\tau)\|\sup_{0\le \tau\le t}\|\zeta_x(\tau)\|^2+\|(\phi,\zeta)\|_{L^{\infty}}^2(\|\phi_x(\tau))\|^2+\|{ U}_{xx}(\tau)\|^2)\right\}d\tau\nonumber\\
      &&+C_{\eta}\delta\int_0^t\int_{\mathbb{R}}|(\phi,\zeta)|^2w^2\,dxd\tau+C_{\eta}\delta^{\frac{3}{2}}\int_0^t(1+\tau)^{-\frac{5}{2}}\,d\tau\nonumber\\
&\le& 2\eta \int_0^t\|\psi_{xx}(\tau)\|^2\,d\tau+C_{\eta}N_{01}^4\int_0^t\|(\phi_x,\psi_x,\zeta_x)(\tau)\|^2\,d\tau\nonumber\\
      &&+C_{\eta}\delta\int_0^t\int_{\mathbb{R}}|(\phi,\zeta)|^2w^2\,dxd\tau+C_{\eta}N_{01}^2\delta^{\frac{3}{2}},
\end{eqnarray}
\begin{eqnarray}\label{3.101}
    \left|\int_0^t\int_{\mathbb{R}}J_5\,dxd\tau\right| &\le&C \int_0^t\int_{\mathbb{R}}\left\{|\psi_x\phi_{xx}^2|+|{ U}_x\phi_{xx}^2|+|v_x^3\psi_{xx}|+|v_xv_{xx}\psi_{xx}|+|v_x^2\theta_x\psi_{xx}|+|v_x\theta_x^3\psi_{xx}|\right\}dxd\tau\nonumber\\
      &\le& C\int_0^t\left\{\|\psi_x(\tau)\|^{\frac{1}{2}}\|\psi_{xx}(\tau)\|^{\frac{1}{2}}\|\phi_{xx}(\tau)\|^2+\|{ U}_x(\tau)\|_{L^{\infty}}\|\phi_{xx}(\tau)\|^2\right\}\,d\tau\nonumber\\
      &&+\eta\int_0^t\|\psi_{xx}(\tau)\|^2\,d\tau+C_{\eta}\int_0^t\int_{\mathbb{R}}\left(\phi_x^6+{V}_x^6+\phi_x^2\phi_{xx}^2+\phi_x^2{ V}_{xx}^2+{ V}_x^2\phi_{xx}^2+{V}_x^2{ V}_{xx}^2\right.\nonumber\\
      &&\left.+\zeta_x^2\phi_x^4+\zeta_x^2{V}_x^4+{ \Theta}_x^2\phi_x^4+{\Theta}_x^2{ V}_x^4+\zeta_x^4\phi_x^2+\zeta_x^4{V}_x^2+{\Theta}_x^4\phi_x^2+{ \Theta}_x^4{ V}_x^2\right)dxd\tau\nonumber\\
&\le& 2\eta\int_0^t\|\psi_{xx}(\tau)\|^2\,d\tau+C_{\eta}\int_0^t\left(\|(\psi_x,\phi_{xx})(\tau)\|^2+\|\phi_{xx}(\tau)\|^4+\|\phi_{x}(\tau)\|^4\|\phi_{xx}(\tau)\|^2\right)\,d\tau\nonumber\\
      &&+C_{\eta}\int_0^t\left\{\|\phi_{x}(\tau)\|\|\phi_{xx}(\tau)\|^3+\|\phi_{x}(\tau)\|^2\|\phi_{xx}(\tau)\|^2\|\zeta_{x}(\tau)\|^2++\|\phi_{x}(\tau)\|^3\|\phi_{xx}(\tau)\|\right.
      \nonumber\\&&\left.+\|\zeta_{x}(\tau)\|^3\|\zeta_{xx}(\tau)\|+\|\zeta_{x}(\tau)\|^2\|\zeta_{xx}(\tau)\|^2\|\phi_{x}(\tau)\|^2\right\}\,d\tau +C_{\eta}\delta^{\frac{7}{2}}\nonumber\\
&\le&2\eta\int_0^t\|\psi_{xx}(\tau)\|^2\,d\tau+C_{\eta}N_{01}^4\int_0^t(\|(\phi_x,\psi_x,\zeta_x,\phi_{xx})(\tau)\|^2\,d\tau\nonumber\\
      &&+C_{\eta}\int_0^t\left(\|\phi_{xx}(\tau)\|^4+\varepsilon^2\|\zeta_{xx}(\tau)\|^2\right)\,d\tau+C_{\eta}\delta^{\frac{7}{2}}.
\end{eqnarray}

For $\displaystyle\int_0^t\int_{\mathbb{R}}J_6\,dxd\tau$, notice that
\begin{equation*}
  \begin{split}
    |\theta_t| & =\left|\frac{1}{C_v-\frac{\theta}{2}\kappa_{\theta\theta}\frac{v_x^2}{v^5}}\right|\left|-pu_x+(\frac{\tilde{\alpha}(v,\theta)\theta_x}{\theta})_x+\frac{\mu(v,\theta)u_x^2}{v}+\frac{v\kappa_{\theta v}-\kappa_{\theta}}{2v^6}\theta u_xv_x^2+\frac{\theta\kappa_{\theta}v_xu_{xx}}{v^5}\right| \\
      & \le C\delta(|u_x|+|\theta_{xx}|+|v_x\theta_x|+|\theta_x^2|+u_x^2+|u_xv_x^2|+|v_xu_{xx}|),
  \end{split}
\end{equation*}
thus we have
\begin{eqnarray}\label{3.102}
    \left|\int_0^t\int_{\mathbb{R}}J_6\ dxd\tau\right| & \le& C\delta\int_0^t\int_{\mathbb{R}}(|u_x|+|\theta_{xx}|+|v_x\theta_x|+|\theta_x^2|+u_x^2+|u_xv_x^2|+|v_xu_{xx}|)\phi_{xx}^2\ dxd\tau \nonumber\\
      & \le& C\delta\int_0^t\int_{\mathbb{R}}(|\psi_x|+|\zeta_{xx}|+|\phi_x\zeta_x|+|\psi_x\phi_x^2|+|\phi_x\psi_{xx}|+|{ \Theta}_{xx}|+|{ \Theta}_x^2|+\psi_x^2+\zeta_x^2 \nonumber\\
      &&\quad+|\phi_x{ \Theta}_x|+|{ V}_x{\Theta}_x|+|\psi_x{ V}_x^2|+|{ U}_x\phi_x^2|+|\phi_x{ U}_{xx}|+|{ V}_x\psi_{xx}|)\phi_{xx}^2\ dxd\tau.
\end{eqnarray}
For (\ref{3.102}),  we only deal with the most difficulty terms $\int_0^t\int_{\mathbb{R}}|\zeta_{xx}\phi_{xx}^2|\,dxd\tau$ and $\int_0^t\int_{\mathbb{R}}|\phi_x\psi_{xx}\phi_{xx}^2|\,dxd\tau$, the other terms can be estimated similarly as (\ref{3.100})-(\ref{3.101}). In fact,
\begin{equation*}
  \begin{split}
    \int_0^t\int_{\mathbb{R}}\left|\zeta_{xx}\phi_{xx}^2\right|\,dxd\tau & \le\int_0^t\|\phi_{xx}(\tau)\|^{\frac{1}{2}}\|\phi_{xxx}(\tau)\|^{\frac{1}{2}}\|\zeta_{xx}(\tau)\|\sup_{0\le \tau\le t}\{\|\phi_{xx}(\tau)\|\} \,d\tau \\
      & \le N_1\int_0^t\left(\|\phi_{xx}(\tau)\|_1^2+\|\zeta_{xx}(\tau)\|^2\right)d\tau,
  \end{split}
\end{equation*}
\begin{equation*}
  \begin{split}
      \int_0^t\int_{\mathbb{R}}\left|\phi_x\psi_{xx}\phi_{xx}^2\right|\,dxd\tau & \le \int_0^t\sup_{0\le \tau\le t}\{\|\phi_{x}(\tau)\|\|\phi_{xx}(\tau)\|\}\|\psi_{xx}(\tau)\|\|\phi_{xxx}(\tau)\|\,d\tau \\
        & \le N_1^2\int_0^t\|(\phi_{xxx},\psi_{xx})(\tau)\|^2\,d\tau.
    \end{split}
\end{equation*}
Consequently, it holds
\begin{equation}\label{3.103}
  \left|\int_0^t\int_{\mathbb{R}}J_6\ dxd\tau\right|\le CN_1\delta\int_0^t\|(\psi_{xx},\zeta_{xx},\psi_x)(\tau)\|^2\ d\tau+CN_1^4\delta\int_0^t\|\phi_{xx}(\tau)\|_1^2\ d\tau.
\end{equation}
Similarly,
\begin{equation}\label{3.104}
  \begin{split}
   \left |\int_0^t\int_{\mathbb{R}}J_7\ dxd\tau\right|\le& \eta\int_0^t\|\psi_{xx}(\tau)\|^2\ d\tau+C_{\eta}\varepsilon \int_0^t\int_{\mathbb{R}}\left(\phi_{xx}^2\zeta_x^2+\phi_{xx}^2{\Theta}_x^2\right.  \\
    &\left.+{V}_{xx}^2\zeta_x^2+{ V}_{xx}^2{ \Theta}_x^2+\phi_x^2\zeta_{xx}^2+\phi_x^2{ \Theta}_{xx}^2+{ V}_x^2\zeta_{xx}^2+{ V}_x^2{\Theta}_{xx}^2\right)\,dxd\tau \\
      \le& \eta\int_0^t\|\psi_{xx}(\tau)\|^2\,d\tau+C_{\eta}\int_0^t\|(\phi_{xx},\phi_{x},\zeta_x,\sqrt{\delta}\zeta_{xx})(\tau)\|^2\,d\tau+C_{\eta}\delta\\
      &+C_{\eta}\varepsilon\int_0^t\left\{\|\zeta_x(\tau)\|\|\zeta_{xx}(\tau)\|\|\phi_{xx}(\tau)\|^2
      +\|\phi_x(\tau)\|\|\phi_{xx}(\tau)\|\|\zeta_{xx}(\tau)\|^2\right\}d\tau\\
      \le&\eta\int_0^t\|\psi_{xx}(\tau)\|^2\,d\tau+C_{\eta}\int_0^t\|(\phi_{xx},\phi_{x},\zeta_x)(\tau)\|^2\,d\tau+C_{\eta}\delta\\
      &+C_{\eta}\varepsilon\int_0^t\|\phi_{xx}(\tau)\|^4\,d\tau+C_{\eta}N_1^2(\delta+\varepsilon)\int_0^t\|\zeta_{xx}(\tau)\|^2\,d\tau,
  \end{split}
\end{equation}
where we have used the assumption $(\ref{1.16-1})_2$.

Combining  (\ref{3.99})-(\ref{3.104}) and using Corollary 3.2, the a priori assumption (\ref{3.10})-(\ref{3.10-1}) and the smallness of $\eta$, we obtain
\begin{equation}\label{3.1041}
  \begin{split}
     & \|(\phi_{xx},\psi_x)(t)\|^2+\int_0^t\|\psi_{xx}(\tau)\|^2\,d\tau \\
     \le & C\|(\psi_{0x},\phi_{0xx})\|^2+N_{01}^4\int_0^t\|(\phi_x,\psi_x,\zeta_x, \phi_{xx})(\tau)\|^2\,d\tau+C_{30}N_1^2(\delta+\varepsilon)\int_0^t\|\zeta_{xx}(\tau)\|^2\,d\tau
     \\
     &+C_{30}N_1^4\delta\int_0^t\|\phi_{xxx}(\tau)\|^2\,d\tau+C_{30}\int_0^t\|\phi_{xx}(\tau)\|^4\,d\tau\\
      \le & C\|(\phi_0,\psi_0,\frac{\zeta_0}{\sqrt{\delta}},\phi_{0x},\psi_{0x},\phi_{0xx})\|^2+N_{01}^6+C_{30}N_1^4(\delta+\varepsilon+N^2_2\delta)
  +C_{30}\int_0^t\|\phi_{xx}(\tau)\|^4\,d\tau\\
   \le &C\|(\phi_0,\psi_0,\frac{\zeta_0}{\sqrt{\delta}},\phi_{0x},\psi_{0x},\phi_{0xx})\|^2+N_{01}^6+C_{28}\int_0^t\|\phi_{xx}(\tau)\|^4\,d\tau
  \end{split}
\end{equation}
 for  all $t\in[0,T]$, provided that
\[\label{3.92-1}
 C_{30}N_1^4\delta^{\frac{1}{4}}<\frac{1}{3},\quad C_{30}N_1^4\varepsilon<\frac{1}{3},\quad N_2^2\delta^{\frac{3}{4}}<1
\]
holds, where  $C_{30}$ is a positive constant depending only on $\underline V,\overline V,\underline \Theta,\overline \Theta$ and $N_{01}$.

Then  Gronwall's inequality implies that
\begin{equation}\label{3.105}
  \|(\phi_{xx},\psi_x)(t)\|^2+\int_0^t\|\psi_{xx}(\tau)\|^2\,d\tau\le C_{31}N_{0}^2\exp(C_{32}N_{01}^2),
\end{equation}
where $C_{31},C_{32}$ are  positive constants depending only on $\underline {V}, \overline { V},\underline {\Theta},\overline {\Theta}$ and $N_{01}$.

Next, we give the estimate of  $\|\zeta_x(t)\|$.
For this, we   multiply $(3.1)_3$ by $-\zeta_{xx}$ to get
\begin{equation}\label{3.106}
 \frac{R}{2\delta}(\zeta_x^2)_t+\frac{\tilde{\alpha}(v,\theta)}{v}\zeta_{xx}^2=J_8+J_9+J_{10},
\end{equation}
where
\begin{equation*}
  \begin{split}
    J_8 =& p(v,\theta)\psi_x\zeta_{xx}+(p(v,\theta)-P({ V},{\Theta})){ U}_x\zeta_{xx}-\left(\frac{\tilde{\alpha}(v,\theta)}{v}\right)_vv_x\zeta_x\zeta_{xx}-\left(\frac{\tilde{\alpha}(v,\theta)}{v}\right)_{\theta}\theta_x\zeta_x\zeta_{xx} \\
      & -\left(\frac{\tilde{\alpha}(v,\theta)}{v}-\frac{\tilde{\alpha}({ V},{\Theta})}{{ V}}\right)_x{ \Theta}_x\zeta_{xx}-\left(\frac{\tilde{\alpha}(v,\theta)}{v}-\frac{\tilde{\alpha}({ V},{ \Theta})}{{ V}}\right){ \Theta}_{xx}\zeta_{xx}+\frac{\mu({ V},{\Theta})}{{ V}}{ U}_x^2\zeta_{xx}\\
      &+R_2\zeta_{xx}-\frac{\mu(v,\theta)}{v}u_x^2\zeta_{xx},\\
      J_9=&-F\zeta_{xx},\quad J_{10}=-\frac{\theta}{2}\kappa_{\theta\theta}\frac{v_x^2}{v^5}\theta_t\zeta_{xx}.
  \end{split}
\end{equation*}
Integrating (\ref{3.106}) over $[0,t]\times \mathbb{R}$ gives
\begin{equation}\label{3.107}
  \left\|\frac{\zeta_x}{\sqrt{\delta}}(t)\right\|^2+\int_0^t\|\zeta_{xx}(\tau)\|^2\,d\tau\leq C_{33}\left\|\frac{\zeta_{0x}}{\sqrt{\delta}}\right\|^2+C_{33}\left|\int_0^t\int_{\mathbb{R}}\sum_{j=8}^{10}J_i\,dxd\tau\right|,
\end{equation}
where  $C_{33}$ is a positive constant depending only on $\underline V,\overline V,\underline \Theta,\overline \Theta$ and $N_{01}$.

Similar to the proof of (\ref{3.105}), we have
\begin{equation*}
  \begin{split}
    \left|\int_0^t\int_{\mathbb{R}}J_8\,dxd\tau\right| & \le \eta\int_0^t\|\zeta_{xx}(\tau)\|^2\,d\tau+C_{\eta}N_{0}^2 \exp(C_{32}N_{01}^2)\left(\int_0^t\|(\phi_x,\psi_{xx}, \phi_{xx},\psi_x,\zeta_x)(\tau)\|^2\,d\tau+\delta^{\frac{1}{2}}\right),\\
     \left|\int_0^t\int_{\mathbb{R}}J_9\,dxd\tau\right| & \le \eta\int_0^t\|\zeta_{xx}(\tau)\|^2\,d\tau+C_{\eta} N_{0}^4 \exp(4C_{32}N_{01}^2)\int_0^t\|(\psi_{xx},\psi_x,\phi_x)(\tau)\|^2\,d\tau+C_{\eta}\delta^{\frac{5}{2}},\\
    \left|\int_0^t\int_{\mathbb{R}}J_{10}\,dxd\tau\right| &\le C_{34}\delta \int_0^t\|\zeta_{xx}(\tau)\|^2\,d\tau+C_{34}N_{0}^8 \exp(8C_{32}N_{01}^2)\delta\left(\int_0^t\|(\psi_{xx},\psi_x,\phi_x,\zeta_x)(\tau)\|^2\,d\tau+\delta^{\frac{1}{2}}\right),
  \end{split}
\end{equation*}
where  $C_{34}$ is a positive constant depending only on $\underline V,\overline V,\underline \Theta,\overline \Theta$ and $N_{01}$.

Putting the estimates of $\int_0^t\int_{\mathbb{R}}J_i\ dxd\tau, i=8,9,10$ into (\ref{3.107}), and using the smallness of $\eta,\delta$ such that $$C_{33}C_{34}\delta<\frac{1}{2}, \quad 2C_{33}\eta<\frac{1}{2},$$
we obtain
\[
  \left\|\frac{\zeta_x}{\sqrt{\delta}}(t)\right\|^2+\int_0^t\|\zeta_{xx}(\tau)\|^2\,d\tau\le C_{35}N_{0}^8 \exp(8C_{32}N_{01}^2),\quad \forall\, t\in[0,T],
\]
where $C_{35}$ is a positive constant depending only on $\underline {V}, \overline { V},\underline {\Theta},\overline {\Theta}$ and $N_{01}$. Letting $C_{28}=\max\{C_{30}, C_{33}C_{34}\}$ and $C_{29}=C_{35}N_{0}^6 \exp(8C_{32}N_{01}^2)$, then we can get (\ref{3.922}) and hence finish  the proof Lemma 3.5.

Finally, we give the estimate of $\displaystyle\int_0^t\|\phi_{xx}(\tau)\|_1^2\,d\tau$.
\begin{Lemma}
 Under the assumptions of Proposition 3.2, there exists a positive constant $C_{36}$ depending only on $\underline V,\overline V,\underline \Theta,\overline \Theta$ and $N_0$ such that
 \begin{equation}\label{3.108}
  \|\phi_{xx}(t)\|^2+\int_0^t\|\phi_{xx}(\tau)\|_1^2\,d\tau\le C_{36}\left(1+\delta^{-\frac{1}{2}}\right).
\end{equation}
  \end{Lemma}
  \noindent{\bf Proof.}~~Differentiating $(\ref{3.1})_2$ with respect to $x$ once, then multiplying the resultant equation by $\frac{\phi_{xx}}{v}$ and using $(\ref{3.1})_1$, we have
\begin{equation}\label{3.109}
  \begin{split}
    \left(\frac{\mu(v,\theta)\phi_{xx}^2}{2v^2}-\psi_x\frac{\phi_{xx}}{v}\right)_t+\frac{p_+\phi_{xx}^2}{v^2}+\frac{\kappa\phi_{xxx}^2}{v^6}
      =J_{11}+J_{12}+J_{13}+J_{14}+\{\cdots\},
  \end{split}
\end{equation}
where
\begin{equation*}
  \begin{split}
    J_{11} =& \frac{\psi_{xx}^2}{v}-\psi_x\psi_{xx}\frac{v_x}{v^2}+\frac{\psi_x\phi_{xx}(\psi_x+{U}_x)}{v^2} ,\\
     J_{12} =& \left (\frac{R\zeta_x}{v}-2(R\zeta_x-p_+\phi_x)\frac{v_x}{v^2}-\frac{(R\zeta-p_+\phi)v_{xx}}{v^2}+2\frac{(R\zeta-p_+\phi)}{v^3}v_x^2\right)\frac{\phi_{xx}}{v},\\
     J_{13}=&\frac{\mu_v}{2v^2}u_x\phi_{xx}^2+\frac{\mu_{\theta}\theta_t\phi_{xx}^2}{2v^2}-\frac{\mu(v,\theta)}{v^2}{ V}_{txx}\phi_{xx}-\frac{\mu(v,\theta)}{v^3}\phi_{xx}^2+{U}_{tx}\frac{\phi_{xx}}{v}\\
     &-\left(\frac{\mu_{vv}v_x^2+2\mu_{\theta v}v_x\theta_x+\mu_{\theta\theta}\theta_x^2+\mu_vv_{xx}+\mu_{\theta}\theta_{xx}}{v}+\frac{2\mu_vv_xv_{xt}+2\mu_{\theta}\theta_xv_{xt}}{v}\right. \\
     &\left. -2\frac{\mu_vv_x^2v_t+\mu_{\theta}\theta_xv_tv_x+\mu v_{xt}v_x}{v^2}-\frac{\mu u_xv_{xx}}{v^2}+\frac{2\mu u_xv_x^2}{v^3}\right)\frac{\phi_{xx}}{v},
  \end{split}
\end{equation*}
\begin{equation*}
  \begin{split}
    J_{14}=&\frac{\kappa{V}_{xxx}\phi_{xxx}}{v^6}+\left(\frac{(\kappa_vv_x+\kappa_{\theta}\theta_x)v_{xx}}{v^5}-\frac{5\kappa v_xv_{xx}}{v^6}\right)\frac{\phi_{xxx}}{v}\\
     &+\left(\frac{v\kappa_v-5\kappa}{2v^6}v_x^2+\frac{\kappa_{\theta}v_x\theta_x}{v^5}\right)_x\frac{\phi_{xxx}}{v}+\left(\frac{\kappa v_{xx}}{v^5}-\frac{5\kappa -v\kappa_v}{2v^6}v_x^2+\frac{\kappa_{\theta}v_x\theta_x}{v^5}\right)\frac{-v_x\phi_{xx}}{v^2}\\
     :=&\frac{\kappa {V}_{xxx}\phi_{xxx}}{v^6}+J_{14}'.
  \end{split}
\end{equation*}

Similar to (\ref{3.100})-(\ref{3.104}), we have
\[
  \left|\int_0^t\int_{\mathbb{R}}J_{11}\,dxd\tau\right|  \le CN_{0}^2 \exp(C_1N_{01}^2)\int_0^t\|(\phi_{xx},\psi_{xx},\psi_x)(\tau)\|^2\,d\tau,\]
 \[\left|\int_0^t\int_{\mathbb{R}}J_{12}\,dxd\tau\right|  \le CN_{0}^4 \exp(2C_{32}N_{01}^2)\left(\int_0^t\|(\phi_{xx},\phi_{x},\zeta_x)(\tau)\|^2\,d\tau+\delta^{\frac{1}{2}}\right),\]
  \[
  \begin{split}
 \left|\int_0^t\int_{\mathbb{R}}\left(J_{13}+J_{14}'\right)\,dxd\tau\right|
   &\le CN_{0}^6 \exp(6C_{32}N_{01}^2)\left(\int_0^t\|(\phi_x,\psi_{x},\zeta_x)(\tau)\|_1^2\,d\tau+\delta^{\frac{1}{2}}\right)\\
  &\quad+\frac{1}{4}\int_0^t\int_\mathbb{R} \left(\frac{p_+\phi_{xx}^2}{v^2}+\frac{\kappa\phi_{xxx}^2}{v^6}\right)dxd\tau.
 \end{split}\]
\begin{equation}\label{3.113}
  \begin{split}
    \left|\int_0^t\int_{\mathbb{R}}\frac{\kappa {V}_{xxx}\phi_{xxx}}{v^6}\,dxd\tau\right| & \le \frac{1}{4}\int_0^t \int_\mathbb{R}\frac{\kappa\phi_{xxx}^2}{v^6}\,dxd\tau+C\int_0^t\|{\Theta}_{xxx}(\tau)\|^2\,d\tau \\
      &\le \frac{1}{4}\int_0^t\int_\mathbb{R} \frac{\kappa\phi_{xxx}^2}{v^6}\,dxd\tau+C\delta^{-\frac{1}{2}}.
  \end{split}
\end{equation}
Integrating (\ref{3.109}) over $[0,t]\times \mathbb{R}$ and using the estimates of $J_{i}, i=11,12,13 $, $J_{14}^\prime$ and  (\ref{3.113}), then  we can obtain (\ref{3.108}) by Corollary 3.2 and Lemma 3.5. The proof of Lemma 3.6 is completed.

Now we are in a position to prove Proposition 3.2.\\
 \noindent{\bf Proof of Proposition 3.2.}~~With the above lemmas in hand,  we set
\begin{equation}\label{3.114}
\left\{\begin{array}{ll}
    \displaystyle \Xi_1(m_0,M_0; \underline V,\overline V,\underline \Theta,\overline \Theta,N_{01}):=\max\{C_{15},C_{28}\},\\[2mm]
 \displaystyle \Xi_2(\underline V,\overline V,\underline \Theta,\overline \Theta, N_{01}):=C_{28},
  \end{array}\right.
 \end{equation}
where $C_{15}$  and $C_{28}$ are positive constants given in (\ref{3.75-8}) and (\ref{3.92}), respectively.
Notice that
\[ \begin{array}{ll}
C_{15}=\max\{C_{10}, C_{16}, C_{20}\},\\[2mm]
C_{10}=\max\{8C_{12}^{-1}C_{13}C_{14},2C_6\left(\min\{2p_+^2,R^2,\underline{\Theta}\}\right)^{-1}\},\\[2mm]
C_{13}=\max\{2C_{5}C_{7}+C_9\},
  \end{array}
 \]
and the constants $C_{i},i=5, 6, 7, 9, 14, 16, 20$ given in (or in the proof of) the previous Lemmas are increasing functions on both $m_0^{-1}$ and $M_0$,
while the constant $C_{12}$ defined in (\ref{3.75-2}) is  decreasing  on both $m_0$ and $M_0^{-1}$. Consequently, the function $\Xi_1(m_0,M_0; \underline V,\overline V,\underline \Theta,\overline \Theta)$ defined in (\ref{3.114}) is  increasing  on both $m_0^{-1}$ and $M_0$. Thus if (\ref{3.10-2}) holds,  then all  the conditions on  $\delta$ and $\varepsilon$ listed in Lemmas 3.1-3.6 and Corollaries 3.1-3.2 satisfy. Letting $C_3:=\max\{C_{24},C_{29}\}$ and $C_4:=C_{36}$, we conclude from  Corollary 3.2 and Lemmas 3.5-3.6 that (\ref{3.4})-(\ref{3.4-1}) hold  for all $(t,x)\in[0,T]\times\mathbb{R}$. Finally,  from (\ref{3.76}) and (\ref{3.12}),  it easy to see that the inequalities in  (\ref{3.3}) holds for all $(t,x)\in[0,T]\times\mathbb{R}$.  This completes the proof of Proposition 3.2.

\section{Proof of Theorem 1.2}
\setcounter{equation}{0}
This section is devoted to proving Theorem 1.2.
Since the viscous contact wave $(V^c,U^c,\Theta^c)$ satisfy (\ref{1.13}), and the rarefaction waves $(V_{\pm}^r,U_{\pm}^r,\Theta_{\pm}^r)(x,t)$ solve the Euler equations
\begin{equation}\label{4.1}
  \left\{\begin{array}{l}
    (V_{\pm}^r)_t-(U_{\pm}^r)_x=0, \\[2mm]
    (U_{\pm}^r)_t+p(V_{\pm}^r,\Theta_{\pm}^r)_x=0, \\[2mm]
      \displaystyle\frac{R}{\gamma-1}(\Theta_{\pm}^r)_t+p(V_{\pm}^r,\Theta_{\pm}^r)(U_{\pm}^r)_x=0,
  \end{array}\right.
\end{equation}
it is easy to check that the composite wave $(V,U,\Theta)(x,t)$ defined in (\ref{1.27}) satisfy
\begin{equation}\label{4.2}
  \left\{\begin{array}{l}
           V_t-U_x=0, \\[2mm]
           U_t+p_x=\displaystyle\left(\frac{\mu(V,\Theta)U_x}{V}\right)_x+G, \\[3mm]
           \displaystyle\frac{R}{\gamma-1}\Theta_t+PU_x=\left(\frac{{\widetilde\alpha(V,\Theta)}\Theta_x}{V}\right)_x+\frac{\mu(V,\Theta)U_x^2}{V}+H,
         \end{array}\right.
\end{equation}
where
\begin{equation}\label{4.3}
    G=(P-P_--P_+)_x+\displaystyle\left(U_t^c-\left(\frac{\mu(V,\Theta) U_x}{V}\right)_x \right)
\end{equation}
and
\begin{equation}\label{4.4}
\aligned
    H=&\left[(P-p^m)U_x^c+(P-p_-)(U_x^r)_x+(P-p_+)(U_+^r)_x\right] \\[2mm]
    \displaystyle&-\mu\frac{U_x^2}{V}+\left(\frac{{\widetilde \alpha}(V^c,\Theta^c)\Theta_x^c}{V^c}-\frac{{\widetilde \alpha}(V,\Theta)\Theta_x}{V}\right)_x
\endaligned
\end{equation}
representing the interactions and  error terms coming from different wave patterns.

Let
\[(\phi,\psi,\zeta)(t,x)=(v-V,u-U,\theta-\Theta)(t,x),\]
 then we have from (\ref{1.1}) and (\ref{4.1}) that
\begin{equation}\label{4.5}
  \left\{\begin{array}{l}
             \phi_x-\psi_t=0, \\[2mm]
               \displaystyle\psi_t +(v^{-1}(R\zeta-P\phi))_x=\left(\frac{\mu(v,\theta)u_x}{v}-\frac{\mu(V,\Theta)U_x}{V}\right)_x+K_x-G,\\[2mm]
                \displaystyle\left(\frac{R}{\gamma-1}-\frac{\theta}{2}\kappa_{\theta\theta}(v,\theta)\frac{v_x^2}{v^5}\right)\zeta_t+p(v,\theta)u_x-P(V,\Theta)U_x\\[2mm]
            = \displaystyle\left(\frac{{\widetilde \alpha}(v,\theta)}{v}\theta_x-  \displaystyle\frac{{\widetilde  \alpha}(V,\Theta)}{V}\Theta_x\right)_x+\frac{\mu(v,\theta)u_x^2}{v}-\frac{\mu(V,\Theta)U_x^2}{V}+F-H+\frac{\theta}{2}\kappa_{\theta\theta}(v,\theta)\frac{v_x^2}{v^5}\Theta_t
         \end{array}\right.
\end{equation}
with the following initial and far field conditions
\begin{equation}\label{4.6}
  \left\{\begin{array}{l}
            (\phi,\psi,\zeta)(0,x)=(\phi_0,\psi_0,\zeta_0)(x),\\[1mm]
           (\phi,\psi,\zeta)(t,\pm \infty)=0,
         \end{array}\right.
\end{equation}
where the nonlinear terms $K, F, G$ and $H$ are given by $(\ref{1.3})_1$, $(\ref{1.3})_2$, (\ref{4.3}) and (\ref{4.4}), respectively.

The local existence result of the Cauchy problem (\ref{4.4})-(\ref{4.5}) is similar to Proposition 3.1. Then in order to
extend the local solution to be a global one, one needs to derive the  a priori estimates as  Proposition 3.2. To do so,  we first  show the temporal decay   estimate of the nonlinear terms $(G,H)$.
\begin{Lemma}
There exists a uniform constant $C>0$ such that
\begin{equation}\label{4.7}
  \|(G,H)(t)\|_{L^1}\le C\delta^{\frac{1}{8}}(1+t)^{-\frac{7}{8}}.
\end{equation}
\end{Lemma}
 \noindent{\bf Proof.}~~ We only give the  detailed proof of the estimate $\|G(t)\|_{L^1}$, the other one can be treated   similarly. For simplicity, we denote the first and second terms on the right hand side of (\ref{4.3}) by $G_1$ and $G_2$ respectively. Then  a direct calculation yields  that
\begin{equation*}
  \begin{split}
     \frac{G_1}{R}= &\left(\frac{\Theta}{V}-\frac{\Theta_-^r}{V_-^r}-\frac{\Theta_+^r}{V_+^r}-\frac{\Theta^c}{V^c}\right)_x  \\
       =&(\Theta_-^r)_x\left(V^{-1}-(V_-^r)^{-1}\right)+(\Theta_+^r)_x\left(V^{-1}-(V_+^r)^{-1}\right)\\
       & +\Theta_x^c\left(V^{-1}-(V^c)^{-1}\right)+(V_-^r)_x\left(\frac{\Theta_-^r}{(V_-^r)^2}-\frac{\Theta}{V^2}\right)\\
       &+(V_+^r)_x\left(\frac{\Theta_+^r}{(V_+^r)^2}-\frac{\Theta}{V^2}\right)+V_x^c\left(\frac{\Theta^c}{(V^c)^2}-\frac{\Theta}{V^2}\right).
   \end{split}
\end{equation*}
We deduce  from (\ref{1.27}) and Lemma 2.4  that
\begin{equation}\label{4.8}
  \begin{split}
      \left|(\Theta_-^r)_x\left(V^{-1}-(V_-^r)^{-1}\right)\right|
    \le  & C|(\Theta_-^r)_x|\left(|V_+^r-v_+^m|+|V^c-v_-^m|\right)\\
    \le  & C|(\Theta_-^r)_x|\left(|V_+^r-v_+^m|+|V^c-v_-^m|\right)|_{\Omega_-}\\
    & +C|(\Theta_-^r)_x|\left(|V_+^r-v_+^m|+|V^c-v_-^m|\right)|_{\Omega_+ \cup \Omega_c}\\
   \le & C\delta\left[\left(|V_+^r-v_+^m|+|V^c-v_-^m|\right)|_{\Omega_-} +|(\Theta_-^r)_x|_{\Omega_+ \cup \Omega_c}\right]\\
   \le & C\delta^2e^{-c_1(|x|+t)}.
  \end{split}
\end{equation}
Similarly, it holds
\begin{equation*}
 \begin{split}
    & \left|(\Theta_+^r)_x\left(V^{-1}-(V_+^r)^{-1}\right)\right|+\left|(V_{\pm}^r)_x\left(\frac{\Theta_{\pm}^r}{(V_{\pm}^r)^2}-\frac{\Theta}{V^2}\right)\right|\le C\delta^2e^{-c_1(|x|+t)} \\
     & \left|\Theta_x^c\left(V^{-1}-(V^c)^{-1}\right)\right|+\left|V_x^c\left(\frac{\Theta^c}{(V^c)^2}-\frac{\Theta}{V^2}\right)\right|\le C\delta^{\frac{3}{2}}e^{-c_1(|x|+t)}.
 \end{split}
\end{equation*}
Thus  we have
\begin{equation}\label{4.9}
  |G_1|\le C\delta^{\frac{3}{2}}e^{-c_1(|x|+t)}.
\end{equation}
For $G_2$, we have
\begin{equation}\label{4.10}
  \begin{split}
    |G_2|\le & C\left(|U_t^c|+|(U_-^r)_{xx}|+|(U_+^r)_{xx}|+|U_{xx}^c|\right) \\
      & +C\left(|(U_-^r)_x|+|(U_+^r)_x|+|U_x^c|\right)\left(|(V_-^r)_x|+|(V_+^r)_x|+|V_x^c|\right)\\
      :=&G_2^1+G_2^2.
  \end{split}
\end{equation}
Then it follows  from Lemmas 2.1 and 2.4, and the inequality: $\min\{a,b\}\le a^{\theta}b^{1-\theta}, \forall\, \theta\in[0,1]$ that
\begin{equation}\label{4.11}
  \left\|G_2^1\right\|_{L^1}\le C\left(\delta^{\frac{1}{8}}(1+t)^{-\frac{7}{8}}+\delta(1+t)^{-1}\right).
\end{equation}
Similar to (\ref{4.8}), we get
\begin{equation*}
 \begin{split}
   \widetilde{G_2^2} & =|(U_+^r)_x||(V_-^r)_x|+|(U_-^r)_x||(V_+^r)_x|+|U_x^c||(V_{\pm}^r)_x|+|V_x^c||(U_{\pm}^r)_x| \\
     & \le C\delta^{\frac{3}{2}}e^{-c_1(|x|+t)}.
 \end{split}
\end{equation*}
On the other hand, by using Lemmas 2.4 and 2.1 again, we have
\begin{equation*}
  \begin{split}
   & \left\|(U_{\pm}^r)_x(V_{\pm}^r)_x\right\|_{L^1}\le \|(U_{\pm}^r)_x\|_{L^{\infty}}\|(V_{\pm}^r)_x\|_{L^1} \le C\delta(1+t)^{-1},\\[2mm]
     &\left\|U_x^cV_x^c\right\|_{L^1} \le C\delta^2(1+t)^{-\frac{5}{2}}.
  \end{split}
\end{equation*}
Consequently,
\begin{equation}\label{4.12}
    \left\|G_2^2\right\|_{L^1}  \le \left\|\widetilde{G_2^2}\right\|_{L^1} +\left\|(U_{\pm}^r)_x(V_{\pm}^r)_x\right\|_{L^1}+\left\|U_x^cV_x^c\right\|_{L^1} \le C\delta(1+t)^{-1}.
\end{equation}
Then the assertion  for $G$ in (\ref{4.7}) follows from (\ref{4.9})-(\ref{4.12}). This completes the proof of Lemma 4.1.

For the $L^2$-estimate of $(\phi,\psi,\zeta)(t,x)$, we have
\begin{Lemma}There exist a positive constant $C(\underline V,\overline V,\underline \Theta,\overline \Theta)$ and a constant $C_{37}>0$ depending only $\underline V,\overline V,\underline \Theta,\overline \Theta$ and $m_0,M_0$ such that
\begin{equation}\label{4.13}
  \begin{split}
     & \int_{\mathbb R}\left[R\Theta\Phi\left(\frac{v}{V}\right)+\frac{\psi^2}{2}+\frac{R}{\gamma-1}\Theta\Phi\left(\frac{\theta}{\Theta}\right)\right]\ dx+\int_{\mathbb R}\frac{\kappa(v,\theta)v_x^2}{v^5}\ dx \\
      & +\int_0^t\int_{\mathbb R}S_0\left((U_-^r)_x+(U_+^r)_x\right)\ dxd\tau+\int_0^t\int_{\mathbb R}\left(\frac{\mu(v,\theta)\Theta\psi_x^2}{\theta v}+\frac{\widetilde \alpha(v,\theta)\Theta\zeta_x^2}{v\theta^2}\right)\ dxd\tau\\
      \le &C(\underline V,\overline V,\underline \Theta,\overline \Theta)\left\|\left(\phi_0,\phi_{0x},\psi_0,\frac{\zeta_0}{\sqrt \delta}\right)\right\|^2+C_{37}\delta\int_0^t\int_{\mathbb R}(1+\tau)^{-1}\left(\phi^2+\frac{\zeta^2}{\delta}\right)e^{-\frac{c_0x^2}{\delta(1+\tau)}}\ dxd\tau \\
      & +C_{37}\left(N_1^6\delta^{\frac{1}{4}}\int_0^t\|(\phi_x,\phi_{0x},\psi_x,\zeta_x)(\tau)\|^2\ d\tau+N_1^5\delta^{\frac{1}{6}}\right),
  \end{split}
\end{equation}
where $S_0$ is positive function defined in (\ref{4.151}) below.
\end{Lemma}
 \noindent{\bf Proof.}~~Multiplying $(\ref{4.5})_1$ by $\displaystyle R\Theta\left(v^{-1}-V^{-1}\right)$, $(\ref{4.5})_2$ by $\psi$, $(\ref{4.5})_3$ by $\displaystyle\zeta\theta^{-1}$, and adding the resultant equations together, similar to the calculations as Lemma 8 of \cite{F.M. Huang-J. Li-A. Matsumura-2010}, we can get
\begin{equation}\label{4.14}
  \begin{split}
     & \left[R\Theta\Phi\left(\frac{v}{V}\right)+\frac{\psi^2}{2}+\frac{R}{\gamma-1}\Theta\Phi\left(\frac{\theta}{\Theta}\right)\right]_t+Q_0\left((U_-^r)_x+(U_+^r)_x\right)\ +\frac{\mu\Theta\psi_x^2}{\theta v}+\frac{\widetilde \alpha(v,\theta)\Theta\zeta_x^2}{v\theta^2}\\
     = & L_x+\frac{\theta}{2}\kappa_{\theta\theta}\frac{v_x^2}{\theta v^5}\theta_t\zeta+(F-H)\frac{\zeta}{\theta}-G\psi+K_x\psi+Q_1+Q_2,
  \end{split}
\end{equation}
where
\begin{eqnarray}\label{4.151}
    L&=& \left(\frac{\mu(v,\theta)u_x}{v}-\frac{\mu(V,\Theta)U_x}{V}\right)\psi-v^{-1}(R\zeta-P\phi)\psi+\left(\frac{{\widetilde \alpha}(v,\theta)}{v}\theta_x-\frac{{\widetilde \alpha}(V,\Theta)}{V}\Theta_x\right)\frac{\zeta}{\theta}, \nonumber\\
     S_0&=& (\gamma-1)P(V,\Theta)\Phi\left(\frac{v}{V}\right)+\frac{P\phi^2}{vV}-P\Phi\left(\frac{\Theta}{\theta}\right)+\frac{\zeta}{\theta}(p-P)\nonumber\\
     &=& P\left(\Phi\left(\frac{\theta V}{\Theta v}\right)+\gamma\Phi\left(\frac{v}{V}\right)\right)
    \ge C_{38}(\underline V,\overline V,\underline \Theta,\overline \Theta,m_0,M_0)(\phi^2+\zeta^2), \\
    S_1&= &-U_x^c\left[\frac{P(V,\Theta)\phi^2}{vV}-p^m\Phi\left(\frac{v}{V}\right)+\frac{p^m}{\gamma-1}\Phi\left(\frac{\Theta}{\theta}\right)+\frac{\zeta}{\theta}(p-P)\right] \nonumber\\
    &&+(\gamma-1)(P-P_-)(U_-^r)_x\left(\Phi\left(\frac{v}{V}\right)-\frac{1}{\gamma-1}\Phi\left(\frac{\Theta}{\theta}\right)\right) \nonumber\\
    &&+(\gamma-1)(P-P_+)(U_+^r)_x\left(\Phi\left(\frac{v}{V}\right)-\frac{1}{\gamma-1}\Phi\left(\frac{\Theta}{\theta}\right)\right),\nonumber\\
 S_2&=& -\left(\frac{\mu(v,\theta)}{v}-\frac{\mu(V,\Theta)}{V}\right)U_x\psi_x+\left(\frac{\mu(v,\theta)(U_x^2+2U_x\psi_x)}{v}-\frac{\mu(V,\Theta)}{V}U_x^2\right)\frac{\zeta}{\theta} \nonumber\\
    &&+\frac{\zeta_x\Theta_x}{\theta}\left(-\frac{{\widetilde \alpha}(v,\theta)}{v}+\frac{{\widetilde \alpha}(V,\Theta)}{V}\right)+\frac{\Theta_x\theta_x\zeta}{\theta^2}\left(\frac{{\widetilde \alpha}(v,\theta)}{v}-\frac{{\widetilde \alpha}(V,\Theta)}{V}\right)\nonumber\\
    &&+\frac{\widetilde \alpha(v,\theta)\zeta_x\zeta\Theta_x}{v\theta^2}.\nonumber
\end{eqnarray}

A direct computation yields
\begin{equation}\label{4.15}
  \begin{split}
    |S_1| \le & C(\underline V,\overline V,\underline \Theta,\overline \Theta,m_0,M_0)\left(\phi^2+\frac{\zeta^2}{\gamma-1}\right)|U_x^{c}| \\
      & +C(\underline V,\overline V,\underline \Theta,\overline \Theta)\{\left(|\Theta^c-\theta_-^m|+|V^c-v_-^m|+|V_+^r-v_+^m|+|\Theta_+^r-\theta_+^r|\right)|(U_-^r)_x| \\
      & +\left(|V^c-v_+^m|+|V_-^r-v_-^m|+|\Theta^c-\theta_+^m|+|\Theta_-^r-\theta_-^m|\right)|(U_+^r)_x|\}(\phi^2+\zeta^2).
  \end{split}
\end{equation}
The second term can be bounded  in a way similar to (\ref{4.8}), thus we have
\begin{equation*}
  |S_1|\le C(\underline V,\overline V,\underline \Theta,\overline \Theta,m_0,M_0)\delta\left(\phi^2+\frac{\zeta^2}{\delta}\right)(1+t)^{-1}e^{-\frac{c_0x^2}{\delta(1+t)}}+C(\underline V,\overline V,\underline \Theta,\overline \Theta)\delta e^{-c_1(|x|+t)}.
\end{equation*}
It follows from Lemmas 2.1 and 2.4, and the Cauchy inequality that
\begin{equation*}
  \begin{split}
    |S_2|\le & \frac{\mu(v,\theta)\Theta\psi_x^2}{4\theta v}+\frac{\widetilde \alpha(v,\theta)\Theta\zeta_x^2}{4v\theta^2}+C(\underline V,\overline V,\underline \Theta,\overline \Theta,m_0,M_0)(\phi^2+\zeta^2)\left(U_x^2+\Theta_x^2\right)  \\
     \le &\frac{\mu(v,\theta)\Theta\psi_x^2}{4\theta v}+\frac{\widetilde \alpha(v,\theta)\Theta\zeta_x^2}{4v\theta^2}+C(\underline V,\overline V,\underline \Theta,\overline \Theta,m_0,M_0)(\phi^2+\zeta^2)\delta(1+t)^{-1}e^{-\frac{c_0x^2}{\delta(1+t)}}\\
     &+C(\underline V,\overline V,\underline \Theta,\overline \Theta,m_0,M_0)(\phi^2+\zeta^2)\delta\left((U_-^r)_x+(U_+^r)_x\right),
  \end{split}
\end{equation*}
where we have used the fact that
\begin{equation*}\aligned
  |(\Theta_-^r)_x|=&|\theta_-(v_-)^{\gamma-1}(1-\gamma)(V_-^r)^{-\gamma}(V_-^r)_x|\\
  =&\left|\theta_-(v_-)^{\gamma-1}(1-\gamma)(V_-^r)^{-\gamma}\frac{(U_-^r)_x}{\lambda_-(V_-^r,s_-)}\right|\le C|(U_-^r)_x|.
\endaligned\end{equation*}
Therefore, it holds
\begin{equation}\label{4.16}
  \begin{split}
    \int_0^t\int_{\mathbb R}(|S_1|+|S_2|)\ dxd\tau \le & \frac{1}{4}\int_0^t\int_{\mathbb R}\left(\frac{\mu(v,\theta)\Theta\psi_x^2}{\theta v}+\frac{\widetilde \alpha(v,\theta)\Theta\zeta_x^2}{v\theta^2}\right)\,dxd\tau \\
      & +C_{39}\delta\int_0^t\int_{\mathbb R}(1+\tau)^{-1}e^{-\frac{c_0x^2}{\delta(1+\tau)}}\left(\phi^2+\frac{\zeta^2}{\delta}\right)\,dxd\tau \\
      & +C_{39}\delta\int_0^t\int_{\mathbb R}\left((U_-^r)_x+(U_+^r)_x\right)\left(\phi^2+\zeta^2\right)\,dxd\tau,
  \end{split}
\end{equation}
where $C_{39}$ is a positive constant depending only on $\underline V,\overline V,\underline \Theta,\overline \Theta$ and $m_0,M_0$.

Similar to (\ref{3.20}), we have
\begin{equation}\label{4.17}
  \begin{split}
    \int_0^t\int_{\mathbb R}K_x\psi\ dxd\tau =&-\int_{\mathbb R}\frac{\kappa(v,\theta)v_x^2}{2v^5}\,dx+\int_{\mathbb R}\frac{\kappa(v_0,\theta_0)}{2v_0^5}v_{0x}^2\,dx+\int_0^t\int_{\mathbb R}\frac{\kappa_{\theta}\theta_tv_x^2}{2v^5}\,dxd\tau  \\
      & +\int_0^t\int_{\mathbb R}\left\{-\frac{\kappa(v,\theta)v_{xx}}{v^5}+\frac{5\kappa-v\kappa_v}{2v^6}v_x^2-\frac{\kappa_{\theta}v_x\theta_x}{v^5}\right\}V_t\,dx\\
      :=&-\int_{\mathbb R}\frac{\kappa(v,\theta)v_x^2}{2v^5}\,dx+\int_{\mathbb R}\frac{\kappa(v_0,\theta_0)}{2v_0^5}v_{0x}^2\,dx+J_1+J_2.
  \end{split}
\end{equation}
For $J_2$, we derive  from integrations by parts, $(\ref{1.13})_1$ and $(\ref{1.1})_1$ that
\begin{equation}\label{4.18}
  \begin{split}
    J_2= & \int_0^t\int_{\mathbb R}\left[v_x\left(\frac{\kappa(v,\theta)}{v^5}V_t\right)_x+\left(\frac{5\kappa-v\kappa_v}{2v^6}v_x^2-\frac{\kappa_{\theta}v_x\theta_x}{v^5}\right)V_t\right]\,dxd\tau \\
     \le & C\int_0^t\int_{\mathbb R}(\phi_x^2+V_x^2+\zeta_x^2+\Theta_x^2)|V_t|\,dxd\tau+C\int_0^t\int_{\mathbb R}|(\phi_x+V_x)V_{tx}|\,dxd\tau \\
     \le & C\int_0^t\int_{\mathbb R}\left(\phi_x^2+\zeta_x^2+\left|(V_+^r)_x\right|^2+
     \left|(V_-^r)_x\right|^2+|(\Theta_+^r)_x|^2+|(\Theta_-^r)_x|^2+|\Theta_x^c|^2\right)\\
     &\qquad\times\left(|(U_+^r)_x|+|(U_-^r)_x|+|U_x^c|\right)\,dxd\tau \\
     &+C\int_0^t\int_{\mathbb R}\left(|\phi_x|+|(V_+^r)_x|+|(V_-^r)_x|+|V_x^c|\right)\left(|(U_+^r)_{xx}|+|(U_-^r)_{xx}|+|U_{xx}^c|\right)\,dxd\tau \\
     =&J_2^1+J_2^2.
  \end{split}
\end{equation}
The Cauchy inequality and Lemmas 2.1 and 2.4 imply that
\begin{equation}\label{4.19}
  \begin{split}
    J_2^1\le & C\int_0^t\int_{\mathbb R}(\phi_x^2+\zeta_x^2)\left\|\left((U_+^r)_x,(U_-^r)_x,U_x^c\right)\right\|_{L^{\infty}}\,dxd\tau \\
      & +C\int_0^t\int_{\mathbb R}\left(|(U_{\pm}^r)_x|^2+|\Theta_x^c|^2\right)\left\|\left((U_+^r)_x,(U_-^r)_x,U_x^c\right)\right\|_{L^{\infty}}\,dxd\tau \\
      \le & C\delta\int_0^t\|(\phi_x,\zeta_x)(\tau)\|^2\,d\tau+C\delta.
  \end{split}
\end{equation}
Similarly,
\begin{equation}\label{4.20}
  J_2^2\le C\delta^{\frac{1}{4}}\int_0^t\|\phi_x(\tau)\|^2\,d\tau+C\delta^{\frac{1}{4}}.
\end{equation}
The term $J_1$ can be controlled in way similar to (\ref{3.22})-(\ref{3.26}) and (\ref{4.19})-(\ref{4.20}), thus
\begin{equation}\label{4.21}
  J_1\le CN_1^5\delta\int_0^t\|(\phi_{xx},\psi_x,\zeta_x,\phi_x)(\tau)\|^2\,d\tau+CN_1^4\delta^2.
\end{equation}
Combining (\ref{4.17})-(\ref{4.21}) yields
\begin{equation}\label{4.22}
  \begin{split}
    \int_0^t\int_{\mathbb R}K_x\psi\,dxd\tau \le&-\int_{\mathbb R}\frac{\kappa(v,\theta)v_x^2}{2v^5}\,dx+\int_{\mathbb R}\frac{\kappa(v_0,\theta_0)}{2v_0^5}v_{0x}^2\,dx \\
      &+CN_1^5\delta^{\frac{1}{4}}\int_0^t\|(\phi_{xx},\psi_x,\zeta_x,\phi_x)(\tau)\|^2\,d\tau+CN_1^4\delta^{\frac{1}{4}}.
  \end{split}
\end{equation}
Using  (\ref{4.21})  and the inequality
$\|\zeta\|_{L^{\infty}_{T,x}}\le \sup_{t\in[0,T]}\|\zeta(t)\|_1\le N_1\sqrt\delta$, we obtain
\begin{equation}\label{4.23}
  \begin{split}
      \left|\int_0^t\int_{\mathbb R}\frac{\kappa_{\theta\theta}(v,\theta)v_x^2}{2v^5}\theta_t\zeta\,dxd\tau \right|
     \le & C\|\zeta(t)\|_{L^{\infty}_{T,x}}\int_0^t\int_{\mathbb R}\left|\frac{\kappa_\theta v_x^2\theta_t}{v^5}\right|\,dxd\tau \\
     \le & CN_1^6\delta^{\frac{3}{2}}\int_0^t\|(\phi_{xx},\psi_x,\zeta_x,\phi_x)(\tau)\|^2\,d\tau+CN_1^5\delta^{\frac{5}{2}}.
  \end{split}
\end{equation}
Moreover, it follows from the Cauchy inequality, the Sobolev inequality and Lemma 4.1 that
\begin{equation}\label{4.24}
  \begin{split}
     &\left|\int_0^t\int_{\mathbb R}G\psi+\frac{\zeta}{\theta}H\,dxd\tau\right| \\
     \le & \frac{1}{4}\int_0^t\int_{\mathbb R}\left(\frac{\mu(v,\theta)\Theta}{v\theta}\psi_x^2+\frac{\widetilde \alpha(v,\theta)\Theta\zeta_x^2}{v\theta^2}\right)\,dxd\tau+C\int_0^t\|(\phi,\zeta)(\tau)\|^{\frac{2}{3}}\|(G,H)(\tau)\|_{L^1}^{\frac{4}{3}}\,d\tau\\
     \le &\frac{1}{4}\int_0^t\int_{\mathbb R}\left(\frac{\mu(v,\theta)\Theta}{v\theta}\psi_x^2+\frac{\widetilde \alpha(v,\theta)\Theta\zeta_x^2}{v\theta^2}\right)\,dxd\tau+CN_1^{\frac{2}{3}}\delta^{\frac{1}{6}}.
  \end{split}
\end{equation}
Finally, similar to (\ref{3.31}), we have
\begin{equation}\label{4.25}
  \left|\int_0^t\int_{\mathbb R}\frac{F}{\theta}\zeta\,dxd\tau\right|\le CN_1^3\delta^{\frac{1}{2}}\int_0^t\|(\phi_{xx},\psi_x,\zeta_x,\phi_x)(\tau)\|^2\,d\tau+CN_1\delta.
\end{equation}
Integrating (\ref{4.14}) in $t$ and $x$ over $[0,t]\times \mathbb R$, and  using (\ref{4.16}), (\ref{4.22})-(\ref{4.25}), we can get (\ref{4.13}). This completes the proof Lemma 4.2.

 For the remainder term $\int_0^t\int_{\mathbb R}(1+\tau)^{-1}\left(\phi^2+\frac{\zeta^2}{\delta}\right)e^{-\frac{c_0x^2}{\delta(1+\tau)}}\, dxd\tau$ in (\ref{4.13}), similar to Lemma 3.2,  one can show that there exist two positive constants $C_{40},C_{41}$ depending only on $\underline V,\overline V,\underline \Theta,\overline \Theta,m_0$ and $M_0$  such that if
 \begin{equation*}
   C_{40}N_1^2\delta^{\frac{1}{2}}<\frac{1}{2}\min\{2P_+,R^2,\underline \Theta\},
 \end{equation*}
 it holds for all $t\in[0,T]$ that
\begin{equation*}
  \begin{split}
   \int_0^t\int_{\mathbb R}\left(\phi^2+\psi^2+\frac{\zeta^2}{\delta}\right)w^2\,dxd\tau\le & C_{41}N_1^6\delta^{-\frac{3}{4}}+C_{41}N_1^7\delta^{-\frac{3}{4}}\int_0^t\|(\phi_{xx},\psi_x,\zeta_x,\phi_x)(\tau)\|^2\,d\tau \\
      & +C_{41}N_1^2\delta\int_0^t\int_ {\mathbb R}(\phi^2+\zeta^2)\left((U_-^r)_x+(U_+^r)_x\right)\,dx.
      \end{split}
\end{equation*}

Moreover, it is easy to check  that some similar estimates as  Lemmas 3.3-3.6 and Corollaries 3.1-3.2 still hold for the solutions to the Cauchy problem (\ref{4.5})-(\ref{4.6}). Thus we can get the desired a priori estimates as Proposition 3.2. Then similar to the proof of  Theorem 3.1,  the global-in-time solutions to problem (\ref{4.5})-(\ref{4.6}) can  also be obtained.  Hence the proof Theorem 1.2 is completed.

\section{Appendix}
\setcounter{equation}{0}
\noindent{\it The proof of Lemma 3.2.}~~The proof of (\ref{3.33}) is divided into the following two parts:
\begin{equation}\label{3.34}
  \begin{split}
    \int_0^t\int_{\mathbb{R}}(R\zeta+\delta p_+\phi)^2w^2\ dxd\tau\le & C_6N_1^2\delta^{\frac{3}{2}}\int_0^t\int_{\mathbb{R}}\left(\phi^2+\frac{\zeta^2}{\delta}\right)w^2\ dxd\tau \\
     &+C_6N_1^7\delta^{\frac{1}{2}}\int_0^t\|(\phi_x,\phi_{xx},\psi_x,\zeta_x)(\tau)\|^2\ d\tau+C_6N_1^6\delta,
  \end{split}
\end{equation}
and
\begin{equation}\label{3.35}
  \int_0^t\int_{\mathbb{R}}[(R\zeta-p_+\phi)^2+\psi^2]w^2\ dxd\tau\le C_8N_1^5\delta^{\frac{1}{4}}+ C_8N_1^6\delta^{\frac{1}{4}}\int_0^t\|(\phi_x,\phi_{xx},\psi_x,\zeta_x)(\tau)\|^2\ d\tau.
\end{equation}
Here $C_6$ and $C_8$ are two positive constants depending only on $\underline V,\overline V,\underline \Theta,\overline\Theta, m_0$ and
$M_0$,  and without loss of generality, we may assume that $C_6\geq c_5m_0^{-1}$ with $c_5$ being a positive constant given in (\ref{3.47}) below.

In fact, notice that
\begin{equation}\label{3.36}
  (R\zeta+\delta p_+\phi)^2+(R\zeta-p_+\phi)^2\geq\delta\left(2p_+^2\phi^2+\frac{\left(R\zeta\right)^2}{\delta}\right)
  \geq\min\{2p_+^2,R^2\}\delta\left(\phi^2+\frac{\zeta^2}{\delta}\right),
\end{equation}
then adding (\ref{3.34}) onto (\ref{3.35}), we can get (\ref{3.33}) by  the assumption (\ref{3.33-1}).

To prove (\ref{3.36}), we denote
\begin{equation*}
  f=\int_{-\infty}^xw^2\ dy,
\end{equation*}
then it holds that
\begin{equation}\label{3.37}
  \|f(t)\|_{L^{\infty}}\le C\delta^{\frac{1}{2}}(1+t)^{-\frac{1}{2}},\quad \|f_t(t)\|_{L^{\infty}}\le C\delta^{\frac{1}{2}}(1+t)^{-\frac{3}{2}}.
\end{equation}\
 We rewrite $(\ref{3.1})_2$ as
\begin{equation}\label{3.38}
  \psi_t+\left(\frac{R\zeta-p_+\phi}{v}\right)_x=\left(\frac{\mu(v,\theta)\psi_x}{v}\right)_x+K_x+G,\quad G=-{U}_t+\left(\frac{\mu(v,\theta){U}_x}{v}\right)_x.
\end{equation}
Multiplying  (\ref{3.38}) by $(R\zeta-p_+\phi)vf$ and integrating the resulting equation over gives
\begin{equation*}
  \begin{split}
      & \int_{\mathbb{R}}\psi_t(R\zeta-p_+\phi)vf\,dx+\int_{\mathbb{R}}\left(\frac{R\zeta-p_+\phi}{v}\right)_x\psi_t(R\zeta-p_+\phi)vf\,dx \\
     =  & \int_{\mathbb{R}}\left(\frac{\mu(v,\theta)\psi_x}{v}\right)_x(R\zeta-p_+\phi)vf\,dx+\int_{\mathbb{R}}K_x(R\zeta-p_+\phi)vf\,dx+\int_{\mathbb{R}}G(R\zeta-p_+\phi)vf\,dx.
   \end{split}
\end{equation*}
Using integrating by parts, we have
\begin{eqnarray}\label{3.39}
  \frac{1}{2}\int_{\mathbb{R}}(R\zeta-p_+\phi)^2w^2\,dx &=&\int_{\mathbb{R}}\psi_t(R\zeta-p_+\phi)vf\,dx-\int_{\mathbb{R}}(R\zeta-p_+\phi)^2v^{-1}v_xf\,dx \nonumber\\
      && +\int_{\mathbb{R}}\mu(v,\theta)\psi_xv^{-1}[(R\zeta-p_+\phi)vf]_x\,dx\nonumber\\
      && -\int_{\mathbb{R}}K_x(R\zeta-p_+\phi)vf\,dx-\int_{\mathbb{R}}G(R\zeta-p_+\phi)vf\,dx\nonumber\\
      &=&\left(\int_{\mathbb{R}}\psi(R\zeta-p_+\phi)vf\ dx\right)_t-\int_{\mathbb{R}}\psi(R\zeta-p_+\phi)_tvf\,dx\nonumber\\
      && -\int_{\mathbb{R}}\psi(R\zeta-p_+\phi)v_tf\,dx-\int_{\mathbb{R}}\psi(R\zeta-p_+\phi)vf_t\,dx\nonumber\\
      && -\int_{\mathbb{R}}(R\zeta-p_+\phi)^2v^{-1}v_xf\,dx+\int_{\mathbb{R}}\mu(v,\theta)\psi_xv^{-1}[(R\zeta-p_+\phi)vf]_x\,dx\nonumber\\
      && -\int_{\mathbb{R}}K_x(R\zeta-p_+\phi)vf\ dx-\int_{\mathbb{R}}G(R\zeta-p_+\phi)vf\,dx\nonumber\\
      &=&\left(\int_{\mathbb{R}}\psi(R\zeta-p_+\phi)vf\ dx\right)_t+\sum_{i=3}^{9}I_i.
\end{eqnarray}
Now we estimate the terms $I_1,I_2,\cdots I_9$ one by one. First, the equations $(\ref{3.1})_1$  and $(\ref{3.1})_2$ imply
\begin{equation}\label{3.40}
  \begin{split}
    \left(\frac{R}{\gamma-1}\zeta+p_+\phi\right)_t =&-\frac{R\zeta-p_+\phi}{v}(\psi_x+{ U}_x)+\frac{\theta}{2}\kappa_{\theta\theta}(v,\theta)\frac{v_x^2}{v^5}\theta_t+\left(\frac{\tilde{\alpha}(v,\theta)\theta_x}{v}-\frac{\tilde{\alpha}({ V},{\Theta}){\Theta}_x}{{V}}\right)_x \\
      & +\frac{\mu(v,\theta)(\psi_x+{U}_x)^2}{v}+F,
  \end{split}
\end{equation}
thus it holds
\begin{equation}\label{3.41}
  \begin{split}
    I_3 =&-\delta\int_{\mathbb{R}}\psi\left(\frac{R}{\delta}\zeta+p_+\phi\right)_tvf \,dx+\gamma p_+\int_{\mathbb{R}}\psi vf\psi_x\,dx \\
      =&\delta\int_{\mathbb{R}}\psi f(R\zeta-p_+\phi)(\psi_x+{U}_x)\, dx-\delta\int_{\mathbb{R}}\psi f\frac{\theta}{2}\kappa_{\theta\theta}(v,\theta)\frac{v_x^2}{v^4}\theta_t\,dx\\
      & -\delta\int_{\mathbb{R}}\psi vf \left(\frac{\tilde{\alpha}(v,\theta)\theta_x}{v}-\frac{\tilde{\alpha}({ V},{\Theta}){ \Theta}_x}{{ V}}\right)_x\,dx-\delta\int_{\mathbb{R}}\psi f\mu(v,\theta)(\psi_x+{U}_x)^2\,dx\\
      & -\delta\int_{\mathbb{R}}\psi vfF\,dx+\gamma p_+\int_{\mathbb{R}}vf\left(\frac{\psi^2}{2}\right)_x\,dx\\
:=&\sum_{i=1}^{6}I_3^i.
  \end{split}
\end{equation}
Using $(\ref{3.37})_1$, the Cauchy inequality, the Young  inequality  and the Sobolev  inequality, we have
\begin{equation}\label{3.42}
  \begin{split}
     |I_4|+\left|I_3^1\right| & \le C\delta^{\frac{1}{2}}(1+t)^{-\frac{1}{2}}\|\psi(t)\|^{\frac{1}{2}}\|\psi_x(t)\|^{\frac{1}{2}}\|(\psi,\zeta)(t)\|\|(\psi_x,{U}_x)(t)\| \\
       & \le CN_1^{\frac{3}{2}}\delta^{\frac{1}{2}}\left(\|\psi_x(t)\|^2+(1+t)^{-\frac{5}{3}}\right).
   \end{split}
\end{equation}
From the estimate of $I_1$, we obtain
\begin{equation}\label{3.43}
  \begin{split}
    \left|I_3^2\right| &\le \delta\left\|\psi f\frac{\theta}{2}\frac{\kappa_{\theta\theta}(v,\theta)}{\kappa\theta}v\right\|_{L^{\infty}}\left|\int_{\mathbb{R}}\frac{\kappa_{\theta}\theta_tv_x^2}{2v^5}\,dx\right| \\
      & \le C\delta^{\frac{5}{2}}(1+t)^{-\frac{1}{2}}N_1^6\|(\phi_{xx},\phi_x,\psi_x,\zeta_x)(t)\|^2+CN_1^5\delta^3(1+t)^{-\frac{3}{2}}.
  \end{split}
\end{equation}
Similarly, it holds that
\begin{equation}\label{3.44}
  \begin{split}
    \left|I_3^3\right| &= \left|-\delta\int_{\mathbb{R}}(\psi vf )_x \left (\frac{\tilde{\alpha}(v,\theta)\zeta_x}{v}+\left(\frac{\tilde{\alpha}(v,\theta)}{v}-\frac{\tilde{\alpha}({ V},{\Theta})}{{ V}}\right){ \Theta}_x\right)\ dx\right| \\
      & \le C\delta\int_{\mathbb{R}}(|\psi_x vf|+|\psi v_xf|+|\psi vf_x|)(|\zeta_x|+|\phi{\Theta}_x|+|\zeta{\Theta}_x|)\ dx\\
      & \le C\delta N_1^2(\|(\phi_x,\psi_x,\zeta_x)(t)\|^2+(1+t)^{-\frac{3}{2}}),
  \end{split}
\end{equation}
\begin{equation}\label{3.45}
     \left|I_3^4\right|\le C\delta\|\psi\|_{L^{\infty}}\|f\|_{L^{\infty}}\|(\psi_x,{ U}_x)\|^2 \le CN_1\delta^{\frac{3}{2}}\|\psi_x(t)\|^2+CN_1\delta^4(1+t)^{-2},
\end{equation}
\begin{equation}\label{3.46}
  \begin{split}
     \left|I_3^5\right| =&\left|\delta\int_{\mathbb{R}}[\psi_xv_xv^{-4}f\theta+\psi f_xv_xv^{-4}\theta\kappa_{\theta}+\psi fv_{xx}v^{-4}\theta\kappa_{\theta}+\right. \\
    &\left. \psi fv_x((v^{-4}\theta\kappa_{\theta})_vv_x
+(v^{-4}\theta\kappa_{\theta})_{\theta}\theta_x)]u_x\ dx \right|+\delta\left|\int_{\mathbb{R}}\psi vf\frac{v\kappa_{\theta v}-\kappa_{\theta}}{2v^6}\theta u_xv_x^2\ dx\right| \\
      \le& C\delta\int_{\mathbb{R}}(|\psi_xv_xf u_x|+|\psi f_xv_xu_x|+|\psi fv_{xx}u_x|+|\psi fv_x^2u_x|+|\psi fv_x\theta_xu_x|)\,dx\\
      \le & C\delta\left[\|v_x\|_{L^{\infty}}\|f\|_{L^{\infty}}\|(\psi_x,{ U}_x)\|^2+\|\psi\|_{L^{\infty}}\|v_x\|_{L^{\infty}}\|w^2\|\|(\psi_x,{U}_x)\|\right.\\
      &\left.+\|\psi\|_{L^{\infty}}\|f\|_{L^{\infty}}\left(\|(\phi_{xx},\psi_x,{ U}_x,{ V}_{xx},{ V}_x^2)\|^2+\|(\psi_x,{ U}_x)\|_{L^{\infty}}\|\phi_x\|^2\right)\right]\\
      \le& CN_1^2\delta^{\frac{5}{4}}\|(\phi_{xx},\psi_x,\phi_x)(t)\|^2+CN_1^2\delta^{\frac{5}{4}}(1+t)^{-\frac{3}{2}},
  \end{split}
\end{equation}
and
\begin{equation}\label{3.47}
 \begin{split}
   I_3^6 & =-\gamma p_+\int_{\mathbb{R}}(vf)_x\frac{\psi^2}{2}\,dx \\
     & \le \frac{\gamma p_+}{2}\int_{\mathbb{R}}(|\phi_x|+|V_x|)|f||\psi|^2\,dx-\frac{\gamma p_+}{2}\int_{\mathbb{R}}vw^2\psi^2\,dx\\
     & \le \frac{\gamma p_+}{2}c_5m_0^{-1}\delta\int_{\mathbb{R}}vw^2\psi^2\,dx-\frac{\gamma p_+}{2}\int_{\mathbb{R}}vw^2\psi^2\,dx+C\delta^{\frac{1}{2}}\|\phi_x\|\|\psi\|^{\frac{3}{2}}\|\psi_x\|^{\frac{1}{2}}(1+t)^{-\frac{1}{2}}\\
     & \le -\frac{\gamma p_+}{4}\int_{\mathbb{R}}vw^2\psi^2\,dx+CN_1^{\frac{3}{2}}\delta^{\frac{1}{2}}(\|(\phi_x,\psi_x)(t)\|^2+(1+t)^{-2}).
 \end{split}
\end{equation}
Here in (\ref{3.47}), the constant $c_5>0$ depends only on $\theta_{\pm}$, and satisfies $|V_xf|\leq c_5w^2$. Moreover, we have used the smallness of $\delta$ such that $c_5m_0^{-1}\delta<\frac{1}{2}$.

Similar to the estimates as above, we also have
\begin{equation}\label{3.48}
    |I_5| \le C\|\psi(t)\|\|(\zeta,\phi)(t)\|\|f_t(t)\|_{L^{\infty}} \le CN_1^2\delta^{\frac{1}{2}}(1+t)^{-\frac{3}{2}},
\end{equation}
\begin{equation}\label{3.50}
  |I_7|\le CN_1\delta^{\frac{1}{4}}\|(\phi_{x},\psi_x,\zeta_x)(t)\|^2+CN_1\delta^{\frac{1}{4}}(1+t)^{-\frac{3}{2}},
\end{equation}
\begin{equation}\label{3.51}
  \begin{split}
    |I_8| & =\left|\int_{\mathbb{R}}K[(R\zeta-p_+\phi)vf]_x\ dx\right| \\
      & \le C\int_{\mathbb{R}} \left(|\phi_{xx}|+|{V}_{xx}|+|\phi_x^2|+|{ V}_x^2|+|\zeta_x^2|+{ \Theta}_x^2\right)\\
      &\qquad\times\left(|(\phi_x,\zeta_x)f|+|(\phi,\zeta)(\phi_x+{ V}_x)f|+|(\phi,\zeta)f_x|\right)\,dx\\
      & \le CN_1^2\delta^{\frac{1}{4}}\|(\phi_{xx},\phi_x,\zeta_x)(t)\|^2+CN_1^2\delta^{\frac{1}{4}}(1+t)^{-\frac{3}{2}},
  \end{split}
\end{equation}
\begin{equation}\label{3.52}
  \begin{split}
    |I_9|&\le \int_{\mathbb{R}}|G||(\phi,\zeta)||f|\ dx \le \int_{\mathbb{R}}(|{ U}_t|+|{ U}_{xx}|+|v_x{ U}_x|+|\theta_x{ U}_x|)|(\phi,\zeta)||f|\,dx\\
      &\le CN_1\delta^{\frac{3}{2}}(\|(\phi_x,\zeta_x)(t)\|^2+(1+t)^{-\frac{3}{2}}).
  \end{split}
\end{equation}
and
\begin{equation}\label{3.49}
  \begin{split}
    |I_6|&\le C\int_{\mathbb{R}}|(\phi,\zeta)|^2|\phi_x||f|\,dx+\left|\int_{\mathbb{R}}v^{-1}(R\zeta-p_+\phi)^2{ V}_xf\,dx\right|\\
      & \le  C\|(\phi,\zeta)(t)\|^{\frac{3}{2}}\|(\phi_x,\zeta_x)(t)\|^{\frac{1}{2}}\|\phi_x(t)\|\delta^{\frac{1}{2}}(1+t)^{-\frac{1}{2}}+c_5m_0^{-1}\delta\int_{\mathbb{R}}|R\zeta-p_+\phi|^2w^2\,dx\\
      &\le CN_1^{\frac{3}{2}}\delta^{\frac{1}{2}}\|(\phi_x,\zeta_x)(t)\|^2+CN_1^{\frac{3}{2}}\delta^{\frac{1}{2}}(1+t)^{-2}+\frac{1}{4}\int_{\mathbb{R}}|R\zeta-p_+\phi|^2w^2\,dx
  \end{split}
\end{equation}
provided that $\delta$ is sufficiently small such that $c_5m_0^{-1}\delta<\frac{1}{4}$.

Combining (\ref{3.39}), (\ref{3.41})-(\ref{3.52})  and integrating the resulting equation over $[0,t]$,  we can get (\ref{3.35}) by the smallness of $\delta$.

Next, we prove (\ref{3.34}).
Let $h=R\zeta+\delta p_+\phi$ and using (\ref{3.40}),  we have
\begin{eqnarray}\label{3.53}
    \frac{1}{\delta}<h_t,hg^2>_{H^{-1}\times H^1} &= & \int_{\mathbb{R}}(\frac{R}{\gamma-1}\zeta+p_+\phi)_thg^2\,dx\nonumber\\
      &=& -\int_{\mathbb{R}}\frac{R\zeta-p_+\phi}{v}\psi_xhg^2\,dx-       \int_{\mathbb{R}}\frac{R\zeta-p_+\phi}{v}{U}_xhg^2\,dx\nonumber\\
      && +\int_{\mathbb{R}}\frac{\theta}{2}\kappa_{\theta\theta}(v,\theta)\frac{v_x^2}{v^5}\theta_thg^2\,dx+\int_{\mathbb{R}} \left(\frac{\tilde{\alpha}(v,\theta)\theta_x}{v}-\frac{\tilde{\alpha}({V},{\Theta}){\Theta}_x}{{V}}\right)_xhg^2\,dx\nonumber\\
      &&+\int_{\mathbb{R}}\frac{\mu(v,\theta)(\psi_x+{U}_x)^2}{v}hg^2\ dx +\int_{\mathbb{R}}Fhg^2\,dx\nonumber\\
      &:=& \sum_{i=10}^{15}I_i.
\end{eqnarray}
Since $|{U}_x|\le C\delta w^2, h=R\zeta+\delta p_+\phi=O(1)\delta^{\frac{1}{2}}(\phi+\frac{\zeta}{\sqrt{\delta}})$, it holds that
\begin{equation}\label{3.54}
    |I_{11}|\le C\int_{\mathbb{R}}\delta^{\frac{1}{2}}\left|\left(\phi,\frac{\zeta}{\sqrt{\delta}}\right)\right||{U}_x||(\zeta,\phi)|\|g\|_{L^{\infty}}^2\,dx \le C\delta^{\frac{5}{2}}\int_{\mathbb{R}}\left|\left(\phi,\frac{\zeta}{\sqrt{\delta}}\right)\right|^2w^2\,dx,
\end{equation}
\begin{equation}\label{3.55}
    |I_{12}| \le C\left|\int_{\mathbb{R}}\frac{\kappa_{\theta}v_x^2}{2v^5}\theta_t\,dx\right|\|hg^2\|_{L^{\infty}} \le C\delta^2N_1^6\int_{0}^{t}\|(\phi_{xx},\phi_x,\psi_x,\zeta_x)\|^2+o(1)N_1^5\delta^{\frac{5}{2}}(1+t)^{-\frac{3}{2}},
\end{equation}
\begin{equation}\label{3.56}
  \begin{split}
    |I_{13}| & =\left|\int_{\mathbb{R}} \left(\frac{\tilde{\alpha}(v,\theta)\theta_x}{v}-\frac{\tilde{\alpha}({V},{\Theta}){\Theta}_x}{{ V}})(hg^2\right)_x\,dx\right| \\
      & \le C\int_{\mathbb{R}}(|\zeta_x|+|(\phi,\zeta){ \Theta}_x|)\left(|(\zeta_x,\phi_x)|\|g\|_{L^{\infty}}^2+\delta^{\frac{1}{2}}\left|\left(\phi,\frac{\zeta}{\sqrt{\delta}}\right)w\right|\|g\|_{L^{\infty}}\right)\,dx\\
      &\le C\delta^{\frac{1}{2}}\|(\phi_x,\zeta_x)(t)\|^2+C\delta^{\frac{3}{2}}\int_{\mathbb{R}}\left|\left(\phi,\frac{\zeta}{\sqrt{\delta}}\right)\right|^2w^2\,dx,
  \end{split}
\end{equation}

\begin{equation}\label{3.57}
  \begin{split}
    |I_{14}|&\le  C\int_{\mathbb{R}}|(\psi_x,{ U}_x)|^2|hg^2|\,dx \le C\int_{\mathbb{R}}|(\psi_x,{U}_x)|^2\delta^{\frac{1}{2}}\left|\left(\phi,\frac{\zeta}{\sqrt{\delta}}\right)\right|\|g\|_{L^{\infty}}^2\,dx\\
       &\le CN_1\delta^{\frac{3}{2}}\|\psi_x(t)\|^2+CN_1\delta^4(1+t)^{-\frac{3}{2}}.
  \end{split}
\end{equation}
Similar to the estimate of $I_3^5$, we obtain
\begin{equation}\label{3.58}
  \begin{split}
    |I_{15}|\le&\left |\int_{\mathbb{R}}\left(\frac{\theta\kappa_{\theta}v_x}{v^5}(R\zeta+\delta p_+\phi)g^2\right)_xu_x\,dx\right|+\left|\int_{\mathbb{R}}\frac{v\kappa_{\theta v}-\kappa_{\theta}}{2v^6}\theta u_xv_x^2(R\zeta+\delta p_+\phi)^2g^2\,dx\right| \\
\le& C\int_{\mathbb{R}}(|\theta_xv_xu_x|+|v_x^2u_x|+|v_{xx}u_x|)|(\zeta,\phi)||g^2|\,dx\\
      &+C\int_{\mathbb{R}}(|v_xu_x(\zeta_x,\phi_x)g^2|+|v_xu_x(\zeta,\phi)gg_x|)\,dx\\
\le& CN_1^2\delta^{\frac{1}{2}}\|(\phi_{xx},\phi_x,\psi_x,\zeta_x)(t)\|^2+CN_1^2\delta(1+t)^{-\frac{3}{2}}.
  \end{split}
\end{equation}
For the estimate of $I_{10}$, we compute from (\ref{3.40}) that
\begin{eqnarray}\label{3.59}
    -2I_{10} &=&2\int_{\mathbb{R}}v^{-1}(h^2-\gamma p_+\phi h)\phi_t g^2\,dx \nonumber\\
      &=& \int_{\mathbb{R}}(2v^{-1}h^2g^2\phi_t-\gamma p_+v^{-1}hg^2(\phi^2)_t)\ dx\nonumber\\
&=&\left(\int_{\mathbb{R}}v^{-1}hg^2\phi(2h-\gamma p_+\phi)\ dx\right)_t-2\int_{\mathbb{R}}v^{-1}hg\phi(2h-\gamma p_+\phi)g_t\,dx\nonumber\\
      &&+\int_{\mathbb{R}}v^{-2}v_tg^2 h\phi(2h-\gamma p_+\phi)\ dx-\int_{\mathbb{R}}v^{-1}g^2\phi(4h-\gamma p_+\phi)h_t\,dx\nonumber\\
&=&\left(\int_{\mathbb{R}}v^{-1}hg^2\phi(2h-\gamma p_+\phi)\,dx\right)_t-\frac{\delta}{2\alpha}\int_{\mathbb{R}}v^{-1}hg\phi(2h-\gamma p_+\phi)w_x\,dx\nonumber\\
      &&+\int_{\mathbb{R}}v^{-2}u_xg^2\phi[h(2h-\gamma p_+\phi)+\delta(4h-\gamma p_+\phi)(R\zeta-p_+\phi)]\,dx\nonumber\\
      &&-\delta\int_{\mathbb{R}}v^{-1}g^2\phi(4h-\gamma p_+\phi)\frac{\theta}{2}\kappa_{\theta\theta}\frac{v_x^2}{v^5}\theta_t\,dx\nonumber\\
      &&-\delta\int_{\mathbb{R}}v^{-1}g^2\phi(4h-\gamma p_+\phi)\left(\frac{\tilde{\alpha}(v,\theta)\theta_x}{v}-\frac{\tilde{\alpha}({V},{ \Theta}){ \Theta}_x}{{V}}\right)_x\,dx\nonumber\\
      &&-\delta\int_{\mathbb{R}}v^{-1}g^2\phi(4h-\gamma p_+\phi)\frac{\mu(v,\theta)(\psi_x+{ U}_x)^2}{v}\,dx-\delta\int_{\mathbb{R}}v^{-1}g^2\phi(4h-\gamma p_+\phi)F\,dx\nonumber\\
      :&=&\left(\int_{\mathbb{R}}v^{-1}hg^2\phi(2h-\gamma p_+\phi)\ dx\right)_t+\sum_{i=1}^{6}I_{10}^i,
\end{eqnarray}
By employing the Sobolev inequality, the Young inequality, the a priori assumption $(\ref{3.10})$ and Lemma 2.1, we can control the terms $I_{10}^i, i=1,\cdots 6$ as follows.
\begin{equation}\label{3.60}
  \begin{split}
    \left|I_{10}^1\right| &\le  C\delta(1+t)^{-1}\int_{\mathbb{R}}|(\phi,\zeta)|^3\,dx  \le  C\delta(1+t)^{-1}\|(\phi,\zeta)(t)\|^{\frac{5}{2}}\|(\phi_x,\zeta_x)(t)\|^{\frac{1}{2}}\\
      &\le  C\delta\left(\|(\phi_x,\zeta_x)(t)\|^2+N_1^{\frac{10}{3}}(1+t)^{-\frac{4}{3}}\right),
  \end{split}
\end{equation}
\begin{equation}\label{3.61}
  \begin{split}
    \left|I_{10}^2\right| & \le  C\int_{\mathbb{R}}\delta|({U}_x,\psi_x)|(|\phi^3|+|\zeta^3|)\,dx \le C\delta\|(\phi,\zeta)(t)\|^2\|(\phi_x,\zeta_x)(t)\|\|({ U}_x,\psi_x)(t)\|\\
      &\le C\delta N_1^2\left(\|(\phi_x,\zeta_x,\psi_x)(t)\|^2+\delta^{\frac{5}{2}}(1+t)^{-\frac{3}{2}}\right),
  \end{split}
\end{equation}
\begin{equation}\label{3.62}
  \begin{split}
    \left|I_{10}^3\right| & \le C\delta^2\int_{\mathbb{R}}\|(\phi,\zeta)(t)\|_{L^{\infty}}^2\left|\frac{\kappa_{\theta}v_x^2}{2v^5}\theta_t\right|\,dx \\
      & \le C N_1^2\delta^2(\delta  N_1^5\|(\phi_{xx},\phi_x,\psi_x,\zeta_x)(t)\|^2+N_1^4\delta^{\frac{3}{2}}(1+t)^{-\frac{3}{2}})\\
      &\le C N_1^7\delta^2\|(\phi_{xx},\phi_x,\psi_x,\zeta_x)(t)\|^2+C N_1^6\delta^{\frac{7}{2}}(1+t)^{-\frac{3}{2}},
  \end{split}
\end{equation}
\begin{equation}\label{3.63}
  \begin{split}
    \left|I_{10}^4\right| & \le C \delta\int_{\mathbb{R}}|(v^{-1}g^2\phi(4h-\gamma p_+\phi))_x|(|\zeta_x|+|(\phi,\zeta){\Theta}_x|)\,dx \\
      & \le C \delta\int_{\mathbb{R}}\left(\delta|(\phi,\zeta)||(\phi_x,\zeta_x)|+\delta^{\frac{1}{2}}|(\phi,\zeta)|^2|w|\right)(|\zeta_x|+|(\phi,\zeta){ \Theta}_x|)\,dx\\
      &\le C\delta^{\frac{3}{2}}\|(\phi_x,\zeta_x)(t)\|^2+CN_1^2\delta^{\frac{3}{2}}\int_{\mathbb{R}}\left|\left(\phi,\frac{\zeta}{\sqrt{\delta}}\right)\right|^2w^2\,dx,
  \end{split}
\end{equation}
\begin{equation}\label{3.64}
    \left|I_{10}^5\right|  \le C \delta^2\int_{\mathbb{R}}|(\phi,\zeta)|^2|(\psi_x,{ U}_x)|^2\,dx \le C N_1^2\delta^2\left(\|(\psi_x)(t)\|^2+\delta^{\frac{5}{2}}(1+t)^{-\frac{3}{2}}\right),
\end{equation}
\begin{eqnarray}\label{3.65}
    \left|I_{10}^6\right|&\le& \left|\delta\int_{\mathbb{R}}(v^{-1}\phi g^2(4h-\gamma p_+\phi)\frac{\theta \kappa_{\theta}v_x}{v^5})_xu_x\,dx\right| \nonumber\\
    &&+\left|\delta\int_{\mathbb{R}}(v^{-1}\phi g^2(4h-\gamma p_+\phi))_x\frac{v\kappa_{\theta v}-\kappa_{\theta}}{2v^6}\theta u_xv_x^2\,dx\right| \nonumber\\
      &\le& C\delta\int_{\mathbb{R}}\left(|(\phi,\zeta)|^2g^2v_x^2|u_x|+|\phi_xg^2||(\phi,\zeta)||v_xu_x|+|\phi gg_x(\phi,\zeta)v_x u_x|\right. \nonumber\\
      &&\left.+|\phi g^2(\phi_x,\zeta_x)v_xu_x|+|\phi g^2(\phi,\zeta)u_x||(\theta_xv_x,v_x^2,v_{xx})|+|\phi g^2(\phi,\zeta)u_xv_x^2|\right)\,dx \nonumber\\
&\le& CN_1^3\delta^{\frac{3}{2}}\|(\phi_{xx},\phi_x,\psi_x,\zeta_x)(t)\|^2+CN_1^2\delta^{\frac{3}{2}}\int_{\mathbb{R}}\left|\left(\phi,\frac{\zeta}{\sqrt{\delta}}\right)\right|^2w^2\,dx+C N_1^2\delta^{\frac{5}{2}}(1+t)^{-\frac{3}{2}},
\end{eqnarray}
where in (\ref{3.60}), we have used $|w_x|\le C\delta^{-\frac{1}{2}}(1+t)^{-1}$.

Combining (\ref{3.53})-(\ref{3.65}) and integrating the resultant equation over $[0,t]$ yields
\begin{equation}\label{3.66}
  \begin{split}
    \frac{1}{\delta}\int_0^t<h_t,hg^2>_{H^{-1}\times H^1}\,d\tau\le&  CN_1^7\delta^{\frac{1}{2}}\int_0^t\|(\phi_{xx},\phi_x,\psi_x,\zeta_x)(\tau)\|^2\,d\tau+ CN_1^6\delta\\
      & +CN_1^2\delta^{\frac{3}{2}}\int_0^t\int_{\mathbb{R}}\left|\left(\phi,\frac{\zeta}{\sqrt{\delta}}\right)\right|^2w^2\,dx.
  \end{split}
\end{equation}
Then (\ref{3.34}) follows from (\ref{3.66}) and Lemma 2.2 immediately. This completes  the proof of Lemma 3.2.

\begin{center}
{\bf Acknowledgement}
\end{center}
This work
was supported by the National Natural Science Foundation of China
(Grant No. 11501003), the  Doctoral Scientific Research Fund of
Anhui University (Grant No. J10113190005), and the Cultivation Fund
of Young Key Teacher at Anhui University.

\end{document}